%% file: ex_article.tex
\documentclass[hidelinks,onefignum,onetabnum]{siamart251216}

\usepackage{adjustbox}
\usepackage{subcaption}
\usepackage{booktabs} 
\usepackage[table]{xcolor}
\usepackage{multirow}
\input{ex_shared}

\ifpdf
\hypersetup{
  pdftitle={Equidistribution of Solution Gradients in Physics-Informed Neural Networks via Coordinate Transformations},
  pdfauthor={T. Yu and I. Oseledets}
}
\fi

\begin{document}

\maketitle

\begin{abstract}
Physics-Informed Neural Networks (PINNs) often struggle to approximate solutions with sharp gradients due to optimization difficulties and spectral bias. Motivated by the principle of gradient equidistribution from the stretching approach in grid generation, which redistributes solution gradients to reduce functional complexity, we propose Coordinate Transformation PINNs (CTPINNs), which incorporate learnable coordinate transformations to map the physical domain into a computational domain where the physical solutions become approximately affine, making them easier for neural networks to learn. We propose two constructions of the transformation: an explicit construction based on Majorants and an implicit construction based on invertable neural networks. Experimentally, the explicit construction can improve accuracy for PDEs with available explicit coordinate transformations, and the implicit construction can be combined with other advanced PINN models as an efficient plug-in module to improve performance across a range of challenging PDE benchmarks.
\end{abstract}



\input{main/introduction}
\input{main/CTPINN}
\input{main/theoretical_analysis}
\input{main/experiments}
\input{main/conclusion}


\bibliographystyle{siamplain}
\bibliography{references}

\appendix
\input{main/appendices}
\end{document}

%% file: ex_shared.tex
\usepackage{lipsum}
\usepackage{amsfonts}
\usepackage{graphicx}
\usepackage{epstopdf}
\usepackage{algorithmic}
\ifpdf
  \DeclareGraphicsExtensions{.eps,.pdf,.png,.jpg}
\else
  \DeclareGraphicsExtensions{.eps}
\fi

\newsiamremark{remark}{Remark}
\newsiamremark{hypothesis}{Hypothesis}
\crefname{hypothesis}{Hypothesis}{Hypotheses}
\newsiamthm{claim}{Claim}
\newsiamthm{assumption}{Assumption}
\crefname{assumption}{Assumption}{Assumptions}
\Crefname{assumption}{Assumption}{Assumptions}
\AddToHook{env/assumption/begin}{\crefalias{theorem}{assumption}}
\newsiamremark{fact}{Fact}
\crefname{fact}{Fact}{Facts}

\headers{CTPINNs}{T. Yu and I. Oseledets}

\title{Equidistribution of Solution Gradients in Physics-Informed Neural Networks via Coordinate Transformations\thanks{Submitted to SIAM Journal on Scientific Computing.}}

\author{Tianchi Yu\thanks{Artificial Intelligence Research Institute and Skolkovo Institute of Science and Technology, Moscow, Russia.
  (\email{tianchi.yu@skoltech.ru}).}
\and Ivan Oseledets\thanks{Artificial Intelligence Research Institute and Skolkovo Institute of Science and Technology, Moscow, Russia.
  (\email{oseledets@airi.net}).}}

\usepackage{amsopn}


%% file: main/introduction.tex
\section{Introduction}
Physics-Informed Neural Networks (PINNs)~\cite{raissi2019physics} have recently emerged as a powerful framework for solving differential equations by embedding physical laws into neural network training. PINNs have been successfully applied across a wide range of domains, including fluid dynamics~\cite{wassing2025physics}, structural mechanics~\cite{schmid2026physics}, power systems~\cite{xu2026physics}, biomedical modeling~\cite{mehmood2026physics}, and environmental science~\cite{di2026physics}. 

Most existing research on PINNs focuses on designing loss functions~\cite{jin2021nsfnets,chiu2022can,yu2025spectral}, network architectures~\cite{shukla2024comprehensive,wang2024piratenets,lai2025hard}, and training strategies~\cite{yu2022gradient,rathore2024challenges,wang2026gradient} to improve the accuracy and efficiency of solving differential equations. In addition, adaptive sampling techniques~\cite{wu2023comprehensive,gao2023failure,tang2023adversarial,lau2024pinnacle} aim to improve the accuracy with which the loss approximates the underlying integral. However, these approaches primarily address the error between the neural network approximation and the target solution from the perspective of the network and its training: they seek to improve the network's ability to represent the solution or to optimize its approximation more effectively. On the other hand, reducing the complexity of the solution itself has received insufficient attention in the PINN literature.

Mesh generation is a fundamental component of numerical methods for partial differential equations, particularly finite element methods, and related ideas have also been explored in the context of PINNs. Adaptive mesh refinement~\cite{berger1984adaptive} has been incorporated into PINNs in the form of the aforementioned adaptive sampling methods, which resample points based on residual-related criteria. On the other hand, variational grid generation methods~\cite{li2001moving} have also been explored in the PINN context. In~\cite{yang2024moving}, the authors first train a neural network to approximate the coordinate transformation (CT) by minimizing a variational energy functional and then use the learned CT to generate transformed grids for subsequent PINN training. In~\cite{hou2023enhancing}, the authors propose an adaptive sampling strategy that employs a monitor function to adjust the random distribution of sampling points. However, these approaches still focus on the network and its training.

In this paper, from the perspective of reducing solution complexity, we propose Coordinate Transformation Physics-Informed Neural Networks (CTPINNs), a novel framework that incorporates learnable coordinate transformations into PINNs. By jointly optimizing the coordinate transformation and solution network, CTPINNs induce an implicit variational objective that promotes gradient equidistribution, stretching the solution toward an approximately affine form and reducing its complexity. This transformation mitigates spectral bias by mapping high-frequency components to lower frequencies, thereby improving accuracy and accelerating optimization.

The main contributions of this paper are as follows:
\begin{itemize}
    \item We propose CTPINNs, a novel framework that integrates learnable coordinate transformations into PINNs, enabling the network to solve PDEs in a transformed domain with reduced solution complexity without explicitly changing the form of loss functions.

    \item We develop explicit and implicit constructions of the coordinate transformation. The explicit construction uses explicit functions derived from known solution features, whereas the implicit construction represents the transformation by a learnable invertible neural network. Both are optimized jointly with the solution network.

    \item We evaluate both constructions on comprehensive PDE benchmarks. The numerical results show that the explicit construction resolves sharp gradients, while the implicit construction, incorporating the implicit coordinate-transformation module into state-of-the-art PINNs, improves accuracy on challenging problems.
\end{itemize}

%% file: main/CTPINN.tex
\section{Methodology}
\subsection{Physics-Informed Neural Networks}
We begin by revisiting the foundational framework of physics-informed neural networks (PINNs) \cite{raissi2019physics, lagaris1998artificial} for solving partial differential equations (PDEs).

Consider a general time-dependent PDE system for the unknown field $\boldsymbol{u}(t, \boldsymbol{x})$:
\begin{equation}
\begin{aligned}
    & \partial_{t}\boldsymbol{u}+\mathcal{N}[\boldsymbol{u}]=0, \quad t \in[0, T],\ \boldsymbol{x} \in \Omega, \\
    & \boldsymbol{u}(0, \boldsymbol{x})=\boldsymbol{g}(\boldsymbol{x}), \quad \boldsymbol{x} \in \Omega, \\
    & \mathcal{B}[\boldsymbol{u}]=0, \quad t \in[0, T],\ \boldsymbol{x} \in \partial \Omega,
\end{aligned}\label{PDE}
\end{equation}
where $\mathcal{N}$ denotes the differential operator, $\Omega$ is the spatial domain, and $\mathcal{B}$ enforces boundary conditions. For time-independent problems, the temporal derivative vanishes ($\partial_{t}\boldsymbol{u} \equiv 0$).

PINNs utilize a neural network $\boldsymbol{u}^{\theta}$ (parameterized by trainable weights and biases) to approximate the unknown solution $\boldsymbol{u}$. This is achieved by minimizing a composite loss that enforces the PDE system \eqref{PDE} at collocation points:
\begin{equation}\label{eq: PINN}
\mathcal{L}(\theta)=\lambda_{ic} \mathcal{L}_{ic}(\theta)+\lambda_{bc} \mathcal{L}_{bc}(\theta)+\lambda_r \mathcal{L}_r(\theta),
\end{equation}
with individual terms defined as
\begin{equation}
\begin{aligned}
& \mathcal{L}_r(\theta)=\frac{1}{N_r} \sum_{i=1}^{N_r}\left|\partial_{t}\boldsymbol{u}^\theta\left(t_r^i, \boldsymbol{x}_r^i\right)+\mathcal{N}\left[\boldsymbol{u}^\theta\right]\left(t_r^i, \boldsymbol{x}_r^i\right)\right|^2,\\
& \mathcal{L}_{ic}(\theta)=\frac{1}{N_{ic}} \sum_{i=1}^{N_{ic}}\left|\boldsymbol{u}^\theta\left(0, \boldsymbol{x}_{ic}^i\right)-\boldsymbol{g}\left(\boldsymbol{x}_{ic}^i\right)\right|^2, \\
& \mathcal{L}_{bc}(\theta)=\frac{1}{N_{bc}} \sum_{i=1}^{N_{bc}}\left|\mathcal{B}\left[\boldsymbol{u}^\theta\right]\left(t_{bc}^i, \boldsymbol{x}_{bc}^i\right)\right|^2, \\
\end{aligned}
\end{equation}
where $\lambda_{ic}, \lambda_{bc}, \lambda_r$ are positive scalar weights balancing the relative importance of the initial condition, boundary condition, and PDE residual terms; $N_{ic}, N_{bc}, N_r$ denote the number of points sampled for the initial condition, boundary condition, and PDE residual, respectively, with corresponding collocation points $\boldsymbol{x}_{ic}^i$, $(t_{bc}^i, \boldsymbol{x}_{bc}^i)$, and $(t_r^i, \boldsymbol{x}_r^i)$.

\subsection{Coordinate Transformations and Majorants}\label{sec: coordinate_transformations_and_majorants}
Using coordinate transformations (CTs) to solve PDEs is one of the grid-generation approaches~\cite{liseikin2017grid}. The main idea is to construct a function $\boldsymbol{T}(\boldsymbol{x})=\boldsymbol{\xi}: \Omega \to \hat{\Omega}$, where $\Omega$ is the physical domain and $\hat{\Omega}$ is the computational domain. Classically, $\boldsymbol{T}$ concentrates resolution in regions with rapid solution variation (e.g., boundary/interior layers and steep gradients), improving accuracy without grid refinement by equidistributing the local truncation error across the domain and reducing the condition number of the resulting linear systems.

Mathematically, the CT function $\boldsymbol{T}$ is a mapping that satisfies the following properties:
\begin{itemize}
\item \textit{Bijectivity.} The CT function should be a bijection between the physical domain $\Omega$ and the computational domain $\hat{\Omega}$, ensuring that each point in $\Omega$ corresponds to exactly one point in $\hat{\Omega}$ and vice versa.
    \item \textit{Smoothness.} The CT function should be at least $C^1$ continuous.
    \item  \textit{Boundary-fittedness.} The CT function should map the boundaries of the physical domain to the boundaries of the computational domain. 
\end{itemize}
In addition to the above properties, CT functions used in numerical methods should satisfy other conditions, including \textit{non-degeneracy}, \textit{gradation} (quasi-uniformity), and \textit{topological compatibility}. These properties are required for the stability of numerical methods. Since PINNs are mesh-free, these additional properties are not necessarily required of CT functions in PINNs. Without loss of generality, in the following discussion, we assume that $\boldsymbol{T}: [0,1]^d\to[0,1]^d$ satisfies the above properties, since the domain can be normalized and transformed to the unit cube $[0,1]^d$.
\subsubsection{Single-variable Majorant Functions}
One direct application of this approach, known as the stretching approach in grid generation, is the use of majorant functions~\cite{mali2013accuracy} to quantify the singular behavior of solutions to singularly perturbed PDEs, which typically contain a small parameter $\epsilon$ multiplying the highest-order derivatives. These majorants provide $\epsilon$-uniform upper bounds on the solution gradients, ensuring that grid clustering is applied only in proportion to the local solution variation.

A function $\psi(x)$, $x\in[0,1]$, is a majorant of $u(x)$ if it is integrable over $[0,1]$ and satisfies $|u'(x)|\leq\psi(x)$ for all $x\in[0,1]$. We then define the coordinate transformation $\xi=T(x)$ by
\begin{equation}\label{eq: coordinate_transformation_def}
    \begin{aligned}
    & \frac{\mathrm{d}T}{\mathrm{d}x}= c\psi(x), \\
    & T(0)=0.        
    \end{aligned}
\end{equation}
Choosing $c$ so that $T(1)=1$ gives
\begin{equation}
    \left|\frac{\mathrm{d}u}{\mathrm{d}\xi}\right| = \left|\frac{\mathrm{d}u}{\mathrm{d}x}\cdot\frac{\mathrm{d}x}{\mathrm{d}\xi}\right| \leq \frac{1}{c}.
\end{equation}
The following \cref{thm:equidistribution_optimality} shows that a unique gradient-equidistributing coordinate transformation exists and minimizes the squared $L^2$ norm of the solution gradient in the computational coordinate.
\begin{theorem}\label{thm:equidistribution_optimality}
Let $u \in W^{1,1}([0,1])$ with $u'(x)\neq 0$ a.e. on $[0,1]$, and define the admissible set
\begin{equation}
    \mathcal{A}:=\left\{T\in W^{1,1}([0,1]): T(0)=0,\ T(1)=1,\ T'(x)>0\ \text{a.e. on }[0,1]\right\}.
\end{equation}
For $T \in \mathcal{A}$, let $v(\xi) := u(T^{-1}(\xi))$. Then
\begin{equation}
    \int_0^1 |v'(\xi)|^2\,\mathrm{d}\xi
    =\int_0^1 \frac{|u'(x)|^2}{T'(x)}\,\mathrm{d}x.
\end{equation}
Moreover, the unique solution to the optimization problem
\begin{equation}
    \min_{T \in \mathcal{A}} \int_0^1 \frac{|u'(x)|^2}{T'(x)}\,\mathrm{d}x
\end{equation}
is
\begin{equation}
    (T')^*(x)=\frac{|u'(x)|}{\int_0^1 |u'(s)|\,\mathrm{d}s}.
\end{equation}
\end{theorem}
Although \cref{thm:equidistribution_optimality} is a classical result that has been widely discussed~\cite{de2006good,dorfi1987simple} and proved for certain PDE classes~\cite{huang1994moving}, we include it for completeness and provide its proof in \cref{appendix: proof_equidistribution_optimality}.

Thus, a function $u(x)$ whose total variation on $[0,1]$ is bounded by a constant independent of $\epsilon$ is transformed into $u(T^{-1}(\xi))$, whose first derivative with respect to $\xi$ is uniformly bounded on $[0,1]$. In this paper, we introduce and use three types of majorant functions, $\psi_1$, $\psi_2$, and $\psi_3$, and their corresponding coordinate transformations, $T_{\psi_1}$, $T_{\psi_2}$, and $T_{\psi_3}$. Define
\begin{equation}
	\frac{\mathrm{d}T_{\psi_i}}{\mathrm{d}x} = \psi_i(x), \quad x \in [0,1], i=1,2,3.
\end{equation}
Thus, integrating these expressions gives explicit formulas for $T_{\psi_i}(x)$. We then have:

\paragraph{Majorant1}
\begin{equation}\label{eq: majorant_1}
\psi_1(x; x_0,b,\alpha)=\alpha^{-1}\exp\left(-\frac{b|x-x_0|}{\alpha}\right).
\end{equation}

\begin{equation}\label{eq: int_majorant_1}
T_{\psi_1}(x; x_0,b,\alpha)=
    \begin{cases}
        \dfrac{1}{b}\left(1-\exp\left(-\dfrac{b(x-x_0)}{\alpha}\right)\right)+C, & x\ge x_0,\\
        -\dfrac{1}{b}\left(1-\exp\left(-\dfrac{b(x_0-x)}{\alpha}\right)\right)+C, & x<x_0.
	\end{cases}
\end{equation}

\paragraph{Majorant2}
\begin{equation}\label{eq: majorant_2}
\psi_2(x; x_0,b,\alpha)=\frac{\alpha^{b}}{\left(\alpha+|x-x_0|\right)^{b+1}}.
\end{equation}

\begin{equation}\label{eq: int_majorant_2}
T_{\psi_2}(x; x_0,b,\alpha)=
    \begin{cases}
        \dfrac{1}{b}\left(1-\left(\dfrac{\alpha}{\alpha+x-x_0}\right)^b\right)+C, & x\ge x_0,\\
        -\dfrac{1}{b}\left(1-\left(\dfrac{\alpha}{\alpha+x_0-x}\right)^b\right)+C, & x<x_0.
	\end{cases}
\end{equation}

\paragraph{Majorant3}
\begin{equation}\label{eq: majorant_3}
\psi_3(x; x_0,\alpha)=\frac{1}{\left(\alpha+|x-x_0|\right)|\log\alpha|}.
\end{equation}

\begin{equation}\label{eq: int_majorant_3}
T_{\psi_3}(x; x_0,\alpha)=
    \begin{cases}
        \dfrac{1}{|\log\alpha|}\log\left(1+\dfrac{x-x_0}{\alpha}\right)+C, & x\ge x_0,\\
        -\dfrac{1}{|\log\alpha|}\log\left(1+\dfrac{x_0-x}{\alpha}\right)+C, & x<x_0.
	\end{cases}
\end{equation}

Here, $\alpha>0$, $b>0$, and $x_0\in[0,1]$ are hyperparameters or trainable parameters. The integration constant $C$ is determined by normalizing the transformation as $T_{\psi_i}(x)\to\frac{T_{\psi_i}(x)-T_{\psi_i}(0)}{T_{\psi_i}(1)-T_{\psi_i}(0)}$, which ensures that $T_{\psi_i}:[0,1]\to[0,1]$, $T_{\psi_i}(0)=0$, and $T_{\psi_i}(1)=1$.

We also summarize the properties and application range of these three majorants in \cref{tab:majorants}, which can serve as a guide for selecting an appropriate majorant for a given problem in CTPINNs.

\begin{table}[ht]
\centering
\caption{Properties and application ranges of three majorants.}
\label{tab:majorants}
\begin{adjustbox}{width=\textwidth}
\begin{tabular}{lcc}
\toprule
\textbf{Majorant} & \textbf{Suitable PDE Type} & \textbf{Physical Scenario Keywords} \\
\midrule
$\psi_1$ & Parabolic PDEs, Navier-Stokes Equations & High Reynolds number, rapidly decaying layer \\
$\psi_2$ & Thin shell equations, Seepage equations & Adsorbed flow, thin structure deformation \\
$\psi_3$ & Potential field equations, Orifice flow equations & Corner singularities, slow diffusion \\
\bottomrule
\end{tabular}
\end{adjustbox}
\end{table}

\subsubsection{Multi-variable Majorant Functions}

Multivariable majorant functions are constructed from the single-variable majorants $\psi$ either by applying majorants independently along each dimension or by using transfinite interpolation along the edges of the domain.

\begin{figure}[ht]
\centering
\subcaptionbox{Uniform \label{fig:grid_uniform}}[0.23\textwidth]{%
\includegraphics[width=\linewidth]{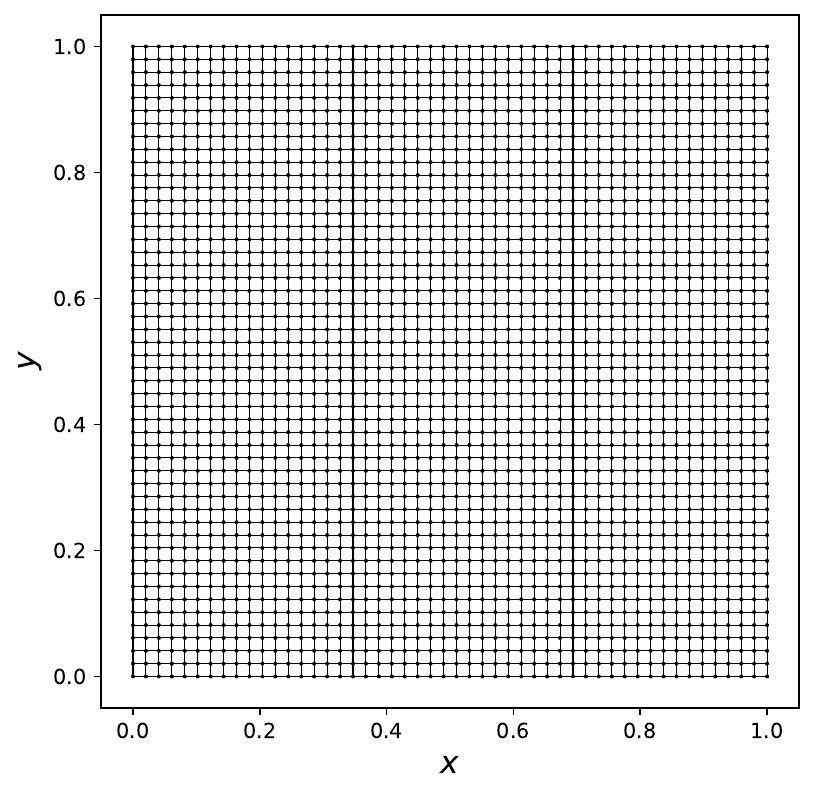}}
\subcaptionbox{$(0.0,0.0)$\label{fig:grid_0_0}}[0.23\textwidth]{%
\includegraphics[width=\linewidth]{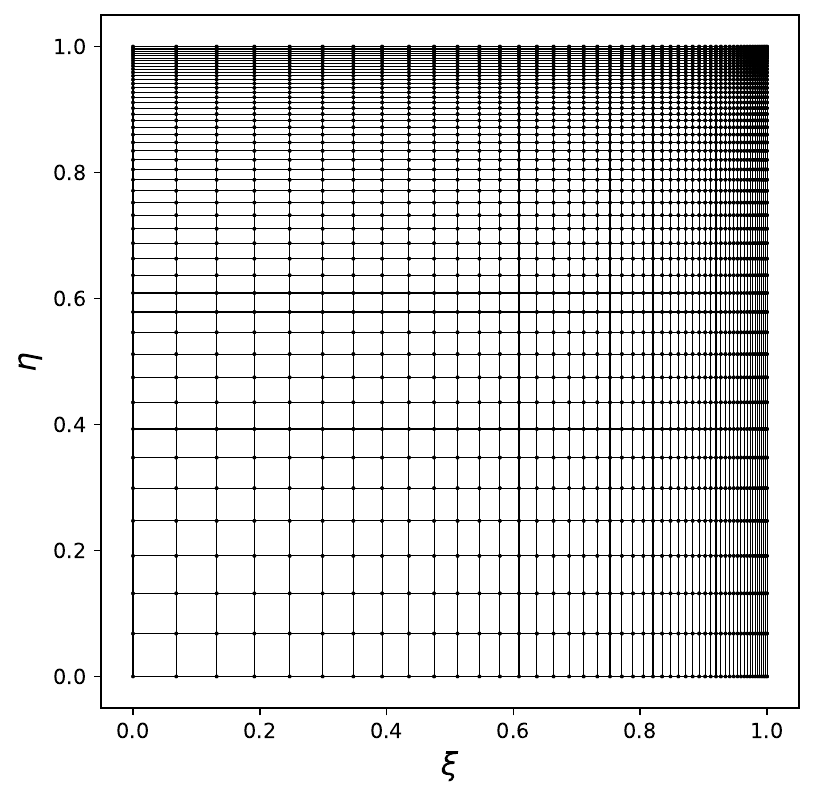}}
\subcaptionbox{$(0.5,0.0)$\label{fig:grid_0_5_0}}[0.23\textwidth]{%
\includegraphics[width=\linewidth]{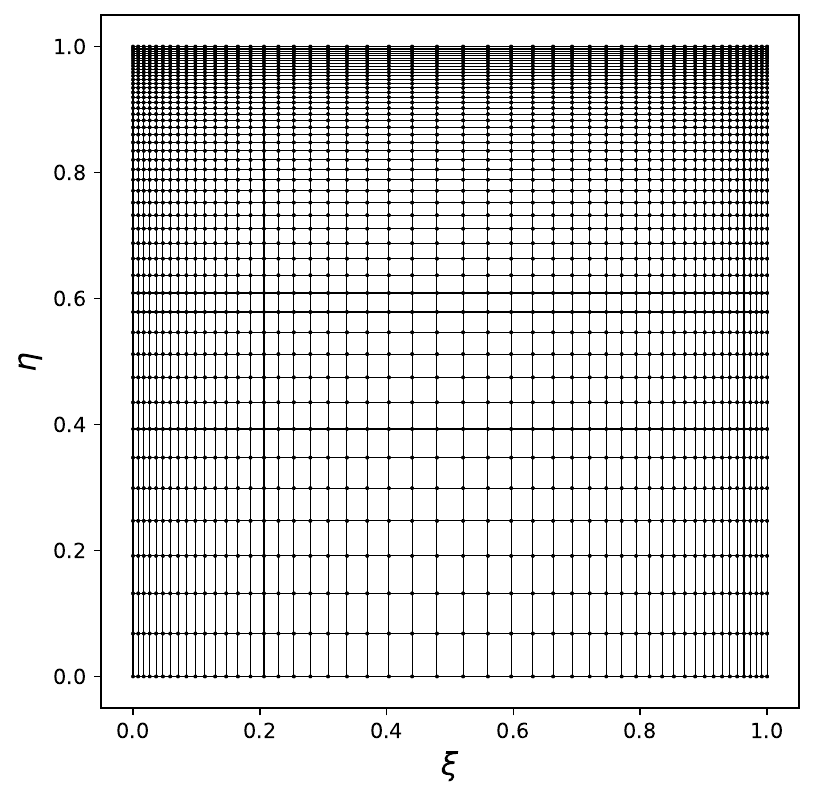}}
\subcaptionbox{$(0.5,0.5)$\label{fig:grid_0_5_0_5}}[0.23\textwidth]{%
\includegraphics[width=\linewidth]{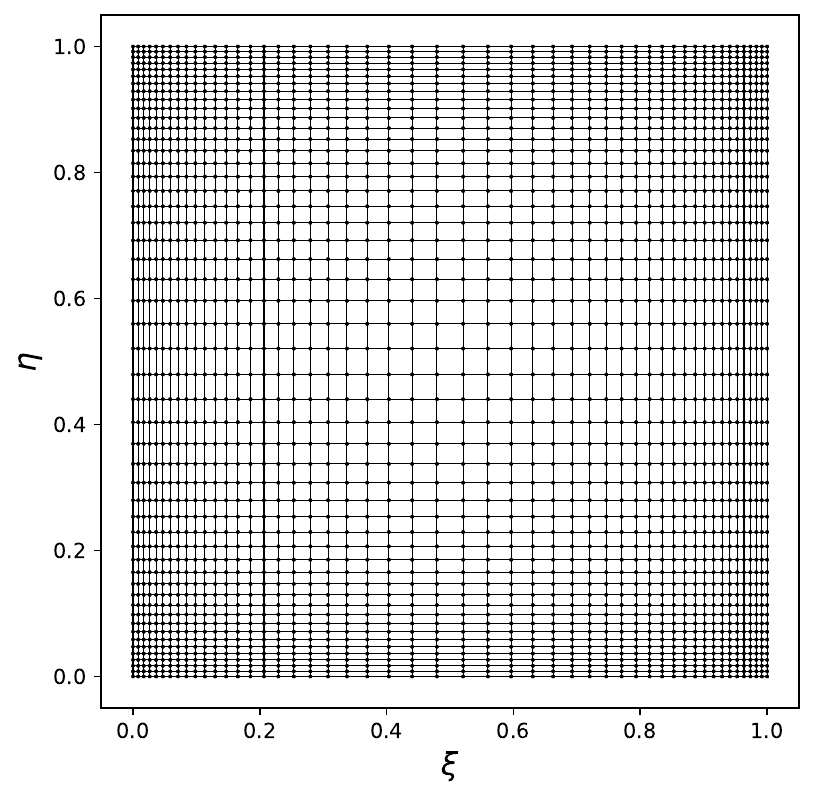}}\\
\subcaptionbox{Transfinite \label{fig:grid_transfinite_0_5_0_2_1_0_5}}[0.23\textwidth]{%
\includegraphics[width=\linewidth]{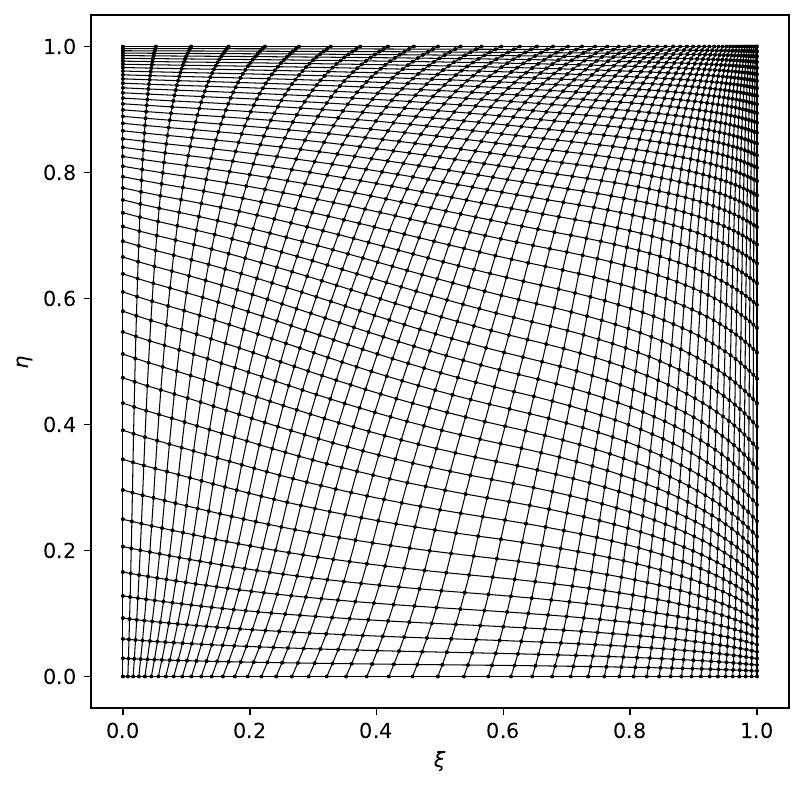}}
\subcaptionbox{Transfinite\label{fig:grid_transfinite_0_6_0_8_0_1_0_2}}[0.23\textwidth]{%
\includegraphics[width=\linewidth]{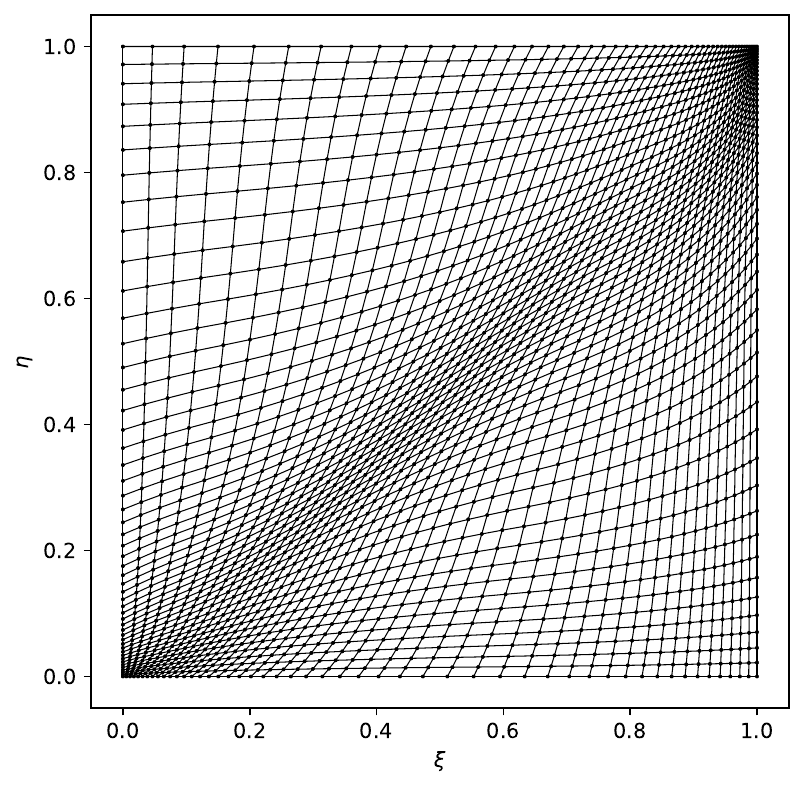}}
\subcaptionbox{Transfinite\label{fig:grid_transfinite_0_6_0_0_0_7_0_5}}[0.23\textwidth]{%
\includegraphics[width=\linewidth]{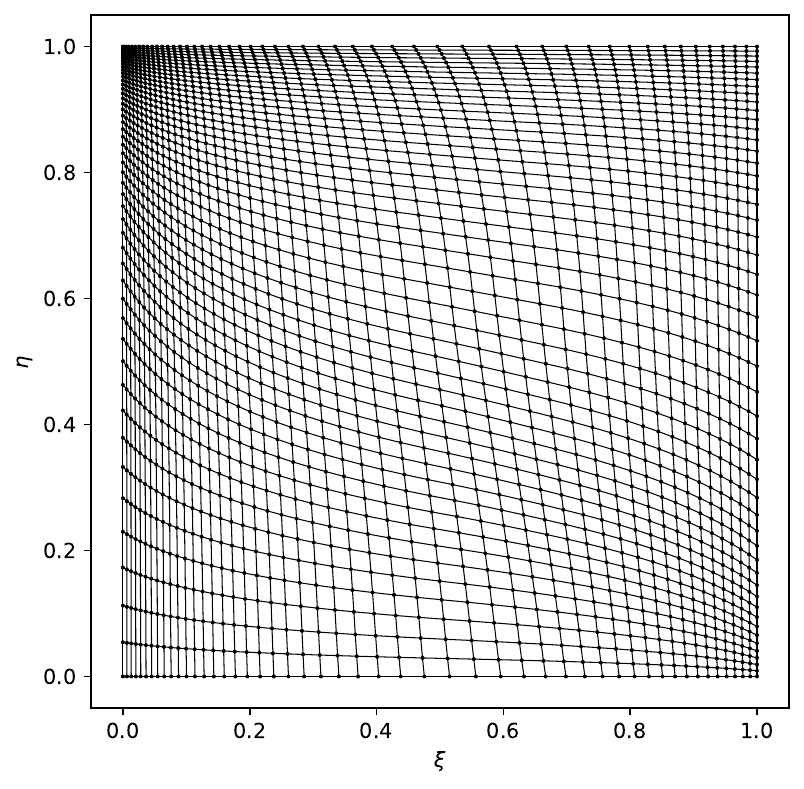}}
\subcaptionbox{Transfinite\label{fig:grid_transfinite_0_9_0_4_0_4_0_7}}[0.23\textwidth]{%
\includegraphics[width=\linewidth]{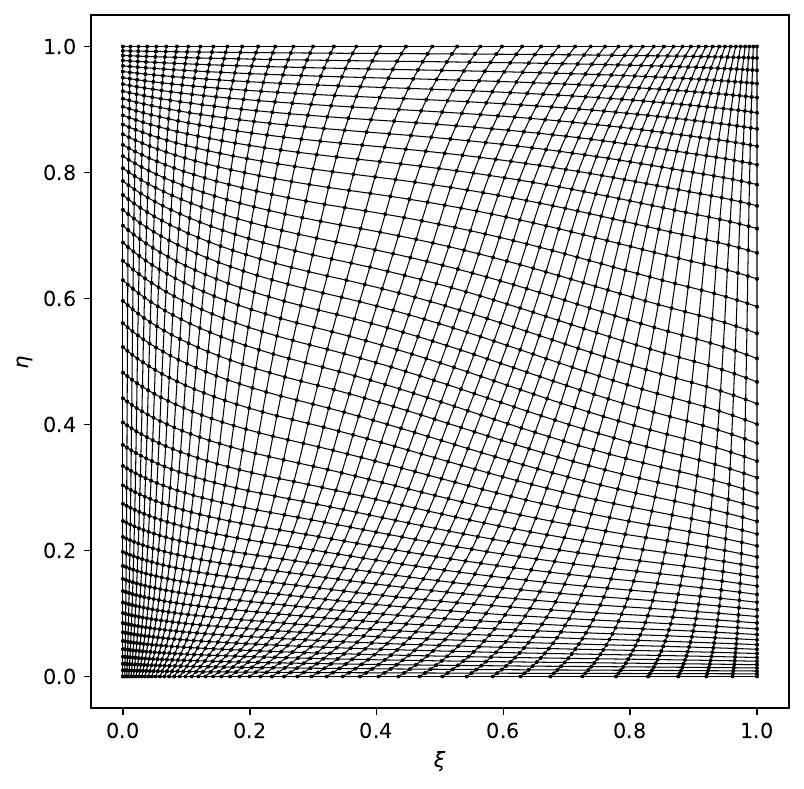}}
\caption{Examples of coordinate transformations for 2D problems. (a): Uniform grid. (b)--(d): Coordinate transformations constructed by applying Majorant1 with $\epsilon=0.3$ along each dimension. The labels indicate $(x_0,y_0)$. (e)--(h): Coordinate transformations constructed by transfinite interpolation with randomly generated parameters.}
\label{fig: grids}
\end{figure}

\paragraph{Applying majorants on each dimension} For a $d$-dimensional function, we can construct the coordinate transformation by applying a majorant independently along each dimension. Specifically, we define $\boldsymbol{T}(\boldsymbol{x})=(T_1(x_1),\dots,T_d(x_d))$, where each $T_i$ is a majorant. Examples of transformations constructed by applying Majorant1 with $\epsilon=0.3$ along each dimension are shown in \cref{fig:grid_0_0,fig:grid_0_5_0,fig:grid_0_5_0_5}.
\paragraph{Transfinite interpolation} The transfinite interpolation method~\cite{gordon1973construction} constructs multidimensional coordinate transformations by blending majorants along the domain edges. For example, when $d=2$, we can define the coordinate transformation as follows:
\begin{equation}
	\begin{aligned}
		T_1(x,y) &= (1-y)T_1^{(0)}(x) + y T_1^{(1)}(x), \\
		T_2(x,y) &= (1-x)T_2^{(0)}(y) + x T_2^{(1)}(y),
	\end{aligned}
\end{equation}
where $T_1^{(0)}$ and $T_1^{(1)}$ are majorants applied along the edges $y=0$ and $y=1$, respectively, and $T_2^{(0)}$ and $T_2^{(1)}$ are majorants applied along the edges $x=0$ and $x=1$, respectively. This method ensures that the coordinate transformation transitions smoothly between the different single-variable transformations along the domain edges and has been shown to be effective in practice~\cite{fanaskov2023general}. Examples of transformations constructed by transfinite interpolation with randomly generated parameters are shown in \cref{fig:grid_transfinite_0_5_0_2_1_0_5,fig:grid_transfinite_0_6_0_8_0_1_0_2,fig:grid_transfinite_0_6_0_0_0_7_0_5,fig:grid_transfinite_0_9_0_4_0_4_0_7}. Implicit constructions are naturally more flexible and can be applied to more complex geometries.

\subsection{CTPINNs: Coordinate Transformation Physics-Informed Neural Networks}
We propose Coordinate Transformation Physics-Informed Neural Networks (CTPINNs). Mathematically, let the network output be $u^\theta(\boldsymbol{x})=\mathcal{NN}(\boldsymbol{x};\theta)$, where $\boldsymbol{x}$ is the network input. After incorporating the CT function $\boldsymbol{T}$, the network output in CTPINNs becomes $u^\theta(\boldsymbol{x})=\mathcal{NN}(\boldsymbol{T}(\boldsymbol{x});\theta)$. The function $\boldsymbol{T}(\boldsymbol{x})$ is designed to satisfy \cref{eq: coordinate_transformation_def}, thereby making the transformed solution in the computational domain simpler and easier for the neural network to learn.

\begin{figure}[ht]
    \centering
    \includegraphics[width=0.8\textwidth]{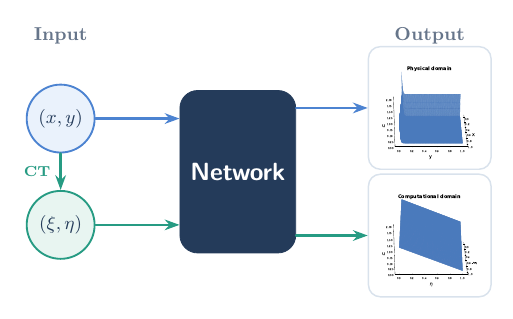}
    \caption{Comparison of PINNs (\textcolor[HTML]{4C83D1}{blue line}) and CTPINNs (\textcolor[HTML]{269B83}{green line}) on the 2D boundary-layer problems in \cref{sec: 2D_boundary_layer_problems}. The CT module maps the physical domain $\boldsymbol{x}=(x,y)$ to the computational domain $\boldsymbol{\xi}=(\xi,\eta)$. In PINNs, the network takes $(x,y)$ as input and approximates the solution in the physical domain. In CTPINNs, the network takes $(\xi,\eta)$ as input and approximates the solution in the computational domain. The transformed solution has lower functional complexity (e.g., smaller gradients and fewer high-frequency components) than the solution in the physical domain and is therefore easier for the neural network to learn.}
    \label{fig:framework}
\end{figure}

\cref{fig:framework} illustrates the benefits of CTPINNs. Learning the solution in the physical domain can be difficult because of sharp gradients and the spectral bias of neural networks~\cite{rahaman2019spectral}. The CT module maps the physical domain $(x,y)\in\Omega$ to the computational domain $(\xi,\eta)\in\hat{\Omega}$, where the solution is simpler and therefore easier for the network to learn.

We next introduce two approaches for constructing $\boldsymbol{T}$. The explicit approach (\cref{section: explicit_construction}) directly employs the explicit majorants, whereas the implicit approach (\cref{section: implicit_construction}) is designed for more general problems in which explicit majorants are unavailable.
\subsubsection{Explicit Construction}\label{section: explicit_construction}
The explicit construction directly defines the CT function $\boldsymbol{T}$ using explicit majorants. For example, consider the following one-dimensional boundary-layer problem, whose analytical solution is $u^*=e^{-x/\epsilon}$:
\begin{equation}
\label{eq: bl1d}
\begin{aligned}
     \epsilon u_{xx} + u_x &= 0, \quad x \in [0,1],\\
     u(0) &= 1, \quad u(1) = e^{-1/\epsilon}.
\end{aligned}
\end{equation}
Then for PINNs, the loss function is defined as:
\begin{equation}
    \mathcal{L}(\theta) = \frac{1}{N_r}\sum_{i=1}^{N_r} |\epsilon u^\theta_{xx}(x_i) + u^\theta_x(x_i)|^2 + \frac{1}{2} |u^\theta(0)-1|^2 + \frac{1}{2} |u^\theta(1)-e^{-1/\epsilon}|^2.
\end{equation}
Consider Majorant1 in \cref{eq: majorant_1} with $\alpha=\epsilon$, $b=1$, and $x_0=0$, which is the optimal choice for this problem (by \cref{thm:equidistribution_optimality}). The CT function for this boundary-layer problem is $T_{\mathrm{bl}}(x)=-\exp(-x/\alpha)$. After normalization, we have $\xi=T(x)=\frac{T_{\mathrm{bl}}(x)-T_{\mathrm{bl}}(0)}{T_{\mathrm{bl}}(1)-T_{\mathrm{bl}}(0)}$, and the PDE in \cref{eq: bl1d} can be transformed into:
\begin{equation}
\label{eq: bl1d_transformed}
\begin{aligned}
     & \epsilon\left(\frac{dT}{dx}u_\xi\right)_\xi\frac{dT}{dx}+u_\xi\frac{dT}{dx}=0, \quad \xi \in [0,1],\\
    \Rightarrow & \frac{\epsilon}{\alpha(1-e^{-1/\alpha})}\Big[1-\xi(1-e^{-1/\alpha})\Big] u_{\xi\xi}
+\left(1-\frac{\epsilon}{\alpha}\right)u_\xi =0,\quad \xi\in[0,1], \\
     & u(0) = 1, \quad u(1) = e^{-1/\epsilon}.
\end{aligned}
\end{equation}
If $\alpha=\epsilon$, the residual equation can be reduced to $u_{\xi\xi}=0$, indicating that the solution is linear in the transformed coordinate $\xi$. However, the transformed residual equation in \cref{eq: bl1d_transformed} is more complex and introduces stronger coupling between $\alpha$ and the higher-order derivatives, which can make optimizing $\alpha$ less stable. This issue becomes even more pronounced for more complex PDEs, such as nonlinear ones. Therefore, for more general PDEs, we retain the original residual equation and leave the PINN loss formulation unchanged in CTPINNs. This is possible because automatic differentiation applies the chain rule to $u^\theta(T(x))$, e.g., $u_x^\theta(T(x))=u^\theta_\xi(\xi)\frac{dT}{dx}$, and similarly for higher-order derivatives. Since automatic differentiation handles the transformation implicitly, CTPINNs can be trained with the same loss function, without explicit modification, as follows:
\begin{equation}
    \mathcal{L}(\theta) = \frac{1}{N_r}\sum_{i=1}^{N_r} |\epsilon u^\theta_{xx}(T(x_i)) + u^\theta_x(T(x_i))|^2 + \frac{1}{2} |u^\theta(T(0))-1|^2 + \frac{1}{2} |u^\theta(T(1))-e^{-1/\epsilon}|^2.
\end{equation}

The explicit construction is simple and effective, but it requires prior knowledge of the solution behavior derived from the formulas of given PDEs to construct the majorant function, even when adaptivity is allowed. Such prior knowledge may be unavailable when the solution exhibits sharp gradients in several distinct regions. Therefore, we also propose an implicit construction method that does not require prior knowledge of the solution behavior and can be applied to more general problems.
\subsubsection{Implicit Construction}\label{section: implicit_construction}
The explicit construction can admit countably many identifiable points at which $u^\prime=0$: the domain can be partitioned at these points, and a majorant constructed on each subinterval. However, when the zero set is uncountable or its structure is too complex or unknown, a learnable explicit construction may not be available. This motivates the implicit construction, which represents $\boldsymbol{T}$ with a learnable neural network $\boldsymbol{T}^\theta$ and trains it jointly with the solution network. Since automatic differentiation also applies to this construction, it can be trained without changing the loss function.

Neural networks have already been investigated as a tool for mesh generation~\cite{lei2023s}. In~\cite{yu2024flow2mesh, zhang2024towards}, the authors propose a neural network that learns mesh-generation mappings from high-accuracy mapping pairs using supervised learning. In~\cite{peng20253dmeshnet}, the authors propose using PINNs to learn both harmonic PDEs (to generate smooth CT functions) and equations describing mesh requirements. In~\cite{rowbottom2025g,wang2025ugm2n}, the authors propose specific loss functions for learning mesh-generation mappings. However, these methods are designed for classical numerical methods, which are not mesh-free and have the additional stability requirements discussed above. They are not designed to train the coordinate transformation jointly with the solution model in PINNs, and their networks are not necessarily differentiable. Accordingly, to construct $\boldsymbol{T}^\theta$ for CTPINNs, we first introduce invertible neural networks (INNs), which are the cornerstone of the implicit construction.

\paragraph{Invertible neural networks} Invertible neural networks (INNs) are a class of neural networks designed to learn bijective mappings, such that the input can be uniquely recovered from the output. Representative approaches include invertible residual networks~\cite{behrmann2019invertible}, normalizing flows~\cite{rezende2015variational}, and continuous normalizing flows~\cite{chen2018neural}. In CTPINNs, we employ invertible networks based on the following formulation for the implicit construction:

\begin{equation}\label{eq: implicit_exp}
    \boldsymbol{T}_{\mathrm{exp}}^\theta(\boldsymbol{x}) = \boldsymbol{x}\odot\exp\left(p\,\boldsymbol{f}_{\mathrm{ct}}^\theta(\boldsymbol{x})\right), \quad \boldsymbol{x}\in [0,1]^d
\end{equation}
where $\boldsymbol{f}_{\mathrm{ct}}^\theta$ is a neural network, $p>0$ is a hyperparameter, and $\odot$ denotes element-wise multiplication. A larger value of $p$ can increase the gradient magnitude of $\boldsymbol{T}_{\mathrm{exp}}^\theta$, allowing the transformation to better handle steep gradients. Notably, $\boldsymbol{T}_{\mathrm{exp}}^\theta$ does not map $[0,1]^d$ to itself; instead, it maps $[0,1]^d$ to $\boldsymbol{\gamma}[0,1]^d$, where $\boldsymbol{\gamma}$ is a scaling vector. However, this issue is negligible in CTPINNs because the scaling factor can be regarded as being absorbed into the weights of the downstream network.

Recall the required properties of CT functions $\boldsymbol{T}$ in \cref{sec: coordinate_transformations_and_majorants}. The \textit{bijectivity} of $\boldsymbol{T}_{\mathrm{exp}}^\theta$ can be ensured by constraining the Lipschitz constant of $\boldsymbol{f}_{\mathrm{ct}}$. The \textit{smoothness} of $\boldsymbol{T}_{\mathrm{exp}}^\theta$ can be ensured by using smooth activation functions. Its \textit{boundary-fittedness} follows from its multiplicative form and injectivity: because each exponential factor is strictly positive, $T_{\mathrm{exp},i}^\theta(\boldsymbol{x})=0$ if and only if $x_i=0$, so each coordinate face $\{x_i=0\}$ is mapped to the corresponding boundary portion $\{\xi_i=0\}$ of the image domain. Moreover, once the Lipschitz constraint below ensures global injectivity, the smooth map is a homeomorphism from $[0,1]^d$ onto its image $\omega:=\boldsymbol{T}_{\mathrm{exp}}^\theta([0,1]^d)$ and maps the boundary of the cube to the boundary of $\omega$.

Next, we explain how to ensure \textit{bijectivity} in INNs by constraining the Lipschitz constant.

The key idea is to control the Lipschitz constant to ensure invertibility. In machine learning, Lipschitz-constant control is commonly used in generative adversarial networks (GANs)~\cite{goodfellow2020generative}, especially in Wasserstein GANs~\cite{arjovsky2017wasserstein} for the Wasserstein-1 distance. Several methods can control the Lipschitz constant of neural networks, including spectral normalization~\cite{miyato2018spectral}, gradient penalty~\cite{gulrajani2017improved}, and weight clipping~\cite{arjovsky2017wasserstein}. In CTPINNs, we use spectral normalization to control the Lipschitz constant of $\boldsymbol{f}_{\mathrm{ct}}^\theta$.

Spectral normalization normalizes layer weight matrices by their largest singular values to constrain the network's Lipschitz constant. Applied to $\boldsymbol{f}_{\mathrm{ct}}^\theta$, it bounds the overall Lipschitz constant. Let $J_{\boldsymbol{f}_{\mathrm{ct}}^\theta}(\boldsymbol{x})$ denote the Jacobian of $\boldsymbol{f}_{\mathrm{ct}}^\theta$. If the weight matrix of each layer is $W_\ell$ and the activation has global Lipschitz constant $C_\sigma$, then a global upper bound on the Lipschitz constant of $\boldsymbol{f}_{\mathrm{ct}}^\theta$ is
\begin{equation}
L_{f_{\mathrm{ct}}}=\prod_{\ell}\|W_\ell\|_2\,C_\sigma\ge \sup_{\boldsymbol{x}\in[0,1]^d}\|J_{\boldsymbol{f}_{\mathrm{ct}}^\theta}(\boldsymbol{x})\|_2.
\end{equation}

Therefore, by choosing $L$ to satisfy the derived lower bound and replacing $\boldsymbol{f}_{\mathrm{ct}}^\theta$ with $\boldsymbol{f}_{\mathrm{ct}}^{\theta}/L$, we can ensure that the corresponding transformation $\boldsymbol{T}^\theta$ is invertible on its domain and, consequently, bijective onto its image. We next derive a suitable lower bound for $L$ for $\boldsymbol{T}_{\mathrm{exp}}^{\theta}$. Its Jacobian can be written as
\begin{equation}
J_{\boldsymbol{T}_{\mathrm{exp}}^{\theta}}
=\operatorname{Diag}\left(\exp\left(\frac{p\,\boldsymbol{f}_{\mathrm{ct}}^{\theta}}{L}\right)\right)
\left(I+\operatorname{Diag}(\boldsymbol{x})\frac{p}{L}J_{\boldsymbol{f}_{\mathrm{ct}}^{\theta}}\right).
\end{equation}
On $[0,1]^d$, $\|\operatorname{Diag}(\boldsymbol{x})\|_2\le 1$. Therefore, for a constant $\eta_{\mathrm{exp}}<1$, a sufficient condition for \textit{bijectivity} of the map from $[0,1]^d$ to $\gamma[0,1]^d$ is
\begin{equation}
\left\|\operatorname{Diag}(\boldsymbol{x})\frac{p}{L}J_{\boldsymbol{f}_{\mathrm{ct}}^{\theta}}\right\|_2
\le\frac{pL_{f_{\mathrm{ct}}}}{L}\le\eta_{\mathrm{exp}}<1,
\end{equation}
which is ensured by
\begin{equation}
L\ge\frac{pL_{f_{\mathrm{ct}}}}{\eta_{\mathrm{exp}}}.
\end{equation}
In implementation, we set $\eta_{\mathrm{exp}}=0.98$ and use
\begin{equation}\label{eq:update_L_exp}
L=\max\left\{1,\frac{pL_{f_{\mathrm{ct}}}}{\eta_{\mathrm{exp}}}\right\}.
\end{equation}

Accordingly, we update $L$ every $N_{\mathrm{upd}}$ iterations using \cref{eq:update_L_exp} to maintain the required \textit{bijectivity} throughout training.

%% file: main/theoretical_analysis.tex
\section{Theoretical Analysis}
Here, we analyze convergence rates using the NTK
(\cref{sec: ntk_convergence_analysis}), investigate the implicit tendency toward
gradient equidistribution (\cref{sec: implicit_tendency_t_gradient_equid}), and
present the error analysis (\cref{sec: error_analysis}). 

The following analysis focuses on training an MLP with $\tanh$ activation to solve the one-dimensional boundary-layer problem in \cref{eq: bl1d} using the CTPINN framework. Let
\begin{equation}
    u^*(x)=e^{-x/\epsilon}
\end{equation}
denote its analytical solution, and let $v^\theta:[0,1]\to\mathbb{R}$ denote the network after the CT function.  The CT function used is
\begin{equation}\label{eq: normalized_exponential_ct}
    T_\alpha(x):=\frac{1-\exp(-x/\alpha)}{1-\exp(-1/\alpha)}, 
\end{equation}
where $\alpha>0$ is a fixed or learnable parameter. We define the CTPINN output as
\begin{equation}\label{eq: shared_boundary_layer_prediction}
    u^{\theta,\alpha}(x):=v^\theta(T_\alpha(x)).
\end{equation}
Given a set of distinct interior collocation points $\mathcal{X}=\{x_1,\ldots,x_{N_r}\}\subset(0,1)$, the
loss\footnote{Compared with the convention in \cref{eq: PINN}, we include an additional factor of $1/2$ in each term to align with the standard NTK gradient-flow formulation.} is
\begin{equation}\label{eq: unit_weight_boundary_layer_loss}
\begin{aligned}
    \mathcal{L}(\theta,\alpha)
    &: =\mathcal{L}_r(\theta,\alpha)
    +\mathcal{L}_{bc}(\theta,\alpha),\\
    \mathcal{L}_r(\theta,\alpha)
    &: =\frac{1}{2N_r}\sum_{i=1}^{N_r}
    \left|\epsilon u^{\theta,\alpha}_{xx}(x_i)
    +u^{\theta,\alpha}_x(x_i)\right|^2,\\
    \mathcal{L}_{bc}(\theta,\alpha)
    &: =\frac{1}{2}\left|u^{\theta,\alpha}(0)-1\right|^2
    +\frac{1}{2}\left|u^{\theta,\alpha}(1)-e^{-1/\epsilon}\right|^2.
\end{aligned}
\end{equation}

\subsection{NTK Convergence Analysis}\label{sec: ntk_convergence_analysis}
We use NTK analysis to compare the local convergence rates of PINNs and CTPINNs. For an $L+1$-layer MLP with $n_0=1$ and hidden widths $n_1,\ldots,n_L$, we use the width-normalized parameterization
\begin{equation}\label{eq: width_normalized_mlp}
\begin{aligned}
    h^{(0)}(\xi)&=\xi,\\
    h^{(\ell)}(\xi)
    &=\tanh\!\left(
    \frac{1}{\sqrt{n_{\ell-1}}}W^{(\ell)}h^{(\ell-1)}(\xi)
    +b^{(\ell)}\right),
    \qquad \ell=1,\ldots,L,\\
    v^\theta(\xi)
    &=\frac{1}{\sqrt{n_L}}W^{(L+1)\top}h^{(L)}(\xi)+b^{(L+1)}.
\end{aligned}
\end{equation}
At initialization, all entries of $W^{(\ell)}$ and $b^{(\ell)}$, for $\ell=1,\ldots,L+1$, are sampled independently from $\mathcal{N}(0,\sigma^2)$, where $\sigma^2>0$ is fixed independently of the layer widths. We analyze the limit $n_1,\ldots,n_L\to\infty$ in the frozen-kernel regime. Define
\begin{equation}
    \begin{aligned}
& \boldsymbol{r}_\alpha(\theta):=
\left(
\left\{\frac{\epsilon u^{\theta,\alpha}_{xx}(x_i)+u^{\theta,\alpha}_x(x_i)}{\sqrt{N_r}}\right\}_{i=1}^{N_r},
u^{\theta,\alpha}(0)-1,
u^{\theta,\alpha}(1)-e^{-1/\epsilon}
\right),\\
& \mathcal{R}_{\alpha,N_r+1}f:=f(0),\\
& \mathcal{R}_{\alpha,N_r+2}f:=f(1).
\end{aligned}
\end{equation}
and for $i=1,\ldots,N_r$, define
\begin{equation}
    \mathcal{R}_{\alpha,i}f: =\frac{1}{\sqrt{N_r}}\left[
    \epsilon\bigl(T_\alpha'(x_i)\bigr)^2
    f''(T_\alpha(x_i))
    +\bigl(\epsilon T_\alpha''(x_i)+T_\alpha'(x_i)\bigr)
    f'(T_\alpha(x_i))\right].
\end{equation}
Let
\begin{equation}
    \Theta(\xi,\zeta)
    :=\lim_{n_1,\ldots,n_L\to\infty}
    \left\langle\nabla_\theta v^\theta(\xi),
    \nabla_\theta v^\theta(\zeta)\right\rangle
\end{equation}
be the frozen infinite-width output NTK at initialization.
The frozen infinite-width residual NTK
$K_\alpha\in\mathbb{R}^{(N_r+2)\times(N_r+2)}$ is then defined by
\begin{equation}\label{eq: explicit_residual_ntk}
    (K_\alpha)_{ij}
    :=\mathcal{R}_{\alpha,i}^{(\xi)}
    \mathcal{R}_{\alpha,j}^{(\zeta)}
    \Theta(\xi,\zeta),
    \qquad i,j=1,\ldots,N_r+2.
\end{equation}
For $T_\alpha(x)$, the identity coordinate transformation is $T_\infty(x):=\lim_{\alpha\to\infty}T_\alpha(x)=x$. Accordingly, we define the residual NTK for the identity transformation, and hence for PINNs, as $K_I:=\lim_{\alpha\to\infty}K_\alpha$. Define
\begin{equation}
    \lambda_I:=\lambda_{\min}(K_I),
    \qquad
    \lambda_\alpha:=\lambda_{\min}(K_\alpha).
\end{equation}
For a fixed coordinate map, continuous-time gradient flow in the frozen-kernel
regime satisfies
\begin{equation}\label{eq:gradient_flow_frozen_kernel}
    \dot{\boldsymbol{r}}_\alpha(t)=-K_\alpha\boldsymbol{r}_\alpha(t),
    \qquad
    \|\boldsymbol{r}_\alpha(t)\|_2^2
    \le e^{-2\lambda_\alpha t}\|\boldsymbol{r}_\alpha(0)\|_2^2.
\end{equation}

\begin{proposition}\label{prop:boundary_layer_ntk_regular}
For the $L+1$-layer MLP described above, $K_\alpha$ is continuously differentiable in
$\alpha\in(0,\infty)$. Consequently, for $\bar\alpha>\epsilon$ and
$\mathcal{I}=[\epsilon,\bar\alpha]$, there exists
$C_{\mathcal{I}}>0$ such that
\begin{equation}
    \|K_\alpha-K_\beta\|_2
    \le C_{\mathcal{I}}|\alpha-\beta|,
    \qquad \alpha,\beta\in\mathcal{I}.
\end{equation}
\end{proposition}
Here, $C_{\mathcal{I}}$ is discussed in the proof.

Based on \cref{thm:equidistribution_optimality}, $\alpha=\epsilon$ is the
optimal choice for the coordinate transformation parameter, and its
corresponding residual NTK is $K_\epsilon$. Let
$N_{\mathrm{ntk}}:=N_r+2$, and define
\begin{equation}
\begin{aligned}
    \underline{\kappa}_\epsilon
    &:=\min_{1\le i\le N_{\mathrm{ntk}}}
    \left((K_\epsilon)_{ii}
    -\sum_{\substack{1\le j\le N_{\mathrm{ntk}}\\j\ne i}}
    |(K_\epsilon)_{ij}|\right),\\
    \eta_\epsilon&:=\underline{\kappa}_\epsilon-\lambda_I.
\end{aligned}
\end{equation}

\begin{proposition}\label{prop:boundary_layer_ntk_gain}
If the collocation points set $\mathcal{X}$ satisfies
\begin{equation}\label{eq: tangent_feature_separation}
    \underline{\kappa}_\epsilon>\lambda_I,
\end{equation}
then $\eta_\epsilon>0$ and $\lambda_\epsilon \ge \lambda_I+\eta_\epsilon >\lambda_I$.
\end{proposition}
Furthermore, for the $i$th
standard basis vector $e_i\in\mathbb{R}^{N_{\mathrm{ntk}}}$, define the
canonical Gram feature associated with the $i$th residual component by
\begin{equation}
    g_{\alpha,i}:=K_\alpha^{1/2}e_i,
    \qquad i=1,\ldots,N_{\mathrm{ntk}}.
\end{equation}
Here, $K_\alpha^{1/2}$ denotes the symmetric positive semidefinite
square root of $K_\alpha$.
\begin{lemma}\label{lemma: unique_K_root}
    For every $\alpha>0$, the residual NTK $K_\alpha$ is real symmetric and positive semidefinite. Hence, by the spectral theorem, if $K_\alpha=Q\Lambda Q^\top$ with $\Lambda\ge0$, its unique symmetric positive semidefinite square root is $K_\alpha^{1/2}=Q\Lambda^{1/2}Q^\top$.
\end{lemma}
Then
\begin{equation}
    (K_\alpha)_{ij}
    =e_i^\top K_\alpha e_j
    =\langle g_{\alpha,i},g_{\alpha,j}\rangle,
    \qquad i,j=1,\ldots,N_{\mathrm{ntk}}.
\end{equation}
Using these features, define
\begin{equation}\label{eq: tangent_feature_collocation_condition}
\begin{aligned}
    m_\epsilon
    &:=\min_{1\le i\le N_{\mathrm{ntk}}}
    \|g_{\epsilon,i}\|_2,\\
    M_\epsilon
    &:=\max_{1\le i\le N_{\mathrm{ntk}}}
    \|g_{\epsilon,i}\|_2,\\
    \mu_\epsilon
    &:=\max_{\substack{1\le i,j\le N_{\mathrm{ntk}}\\i\ne j}}
    \frac{|\langle g_{\epsilon,i},g_{\epsilon,j}\rangle|}
    {\|g_{\epsilon,i}\|_2\|g_{\epsilon,j}\|_2}.
\end{aligned}
\end{equation}
Then we have the following corollary.
\begin{corollary}
\label{cor:boundary_layer_ntk_coherence}
If $(K_\epsilon)_{ii}>0$ for every $i=1,\ldots,N_{\mathrm{ntk}}$, equivalently
$g_{\epsilon,i}\ne0$ and
\begin{equation}\label{eq: tangent_feature_coherence_condition}
    m_\epsilon^2-(N_{\mathrm{ntk}}-1)
    \mu_\epsilon M_\epsilon^2>\lambda_I,
\end{equation}
then \cref{eq: tangent_feature_separation} holds, and consequently
    $\lambda_\epsilon>\lambda_I$.
\end{corollary}
The rowwise condition in \cref{eq: tangent_feature_separation} permits the
feature norms and correlations to vary across collocation points. The more
conservative condition in \cref{eq: tangent_feature_coherence_condition}
instead requires uniformly large norms and low mutual coherence.

\begin{theorem}\label{thm:boundary_layer_ntk_acceleration}
Let $\mathcal{I}=[\epsilon,\bar\alpha]$ for some
$\bar\alpha>\epsilon$. Then, for every $\alpha\in\mathcal{I}$,
\begin{equation}\label{eq: boundary_layer_ntk_rate}
    \lambda_\alpha
    \ge \lambda_I+\eta_\epsilon
    -C_{\mathcal{I}}(\alpha-\epsilon).
\end{equation}
Furthermore, if the collocation points condition in
\cref{eq: tangent_feature_separation} holds, then $\eta_\epsilon>0$ and,
for any $\alpha\in\mathcal{I}$ satisfying
$\alpha-\epsilon<\eta_\epsilon/C_{\mathcal{I}}$, we have
$\lambda_\alpha>\lambda_I$. 
\end{theorem}
By \cref{eq:gradient_flow_frozen_kernel}, the inequality $\lambda_\alpha>\lambda_I$ means that CTPINNs have a strictly faster NTK convergence rate than PINNs for all $\alpha\in\mathcal{I}$, as well as means that CTPINNs mitigate spectral bias effects in the training dynamics.

\subsection{Implicit Tendency towards Gradient Equidistribution}\label{sec: implicit_tendency_t_gradient_equid}
First, define
$u^{\theta(t),\alpha(t)}(x)=v^{\theta(t)}(T_{\alpha(t)}(x))$ and the
physical-domain PDE residual
\begin{equation}\label{eq: trajectory_pde_residual}
    R(t,x):=\epsilon u^{\theta(t),\alpha(t)}_{xx}(x)
    +u^{\theta(t),\alpha(t)}_x(x).
\end{equation}
Since $T_\alpha''=-T_\alpha'/\alpha$, the chain rule gives
\begin{equation}\label{eq: residual_curvature_decomposition}
\begin{aligned}
    R(t,x)
    ={}&\epsilon\left(T_{\alpha(t)}'(x)\right)^2
    v_{\xi\xi}^{\theta(t)}(T_{\alpha(t)}(x))\\
    &+\left(1-\frac{\epsilon}{\alpha(t)}\right)
    T_{\alpha(t)}'(x)v_\xi^{\theta(t)}(T_{\alpha(t)}(x)).
\end{aligned}
\end{equation}

We use the following successful-training assumption for the coordinate
convergence analysis.

\begin{assumption}\label{assump:stable_training}
The original training objective converges to zero:
\begin{equation}
    \mathcal L(\theta(t),\alpha(t))\to0
    \qquad\text{as }t\to\infty.
\end{equation}
\end{assumption}

However, $\mathcal L(\theta(t),\alpha(t))$ need not equal zero at any finite time. If it did, there would be no gradient signal to drive further updates of $\alpha(t)$.

For the analytical solution $u^*(x)=e^{-x/\epsilon}$, the unique coordinate transformation that equidistributes the gradient, as given by
\cref{thm:equidistribution_optimality}, is
\begin{equation}
    T^*(x)
    =\frac{\int_0^x |(u^*)'(s)|\,\mathrm ds}
    {\int_0^1 |(u^*)'(s)|\,\mathrm ds}
    =\frac{1-e^{-x/\epsilon}}{1-e^{-1/\epsilon}}
    =T_\epsilon(x).
\end{equation}
Consequently, the exact transformed solution is affine:
\begin{equation}\label{eq: target_affine_profile}
    v_{\mathrm{aff}}(\xi)
    :=u^*((T^*)^{-1}(\xi))
    =1-(1-e^{-1/\epsilon})\xi.
\end{equation}

Because the target solution in computational coordinates is affine, it is represented
exactly by the one-layer linear MLP
\begin{equation}\label{eq: affine_one_layer_mlp}
    v^{a,b}(\xi)=a+b\xi.
\end{equation}

\cref{thm:affine_alpha_convergence} shows that the coordinate parameter $\alpha(t)$ converges to $\epsilon$ when the CTPINN is trained with the one-layer linear MLP.

\begin{theorem}
\label{thm:affine_alpha_convergence}
Consider gradient flow for the full loss in \cref{eq: unit_weight_boundary_layer_loss} with the one-layer MLP \cref{eq: affine_one_layer_mlp}. If $\alpha(0)\ge\epsilon$ and
\cref{assump:stable_training} holds, then
\begin{equation}
    \alpha(t)\to\epsilon
    \qquad\text{as }t\to\infty.
\end{equation}
\end{theorem}
\cref{thm:affine_alpha_convergence} relies on prior information that the transformed solution is affine and that $\alpha=\epsilon$. However, if it is not known
a priori whether training will drive the learned CT to the exact optimal
transformation, the transformed solution need not be
affine. We therefore follow standard deep-learning practice and use a
multilayer MLP as the solution model rather than imposing an affine profile.
The following convergence result applies under successful training when the
non-affine defect is controlled relative to the coordinate mismatch. Unlike the MLP in \cref{sec: ntk_convergence_analysis}, this MLP is defined as
\begin{equation}\label{eq: trajectory_tanh_mlp}
\begin{aligned}
    h^{(0)}(\xi)&=\xi,\\
    h^{(\ell)}(\xi)
    &=\tanh\!\left(B^{(\ell)}h^{(\ell-1)}(\xi)+c^{(\ell)}\right),
    \qquad \ell=1,\ldots,L,\\
    v^\theta(\xi)&=\beta^\top h^{(L)}(\xi)+c^{(L+1)},
\end{aligned}
\end{equation}
and decompose its output into the affine part and the non-affine part:
\begin{equation}\label{eq: mlp_affine_decomposition}
\begin{aligned}
    \ell^\theta(\xi)
    &:=(1-\xi)v^\theta(0)+\xi v^\theta(1),\\
    w^\theta(\xi)&:=v^\theta(\xi)-\ell^\theta(\xi),\\
    \delta_3(\theta)&:=\|w^\theta\|_{C^3([0,1])}.
\end{aligned}
\end{equation}
Here $\delta_3$ measures the deviation from the affine regime. 

\begin{theorem}
\label{thm:alpha_descent}
If there exist $C_{\rho}>0$ and $t_0\ge0$ such that \cref{assump:stable_training} holds and, for all $t\ge t_0$, $\alpha(t)\ge\epsilon$ and $\delta_3(\theta(t))\le C_{\rho}\left|1-\frac{\epsilon}{\alpha(t)}\right|$,
then
\begin{equation}
    \alpha(t)\to\epsilon
    \qquad\text{as }t\to\infty.
\end{equation}
\end{theorem}
\begin{remark}
The threshold $C_{\rho}$ depends on $\epsilon$ and the collocation points.
An explicit sufficient choice is provided in the proof.
\end{remark}
The proportional condition uses the normalized mismatch
$|1-\epsilon/\alpha(t)|=|\alpha(t)-\epsilon|/\alpha(t)$; it requires neither a local
neighborhood nor an upper bound on $\alpha$. The one-sided trajectory condition $\alpha(t)\ge\epsilon$ excludes the under-resolved limit $\alpha(t)\rightarrow 0$. For a non-affine MLP, this condition is assumed rather than deduced from the initialization.

\subsection{Error Analysis}\label{sec: error_analysis}
We finally analyze how the approximation error in the computational domain transfers back to the physical domain for the same boundary-layer problem.
Define the exact transformed solution
\begin{equation}
    v^*(\xi):=u^*(T_\alpha^{-1}(\xi)), \quad \xi \in [0,1].
\end{equation}
Define $m_\alpha:=\min_{x\in[0,1]}T_\alpha'(x)$ and $M_\alpha:=\max_{x\in[0,1]}T_\alpha'(x)$.
\begin{theorem}\label{thm:error_transfer}
Fix $\epsilon>0$ and $\alpha>0$. 
Since $T_\alpha$ is smooth and strictly increasing on $[0,1]$, $0<m_\alpha< M_\alpha<\infty$. Then
\begin{equation}
    \|u^*-u^{\theta,\alpha}\|_{L^2(0,1)}
    \le C_{1,\alpha}\|v^*-v^\theta\|_{L^2(0,1)},
\end{equation}
and
\begin{equation}
    \|(u^*-u^{\theta,\alpha})'\|_{L^2(0,1)}
    \le C_{2,\alpha}\|(v^*-v^\theta)'\|_{L^2(0,1)},
\end{equation}
where $C_{1,\alpha}=m_\alpha^{-1/2}$ and
$C_{2,\alpha}=M_\alpha^{1/2}$. Consequently,
\begin{equation}
    \|u^*-u^{\theta,\alpha}\|_{H^1(0,1)}
    \le C_{T_\alpha}\|v^*-v^\theta\|_{H^1(0,1)},
\end{equation}
where $C_{T_\alpha}:=
    \min\left\{
    \max\{C_{1,\alpha},C_{2,\alpha}\},
    \sqrt{2+M_\alpha}
    \right\}$.
\end{theorem}

\begin{remark}
For the prescribed normalized exponential transformation in
\cref{eq: normalized_exponential_ct},
\begin{equation}
    T_\alpha'(x)
    =\frac{\exp(-x/\alpha)}{\alpha(1-\exp(-1/\alpha))},
\end{equation}
which is strictly decreasing on $[0,1]$. Hence the constants in
\cref{thm:error_transfer} are explicitly
\begin{equation}
    m_\alpha
    =\frac{1}{\alpha(\exp(1/\alpha)-1)},
    \qquad
    M_\alpha
    =\frac{1}{\alpha(1-\exp(-1/\alpha))}.
\end{equation}
Consequently,
\begin{equation}
    C_{1,\alpha}
    =\sqrt{\alpha(\exp(1/\alpha)-1)},
    \qquad
    C_{2,\alpha}
    =\frac{1}
    {\sqrt{\alpha(1-\exp(-1/\alpha))}}.
\end{equation}
Moreover, $2\alpha\sinh(1/(2\alpha))>1$ implies
$C_{1,\alpha}>C_{2,\alpha}$. Hence
\begin{equation}
    C_{T_\alpha}
    =\min\left\{
    \sqrt{\alpha(\exp(1/\alpha)-1)},
    \sqrt{2+\frac{1}{\alpha(1-\exp(-1/\alpha))}}
    \right\}.
\end{equation}
In particular, the second bound grows only as $O(\alpha^{-1/2})$ as
$\alpha\rightarrow 0$, avoiding the exponential growth produced by a direct
lower-Jacobian estimate. 
\end{remark}

\cref{thm:error_transfer} shows that once the network reduces the approximation error in the computational domain, the physical-domain error remains controlled by a constant determined by $T_\alpha$. Since this constant is independent of the number of collocation points $N_r$, any convergence rate in $N_r$ achieved by $v^\theta$ is preserved by the physical-domain CTPINN prediction.

%% file: main/experiments.tex
\section{Experiments}
We first study explicit coordinate transformations in \cref{sec: explicit_construction_exp}. To evaluate their adaptiveness (\textit{i.e.}, their implicit tendency toward gradient equidistribution), we consider 1D and 2D boundary-layer problems in \cref{sec: adaptiveness_ctpinns} and \cref{sec: 2D_boundary_layer_problems}, respectively, and then compare CTPINNs with adaptive sampling methods on Poisson equations in \cref{sec: poisson_equations}. We next investigate implicit constructions in \cref{sec: implicit_construction_exp}, focusing on the effectiveness of the CT module when integrated into advanced models. Specifically, we combine the CT module with PirateNet~\cite{wang2024piratenets}, SPINN~\cite{cho2023separable}, and QRPINN~\cite{yu2026quasi} and evaluate its impact on their respective benchmarks. Finally, we analyze the computational cost of CTPINNs in \cref{sec: cost_analysis}, showing that the lightweight module is efficient and introduces no significant overhead. The hyperparameter details are provided in \cref{appendix: experimental_details} and the repository \url{https://github.com/DUCH714/CTPINN}, which will be made public upon publication.

We use the relative $L^2$ error as the primary evaluation metric and highlight the minimum error in blue in the tables. For a reference solution $u$ and a model prediction $u^\theta$, the relative $L^2$ error is defined as
\begin{equation}
\text{Relative } L^2 = \frac{\|u-u^\theta\|_{L^2(\Omega)}}{\|u\|_{L^2(\Omega)}}.
\end{equation}
In discrete form, this corresponds to replacing the $L^2$ norms by Euclidean norms on evaluation points. Smaller values indicate better approximation accuracy.

\subsection{Explicit Construction}\label{sec: explicit_construction_exp}
To verify the adaptiveness of CTPINNs and evaluate the performance of explicit constructions, we conduct experiments on 1D and 2D boundary-layer problems, for which optimal coordinate transformations are known. Additionally, we conduct experiments on 2D Poisson equations to evaluate the performance and robustness of CTPINNs when optimal coordinate transformations are unavailable and to compare them with adaptive sampling methods. For 2D problems, we construct the coordinate transformation by applying majorants independently along each dimension.
\subsubsection{1D Boundary Layer Problems}\label{sec: adaptiveness_ctpinns}
This section focuses on the 1D boundary-layer problem in \cref{eq: bl1d}, a classical benchmark for evaluating methods that resolve sharp gradients near boundaries. By \cref{thm:equidistribution_optimality}, the optimal coordinate transformation for \cref{eq: bl1d} is given by \cref{eq: majorant_1} with $\alpha=\epsilon$.

The results are presented in \cref{tab:result_bl1d}, where 'N' and 'A' denote results with fixed and adaptive $\alpha$, respectively. For $\epsilon \geq 10^{-3}$, the PINN performs well and successfully solves this simple problem. However, when $\epsilon < 10^{-3}$, the PINN fails to produce accurate predictions. This challenge has also been reported in prior studies on boundary-layer problems~\cite{gie2024semi,yu2024sinc,shen2026matched}, which propose methods specifically designed to address this class of problems. Because this experiment is intended to evaluate the performance and adaptiveness of CTPINNs, and our objective is to develop a problem-agnostic framework rather than a method targeting boundary-layer problems specifically, we do not compare against these methods. Nevertheless, we argue that our explicit construction is a competitive approach for solving these problems compared with the aforementioned methods. Additionally, our coordinate transformation module is orthogonal to some of these approaches and can be combined with them to further improve performance.

\begin{figure}[ht]
\centering
\subcaptionbox{$\epsilon=10^{-2}$ \label{fig:bl1d_1e-2}}[0.32\textwidth]{%
\includegraphics[width=\linewidth]{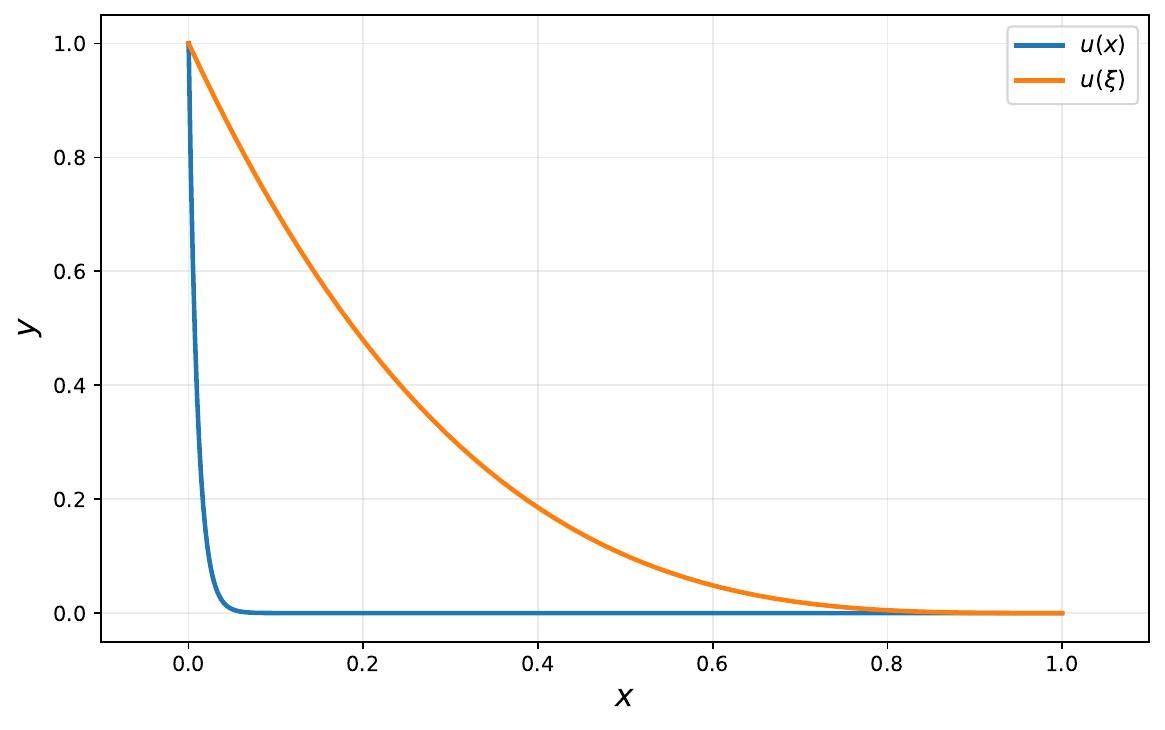}}
\subcaptionbox{$\epsilon=10^{-4}$ \label{fig:bl1d_1e-4}}[0.32\textwidth]{%
\includegraphics[width=\linewidth]{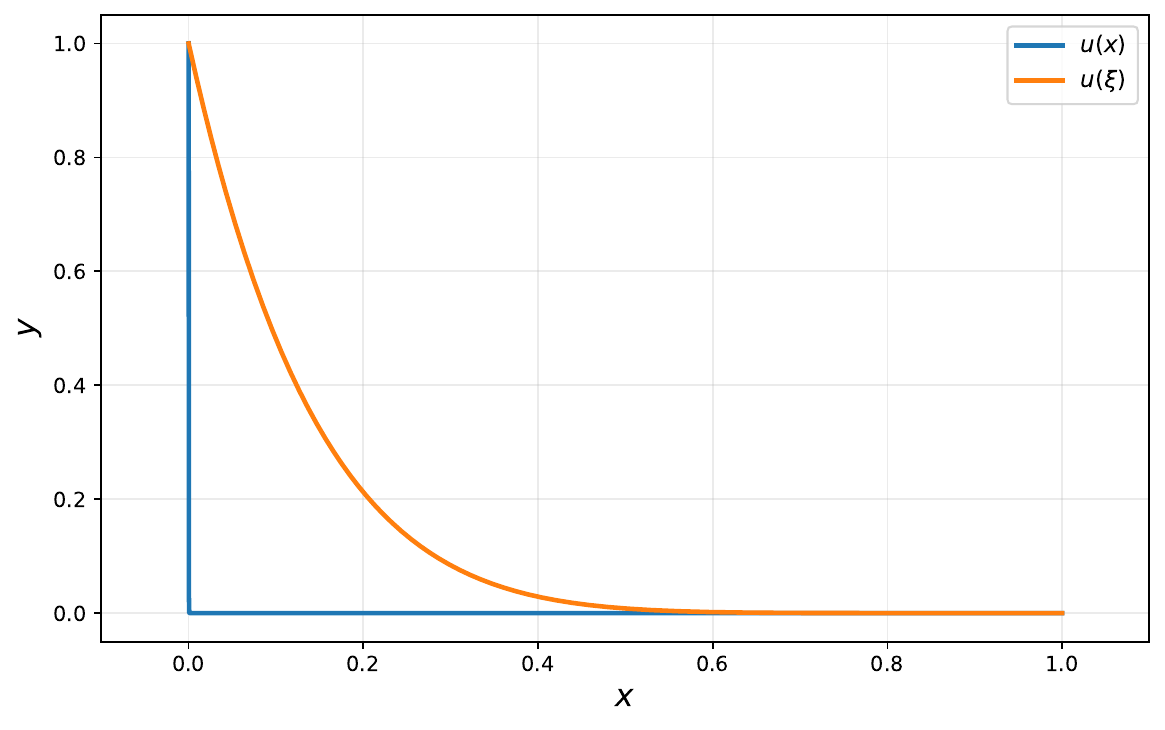}}
\subcaptionbox{$\epsilon=10^{-5}$\label{fig:bl1d_1e-5}}[0.32\textwidth]{%
\includegraphics[width=\linewidth]{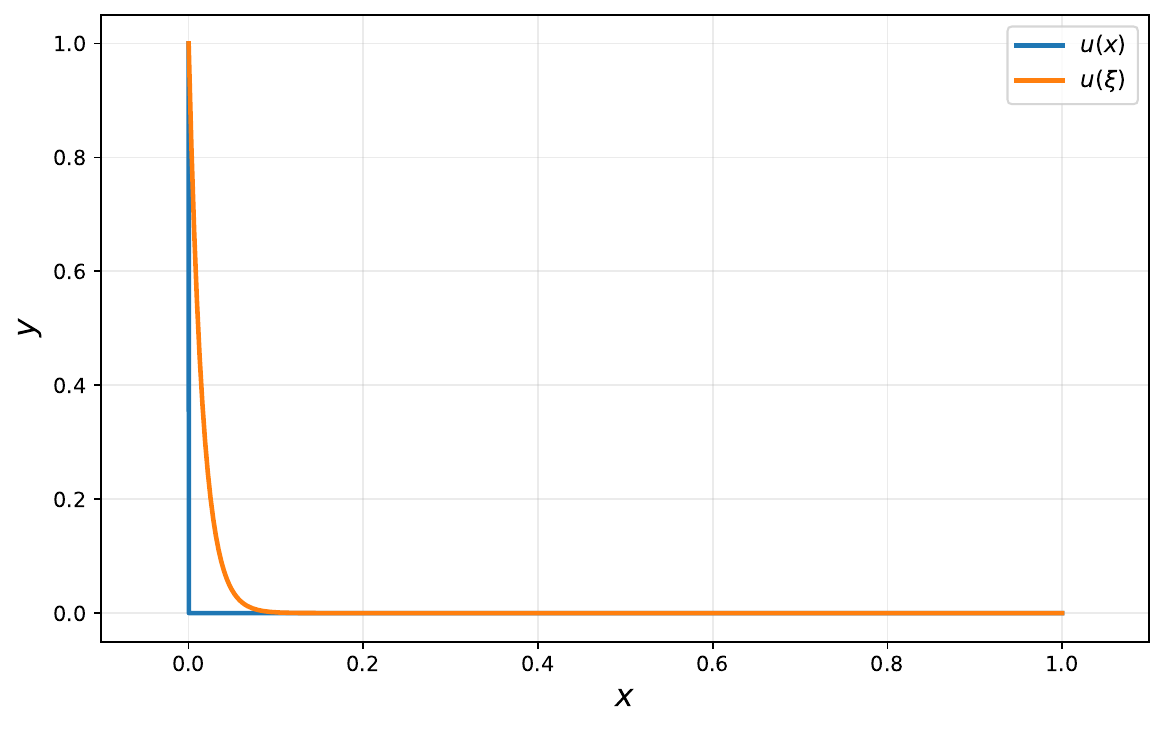}}
\caption{Reference solutions and reference transformed solutions (adaptive $T_{\psi_1}$) of the 1D boundary layer problem for different values of $\epsilon$. The final $\alpha$ is $3.33\times 10^{-2}$, $6.30\times 10^{-4}$, and $6.94\times 10^{-4}$ for $\epsilon=10^{-2}$, $10^{-4}$, and $10^{-5}$, respectively. The transformed solutions stretch the boundary layer region to make the gradient more equidistributed and facilitate learning.}
\label{fig: bl1d_solution}
\end{figure}

\begin{figure}[ht]
\centering
\subcaptionbox{Dynamic $\alpha$ \label{fig:bl1d_dynamic_alpha}}[0.32\textwidth]{%
\includegraphics[width=\linewidth]{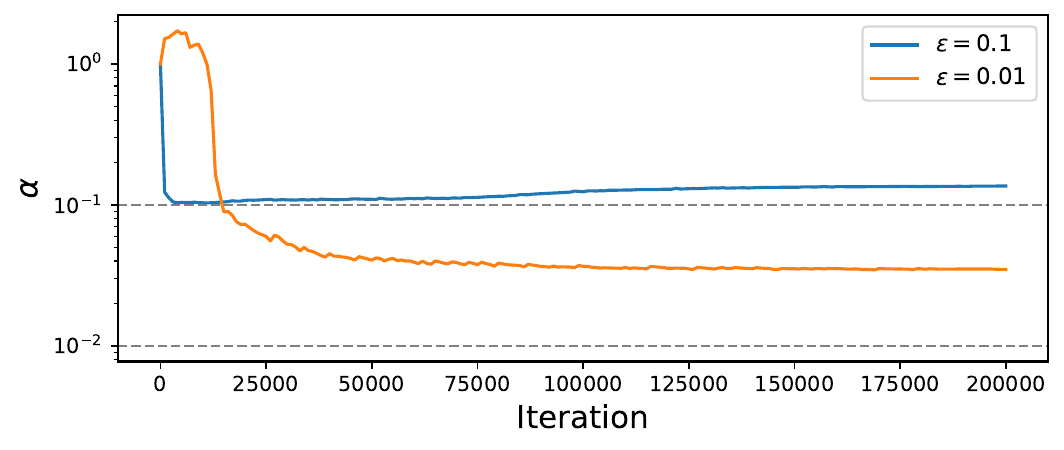}}
\subcaptionbox{Loss for $\epsilon=0.1$ \label{fig:bl1d_loss_0.1}}[0.32\textwidth]{%
\includegraphics[width=\linewidth]{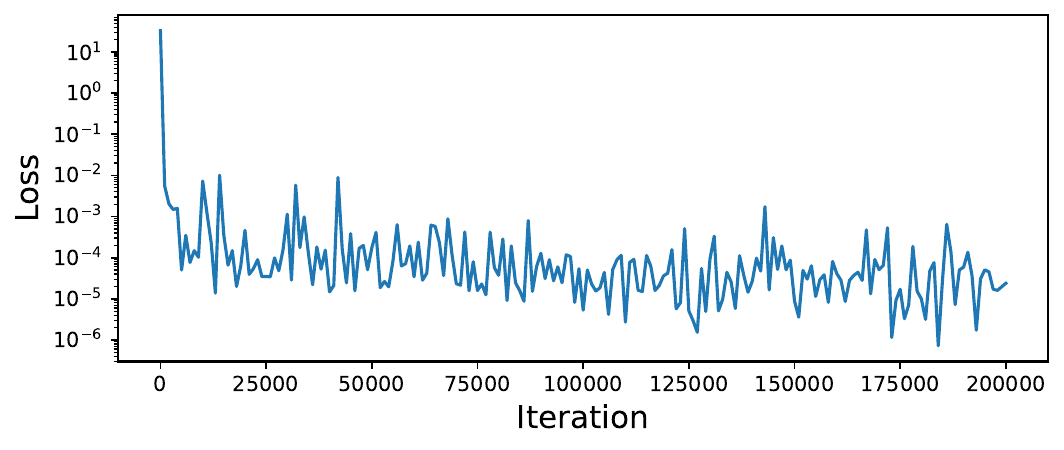}}
\subcaptionbox{Loss for $\epsilon=0.01$ \label{fig:bl1d_loss_0.01}}[0.32\textwidth]{%
\includegraphics[width=\linewidth]{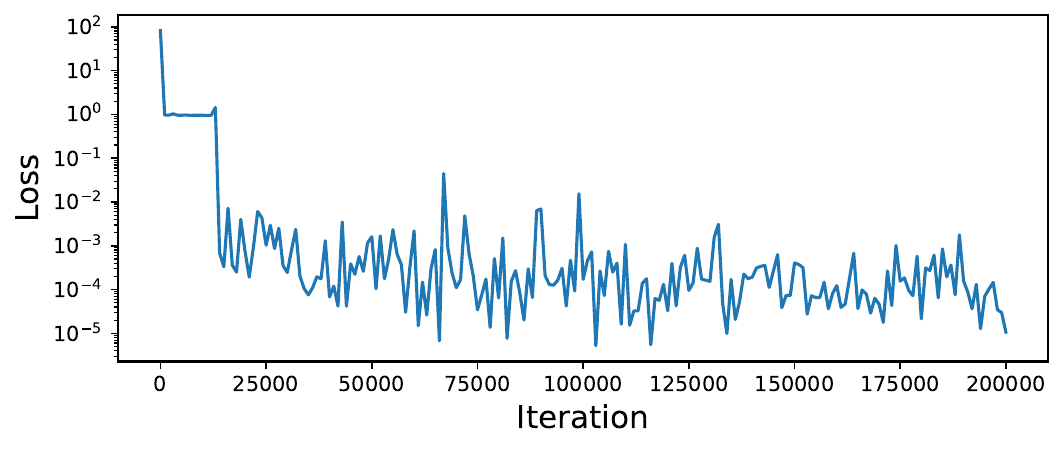}}
\caption{Training details for the 1D boundary-layer problem with $\epsilon=0.1$ and $0.01$. (a) shows the evolution of $\alpha$ for both values of $\epsilon$, indicating that $\alpha$ tends toward $\epsilon$ during training. (b) and (c) show the loss curves for $\epsilon=0.1$ and $\epsilon=0.01$, respectively.}
\label{fig: bl1d}
\end{figure}

\cref{tab:result_bl1d} shows that CTPINNs can successfully solve the problem, even when $\epsilon$ is as small as $10^{-6}$, with both fixed and adaptive $\alpha$, demonstrating their strong performance. For this boundary-layer problem, Majorant1 \cref{eq: majorant_1} is the theoretically optimal choice. However, Majorant2 \cref{eq: majorant_2} and Majorant3 \cref{eq: majorant_3} also achieve comparable results, demonstrating that CTPINNs do not necessarily require optimal CT functions. Notably, when $\epsilon=10^{-6}$, the best results are achieved by adaptive Majorant1, for which the final $\alpha=6.94\times 10^{-4}$ rather than the theoretical optimum $\alpha^*=10^{-6}$. Computing second-order derivatives with Majorant1 and $\alpha=10^{-6}$ causes numerical overflow, so using a larger but still effective $\alpha$ is more stable for CTPINNs. This also demonstrates the robustness of our adaptive method, which can find a numerically effective $\alpha$ rather than relying on a theoretical optimum that is unstable to evaluate. The reference solutions and their transformed counterparts are compared in \cref{fig: bl1d_solution}, illustrating the implicit tendency. In addition, \cref{fig: bl1d} shows the evolution of $\alpha(t)$ and the loss during training.
\begin{table}[ht]
\centering
\caption{Results on 1D boundary layer problems}
\label{tab:result_bl1d}
\begin{adjustbox}{width=\textwidth}
\begin{tabular}{lccccccc}
\toprule
$\epsilon$ & PINN & \multicolumn{2}{c}{PINN+$T_{\psi_1}$} & \multicolumn{2}{c}{PINN+$T_{\psi_2}$} & \multicolumn{2}{c}{PINN+$T_{\psi_3}$} \\
 &  & N & A & N & A & N & A \\
\midrule
$10^{-3}$ & $1.60 \times 10^{-2}$ & $8.70 \times 10^{-4}$ & \cellcolor{blue!25}$5.19 \times 10^{-4}$ & $3.33 \times 10^{-3}$ & $4.71 \times 10^{-3}$ & $9.06 \times 10^{-3}$ & $9.18 \times 10^{-3}$ \\
$10^{-4}$ & $1.05 \times 10^{1}$ & \cellcolor{blue!25}$2.83 \times 10^{-4}$ & $8.71 \times 10^{-4}$ & $9.01 \times 10^{-3}$ & $4.56 \times 10^{-3}$ & $5.88 \times 10^{-3}$ & $7.43 \times 10^{-2}$ \\
$10^{-5}$ & $>1$ & \cellcolor{blue!25}$1.98 \times 10^{-4}$ & $1.44 \times 10^{-3}$ & $1.88 \times 10^{-2}$ & $6.71 \times 10^{-3}$ & $2.88 \times 10^{-2}$ & $2.17 \times 10^{0}$ \\
$10^{-6}$ & $>1$ & $1.42 \times 10^{-2}$ & \cellcolor{blue!25}$1.74 \times 10^{-3}$ & $3.96 \times 10^{-3}$ & $5.54 \times 10^{-3}$ & $1.98 \times 10^{-1}$ & $>1$ \\
\bottomrule
\end{tabular}
\end{adjustbox}
\end{table}

\subsubsection{2D Boundary Layer Problems}\label{sec: 2D_boundary_layer_problems}
We next extend our study to the 2D boundary-layer problem defined by the following PDE:
\begin{equation}\label{eq: bl2d}
\begin{aligned}
     \epsilon (u_{xx} + u_{yy}) + u_x + u_y &= 0, \quad (x,y) \in [0,1]^2,\\
\end{aligned}
\end{equation}
with Dirichlet boundary conditions and the analytical solution $u^*(x,y)=e^{-(x+y)/\epsilon}$. The solution and its transformed counterpart are shown in \cref{fig:framework}. \cref{tab:result_bl2d} summarizes the results for different values of $\epsilon$. Majorant1 with fixed $\alpha$, which is also the theoretically optimal choice for this problem, consistently achieves the best results. However, the other approaches also achieve comparable results, and all CTPINN variants outperform PINNs.
\begin{table}[ht]
\centering
\caption{Results on 2D boundary layer problems}
\label{tab:result_bl2d}
\begin{adjustbox}{width=\textwidth}
\begin{tabular}{lccccccc}
\toprule
$\epsilon$ & PINN & \multicolumn{2}{c}{PINN+$T_{\psi_1}$} & \multicolumn{2}{c}{PINN+$T_{\psi_2}$} & \multicolumn{2}{c}{PINN+$T_{\psi_3}$} \\
 &  & N & A & N & A & N & A \\
\midrule
$1$ & $4.48 \times 10^{-6}$ & \cellcolor{blue!25}$2.98 \times 10^{-6}$ & $4.32 \times 10^{-6}$ & $7.62 \times 10^{-6}$ & $7.72 \times 10^{-6}$ & $4.95 \times 10^{-6}$ & $3.79 \times 10^{-6}$ \\ 
$10^{-1}$ & $5.74 \times 10^{-5}$ & \cellcolor{blue!25}$3.99 \times 10^{-5}$ & $1.90 \times 10^{-4}$ & $1.45 \times 10^{-4}$ & $1.74 \times 10^{-4}$ & $1.30 \times 10^{-4}$ & $9.95 \times 10^{-5}$ \\ 
$10^{-2}$ & $9.71 \times 10^{-4}$ & \cellcolor{blue!25}$4.80 \times 10^{-5}$ & $1.69 \times 10^{-4}$ & $1.18 \times 10^{-3}$ & $6.06 \times 10^{-4}$ & $5.32 \times 10^{-4}$ & $6.26 \times 10^{-4}$ \\ 
$10^{-3}$ & $3.70 \times 10^{-2}$ & \cellcolor{blue!25}$1.05 \times 10^{-4}$ & $1.92 \times 10^{-4}$ & $5.10 \times 10^{-3}$ & $2.46 \times 10^{-3}$ & $7.07 \times 10^{-3}$ & $8.01 \times 10^{-3}$ \\ 
$10^{-4}$ & $8.10 \times 10^{0}$ & \cellcolor{blue!25}$2.01 \times 10^{-4}$ & $2.48 \times 10^{-3}$ & $9.71 \times 10^{-3}$ & $2.28 \times 10^{-2}$ & $6.69 \times 10^{-3}$ & $5.49 \times 10^{-1}$ \\ 
$10^{-5}$ & $>1$ & \cellcolor{blue!25}$8.17 \times 10^{-4}$ & $3.34 \times 10^{-3}$ & $4.69 \times 10^{-2}$ & $2.42 \times 10^{-2}$ & $1.45 \times 10^{-2}$ & $1.51 \times 10^{0}$ \\ 
\bottomrule
\end{tabular}
\end{adjustbox}
\end{table}

\subsubsection{Poisson's Equations}\label{sec: poisson_equations}
In this section, we consider the following Poisson equation:
\begin{equation}
\label{eq: poisson}
\begin{aligned}
	 \Delta u(\boldsymbol{x}) &= f(\boldsymbol{x}), \quad  \boldsymbol{x} \in \Omega \subset \mathbb{R}^d,\\
	 u(\boldsymbol{x}) &= g(\boldsymbol{x}), \quad \boldsymbol{x} \in \partial\Omega.
\end{aligned}
\end{equation}
Here, we set $\Omega=[0,1]^2$, define the solution as $u(\boldsymbol{x})=e^{-\epsilon\|\boldsymbol{x}-0.5\|_2^2/2}$ with $g=u|_{\partial\Omega}$, and set $f(\boldsymbol{x})=\left(\epsilon^2\|\boldsymbol{x}-0.5\|_2^2-2\epsilon\right)e^{-\epsilon\|\boldsymbol{x}-0.5\|_2^2/2}$. \cref{eq: poisson} is a classical benchmark problem for adaptive sampling methods for PINNs. Thus, in these experiments, we compare CTPINNs with the adaptive sampling methods RAD~\cite{wu2023comprehensive} and ACLE~\cite{lau2024pinnacle}. We use adaptive $T_{\psi_1}$ (denoted by $T_{\psi_1}$(A)) as the CT function.

\begin{table}[ht]
\centering
\caption{Results on 2D Poisson's equations}
\label{tab:result_poisson2d}
\begin{tabular}{lccc}
\toprule
& $\epsilon=1000$ & $\epsilon=2000$ & $\epsilon=5000$ \\
\midrule
PINN & $3.91 \times 10^{-2}$ & $1.71 \times 10^{-2}$ & $8.14 \times 10^{-2}$  \\
RAD & $2.27 \times 10^{-2}$ & $4.15 \times 10^{-2}$ &  $2.85 \times 10^{-1}$ \\
ACLE & $1.62 \times 10^{-2}$ & $3.20 \times 10^{-2}$ & $7.73 \times 10^{-2}$  \\
PINN+$T_{\psi_1}$(A) & \cellcolor{blue!25}$1.45 \times 10^{-2}$ & \cellcolor{blue!25}$1.49 \times 10^{-2}$ & \cellcolor{blue!25}$2.74 \times 10^{-2}$\\
\bottomrule
\end{tabular}
\end{table}

\begin{figure}[ht]
\centering
\subcaptionbox{Reference $u(x,y)$ \label{fig:poisson_rel}}[0.32\textwidth]{%
\includegraphics[width=\linewidth]{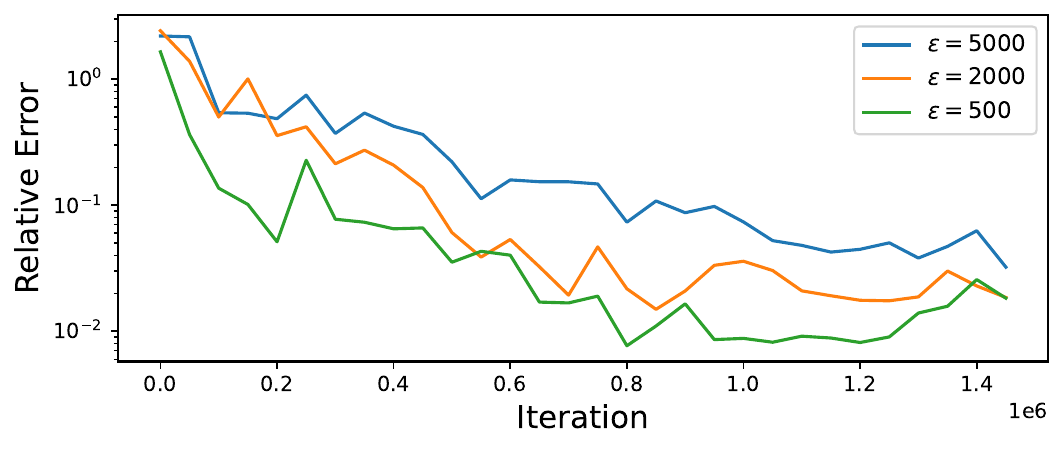}}
\subcaptionbox{Predicted $u(x,y)$ \label{fig:poisson_eps_x}}[0.32\textwidth]{%
\includegraphics[width=\linewidth]{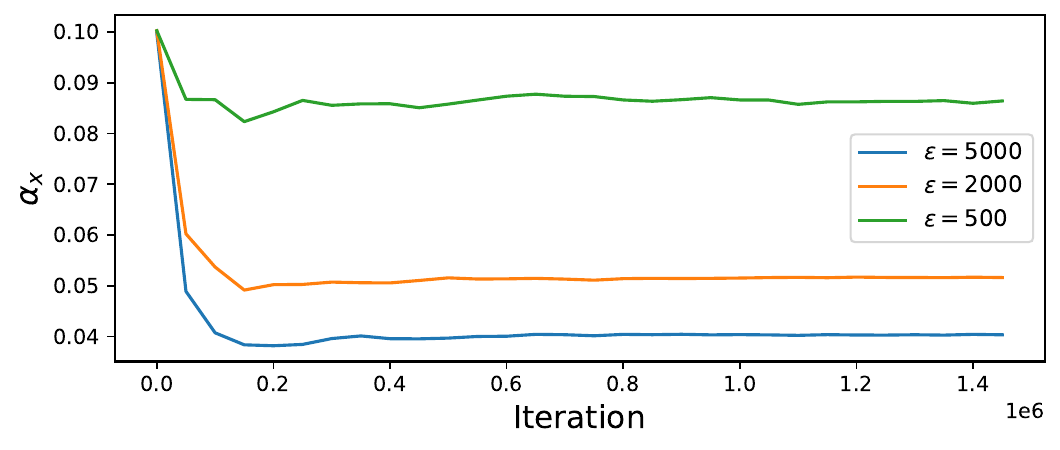}}
\subcaptionbox{Reference $u(\xi,\eta)$\label{fig:poisson_eps_y}}[0.32\textwidth]{%
\includegraphics[width=\linewidth]{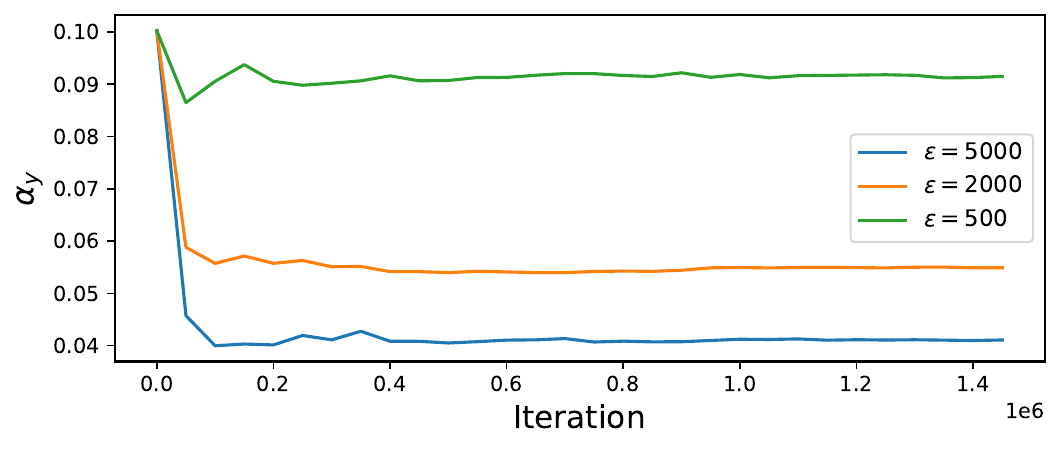}}
\caption{Training dynamics for solving the 2D Poisson equation in \cref{eq: poisson} with CTPINNs and adaptive $T_{\psi_1}$. (a) shows the relative error. (b) and (c) show the evolution of $\alpha_x$ and $\alpha_y$, respectively. The parameter $\alpha_x$ is the learnable parameter of Majorant1 \cref{eq: majorant_1} in the $x$-direction, and $\alpha_y$ is the corresponding parameter in the $y$-direction. As $\epsilon$ increases, the solution gradients become steeper. Accordingly, $\alpha_x$ and $\alpha_y$ converge to smaller values, yielding coordinate transformations that promote gradient equidistribution.}
\label{fig: poisson_2d_alpha}
\end{figure}

\begin{figure}[ht]
\centering
\subcaptionbox{Reference $u(x,y)$ \label{fig:poisson_2d_ana_physical}}[0.24\textwidth]{%
\includegraphics[width=\linewidth]{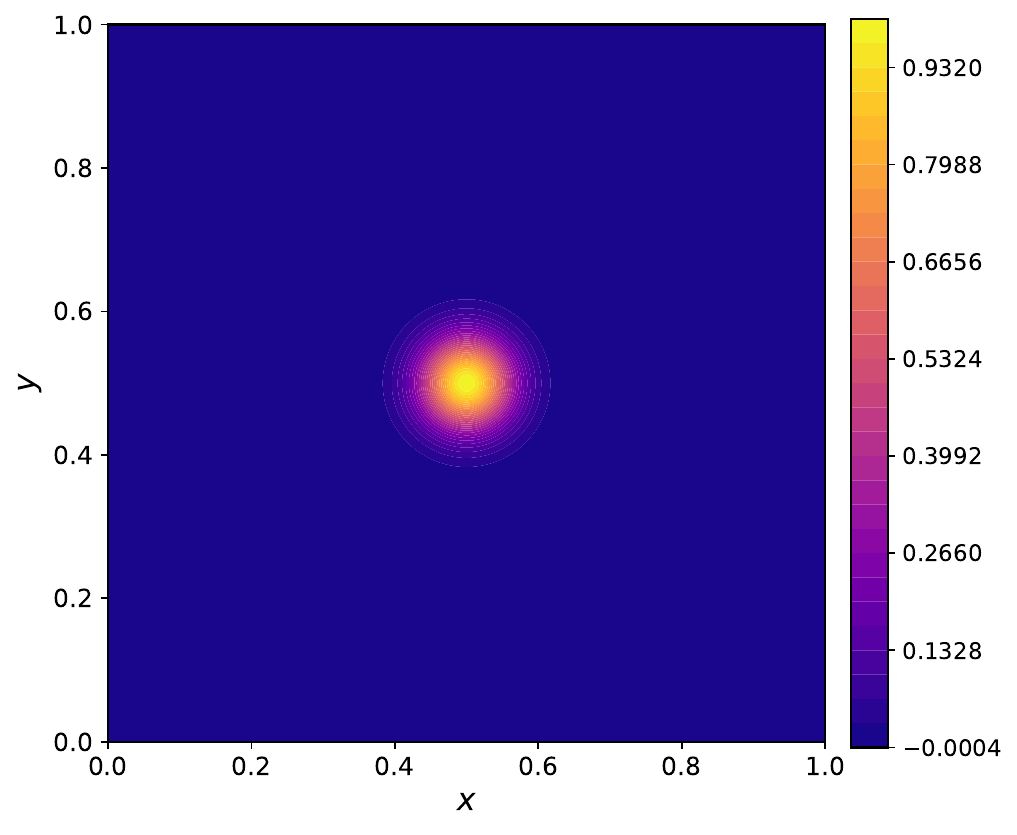}}
\subcaptionbox{Predicted $u(x,y)$ \label{fig:poisson_2d_pred_physical}}[0.24\textwidth]{%
\includegraphics[width=\linewidth]{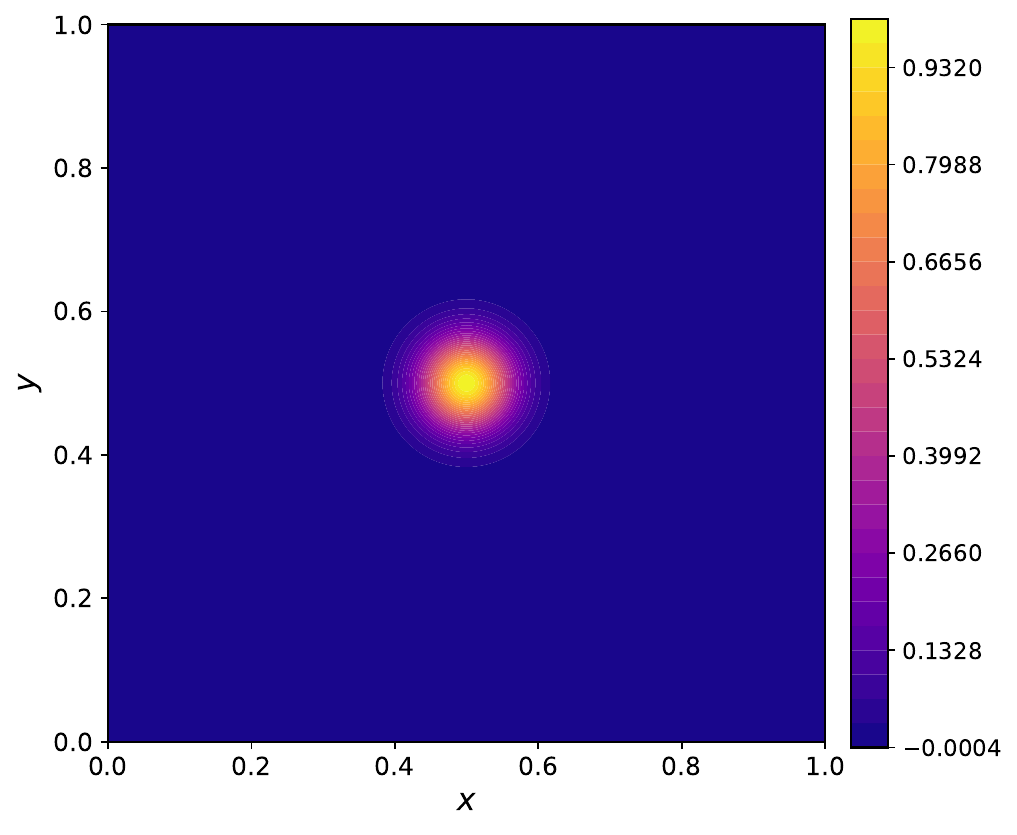}}
\subcaptionbox{Reference $u(\xi,\eta)$\label{fig:poisson_2d_ana_com}}[0.24\textwidth]{%
\includegraphics[width=\linewidth]{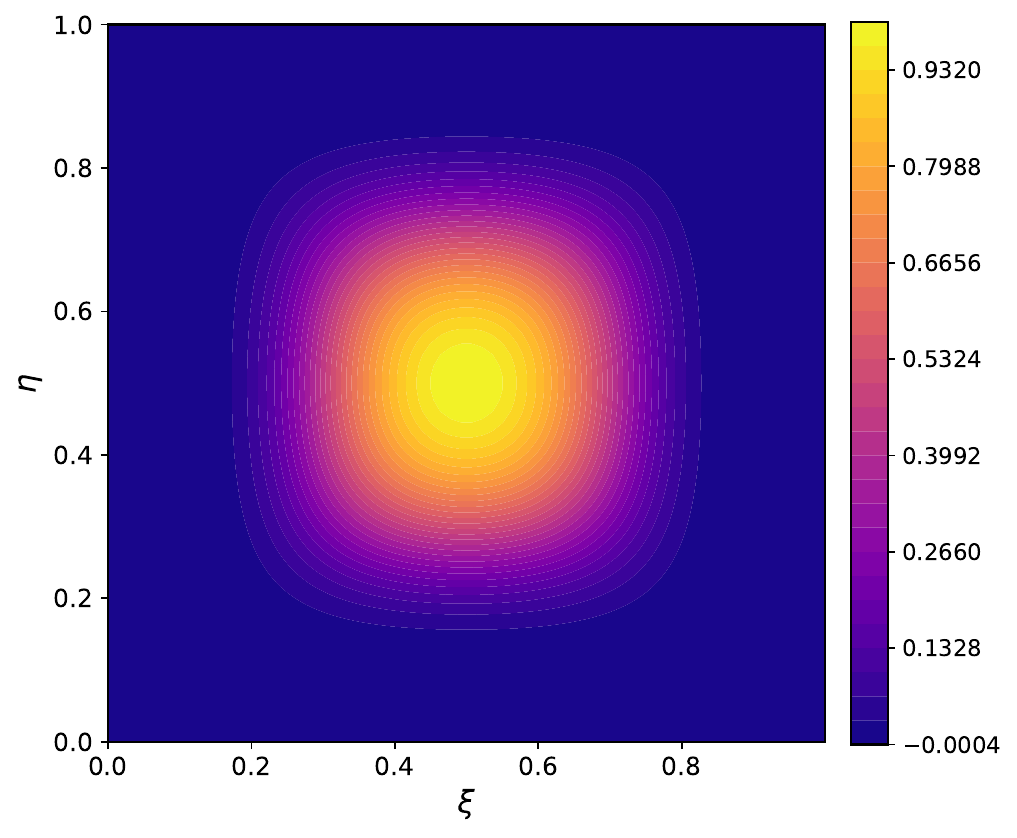}}
\subcaptionbox{Predicted $u(\xi,\eta)$ \label{fig:poisson_2d_pred_com}}[0.24\textwidth]{%
\includegraphics[width=\linewidth]{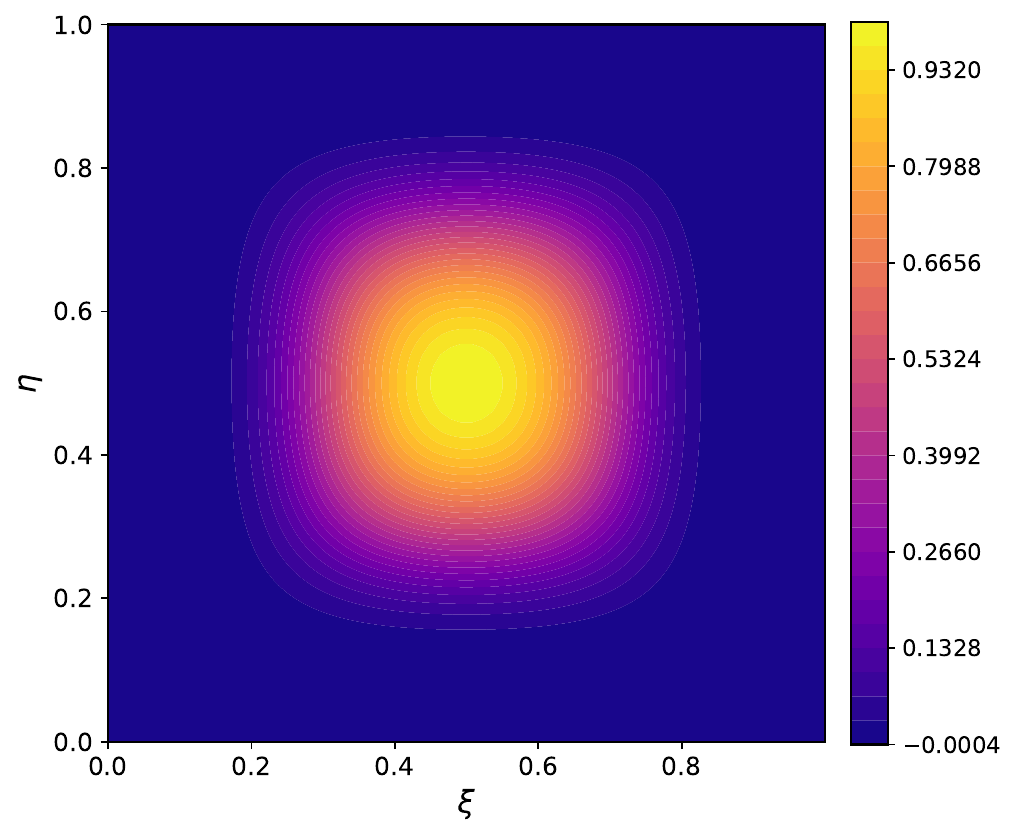}}\\
\subcaptionbox{Predicted $u_x(x,y)$ \label{fig:poisson_2d_ux}}[0.24\textwidth]{%
\includegraphics[width=\linewidth]{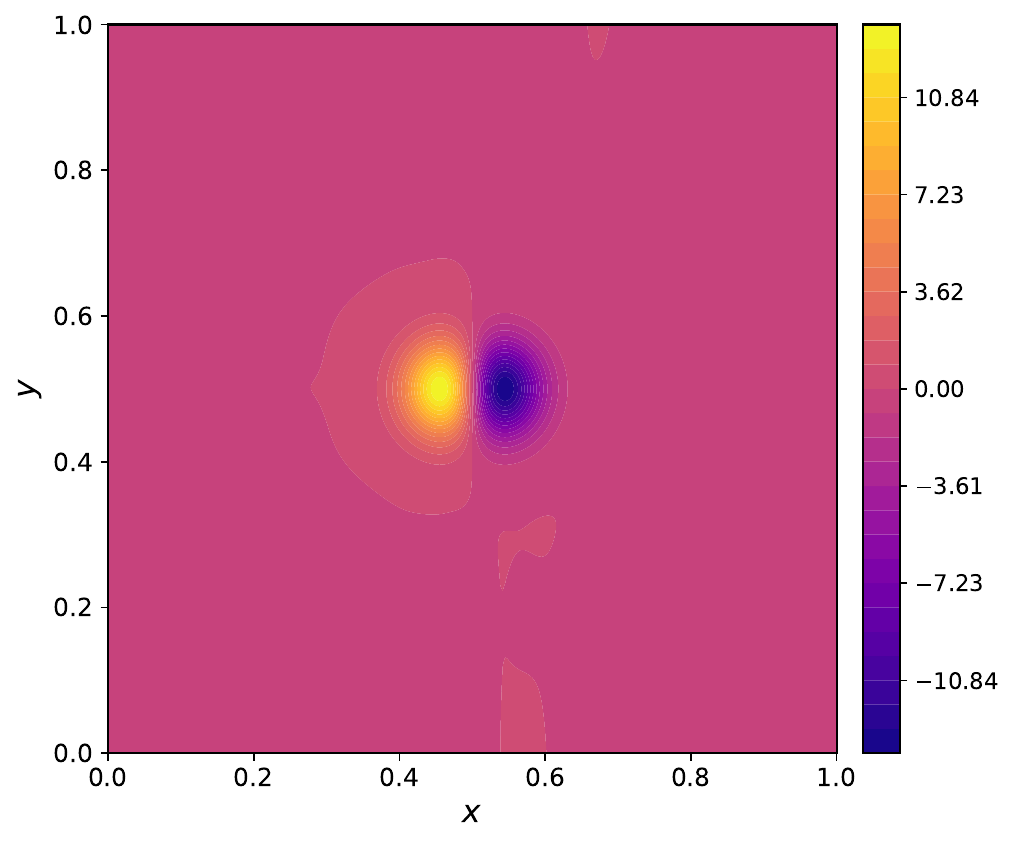}}
\subcaptionbox{Predicted $u_y(x,y)$ \label{fig:poisson_2d_uy}}[0.24\textwidth]{%
\includegraphics[width=\linewidth]{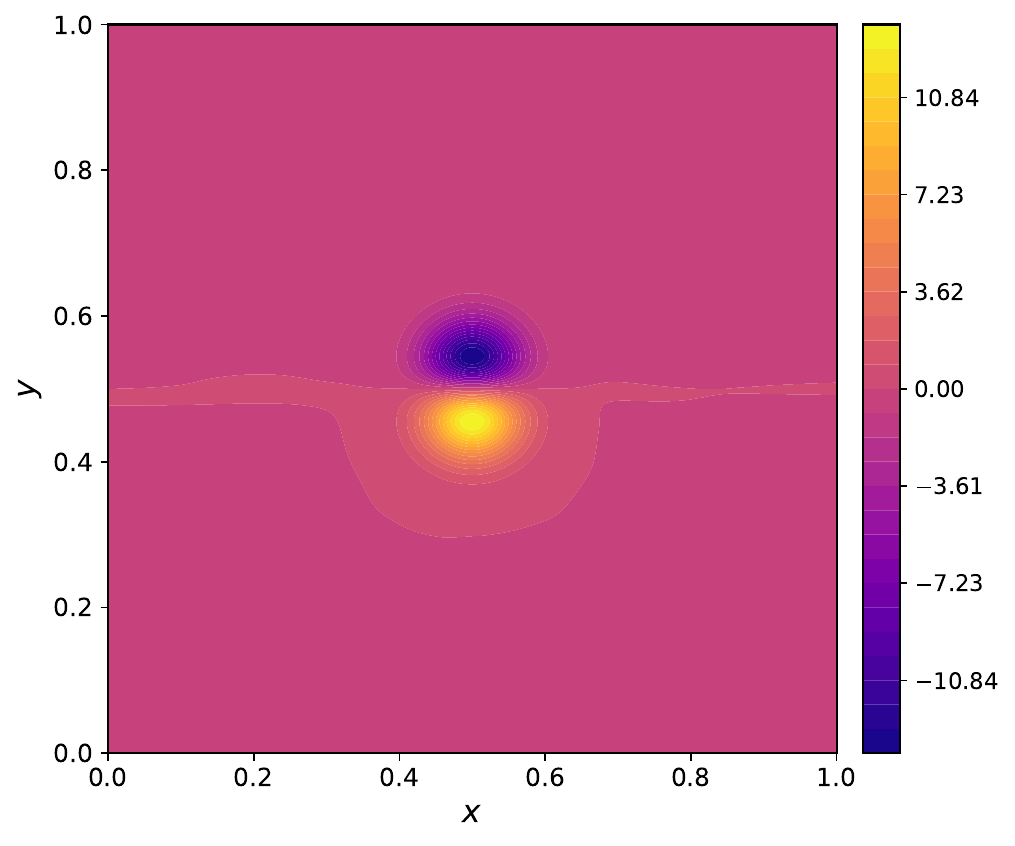}}
\subcaptionbox{Predicted $u_\xi(\xi,\eta)$\label{fig:poisson_2d_uxi}}[0.24\textwidth]{%
\includegraphics[width=\linewidth]{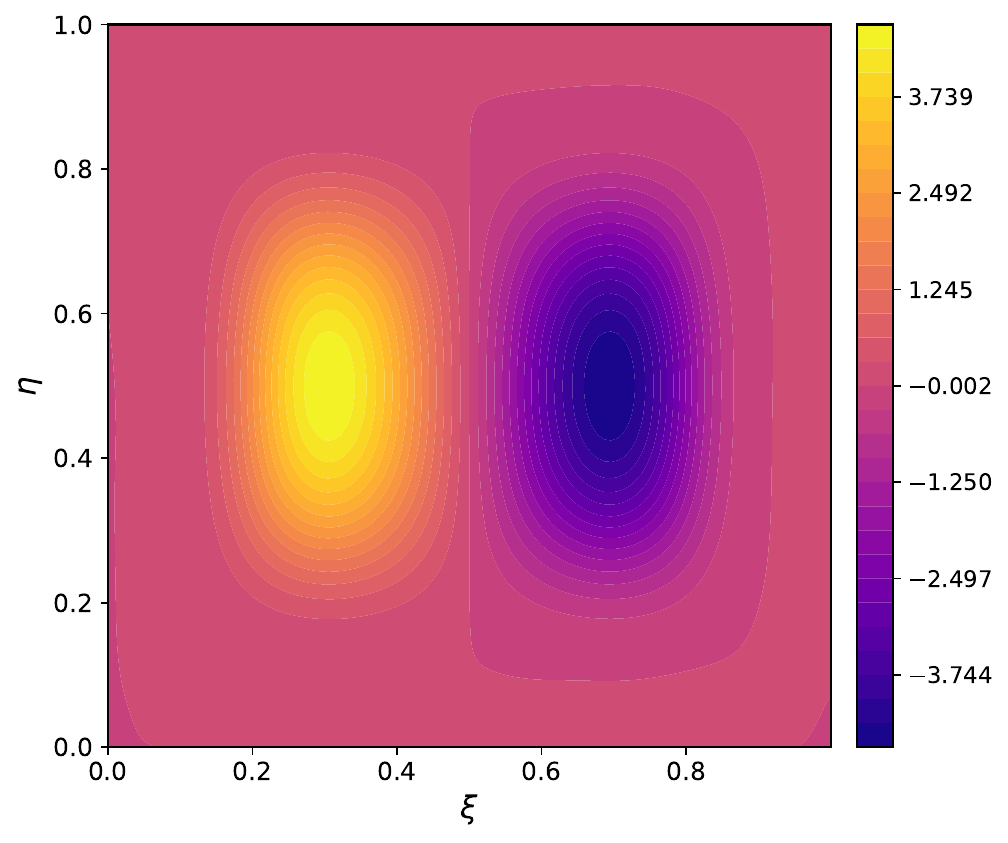}}
\subcaptionbox{Predicted $u_\eta(\xi,\eta)$ \label{fig:poisson_2d_ueta}}[0.24\textwidth]{%
\includegraphics[width=\linewidth]{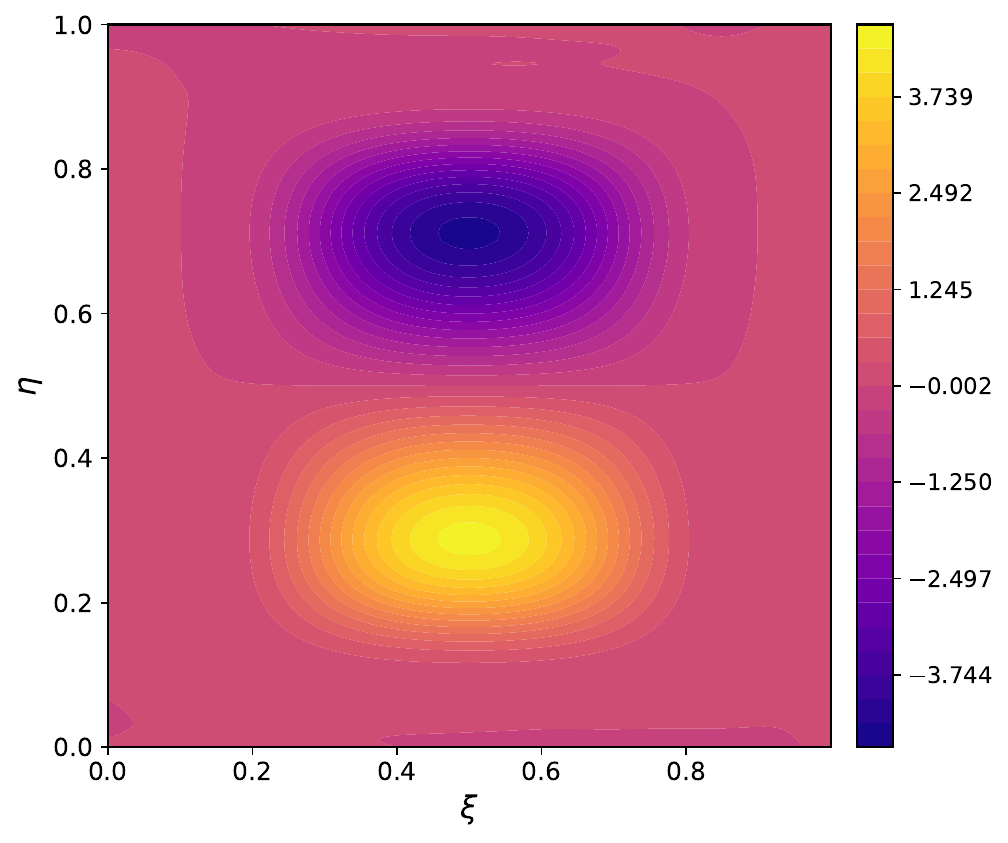}}
\caption{Solution of the 2D Poisson equation in \cref{eq: poisson} using CTPINNs with adaptive $T_{\psi_1}$. The variables $(x,y)$ and $(\xi,\eta)$ belong to the physical and computational domains, respectively. CTPINNs successfully solve the problem in the computational domain, and the transferred solution in the physical domain is also accurate. (c) and (d) show that the CT function stretches regions with large solution gradients and compresses regions with small gradients to promote gradient equidistribution. (e)--(h) show the gradients with respect to $x$, $y$, $\xi$, and $\eta$, respectively. The gradients in the computational coordinates are more uniform and have smaller maximum magnitudes, which is beneficial for training the network.}
\label{fig: poisson_2d}
\end{figure}

\cref{fig: poisson_2d_alpha} shows the relative error and the learnable parameters $\alpha_x$ and $\alpha_y$ during training as CTPINNs with $T_{\psi_1}$(A) solve the 2D Poisson equation. The dynamics of $\alpha_x$ and $\alpha_y$ show how the coordinate transformation adapts during training to promote gradient equidistribution. In particular, steeper gradients correspond to smaller values of $\alpha_x$ and $\alpha_y$, which in turn lead to coordinate transformations that stretch regions with large gradients. We also show the solution and its gradients in the physical and computational domains in \cref{fig: poisson_2d}. These results indicate that CTPINNs with adaptive $T_{\psi_1}$ can stretch regions with large gradients and compress regions with small gradients to promote gradient equidistribution. Notably, the fixed CT function $T_{\psi_1}$ can be combined with adaptive sampling methods, providing additional flexibility to improve the training efficiency and accuracy of PINNs for problems with steep gradients. The adaptive CT function $T_{\psi_1}$ may conflict with an adaptive sampling strategy by changing the solution distribution during training.

\subsection{Implicit Construction}\label{sec: implicit_construction_exp}
In this section, we focus on the implicit construction of the coordinate transformation. For many practical PDEs, some features of the solution are unknown, making explicit construction challenging or infeasible. This motivates learning the transformation implicitly, as described in \cref{section: implicit_construction}. The network $\boldsymbol{f}^\theta$ is a lightweight, two-layer MLP with $\tanh$ activation.
\subsubsection{Challenging PDEs}\label{sec: challenging_benchmark}
First, we conduct experiments on three PDEs from the JAX-PI benchmark~\cite{wang2022and}: the Allen-Cahn equation, the Korteweg-De Vries equation, and the lid-driven cavity flow problem. Their formulations and boundary conditions are as follows:

\textbf{Allen-Cahn.} We consider the following Allen-Cahn equation with periodic boundary conditions:
\begin{equation}
\begin{aligned}
\frac{\partial u}{\partial t} - 0.0001 \frac{\partial^2 u}{\partial x^2} = 0, & \quad t \in [0,1], x \in [-1,1], \\
u(0,x) = x^2 \cos(\pi x), & \quad x \in [-1,1].
\end{aligned}
\end{equation}
\textbf{Korteweg-De Vries.} We consider the following Korteweg-De Vries equation:
\begin{equation}
\begin{aligned}
\frac{\partial u}{\partial t} + u \frac{\partial u}{\partial x} + (0.022)^2 \frac{\partial^3 u}{\partial x^3} = 0, & \quad t \in [0,1], x \in [-1,1], \\
u(x,0) = \cos(\pi x), & \quad x \in [-1,1], \\
u(t,-1) = u(t,1), & \quad t \in [0,1].
\end{aligned}
\end{equation}
\textbf{Lid-driven cavity flow.} We consider the following Navier-Stokes equation for lid-driven cavity flow with Dirichlet boundary conditions:
\begin{equation}
\begin{aligned}
\boldsymbol{u} \cdot \nabla \boldsymbol{u} + \nabla p - \frac{1}{3200} \Delta \boldsymbol{u} &= 0, \quad (x,y) \in [0,1]^2, \\
\nabla \cdot \boldsymbol{u} &= 0, \quad (x,y) \in [0,1]^2, \\
\boldsymbol{u}(x,1) = \left(1 - \frac{\cosh(50(x-0.5))}{\cosh(25)}, \, 0\right) &,  \quad (x,y) \in \partial [0,1]^2, \\
\boldsymbol{u}(x,0) = \boldsymbol{u}(0,y) & = \boldsymbol{u}(1,y) = (0,0), \quad (x,y) \in \partial [0,1]^2.
\end{aligned}
\end{equation}

\begin{figure}[ht]
\centering
\subcaptionbox{Allen-Cahn\label{fig:ac_loss}}[0.32\textwidth]{%
\includegraphics[width=\linewidth]{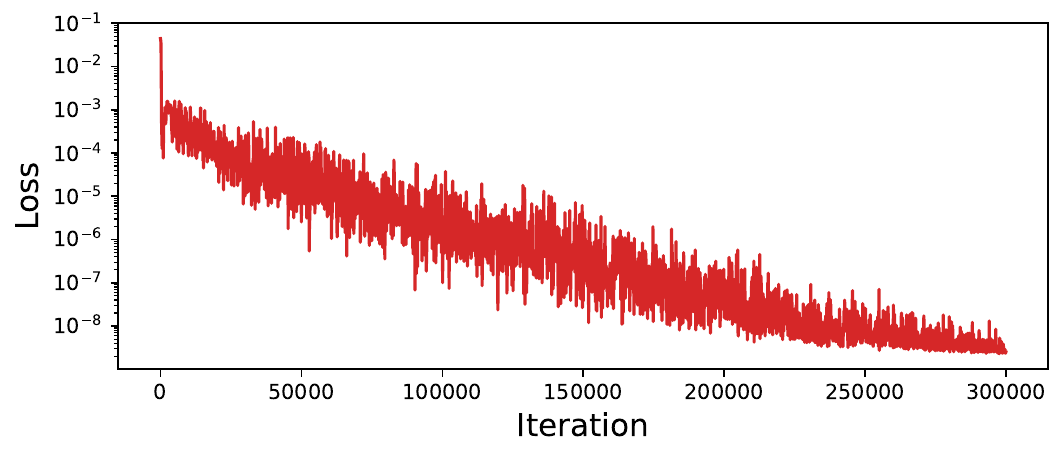}}
\hfill
\subcaptionbox{Korteweg-De Vries\label{fig:kdv_loss}}[0.32\textwidth]{%
\includegraphics[width=\linewidth]{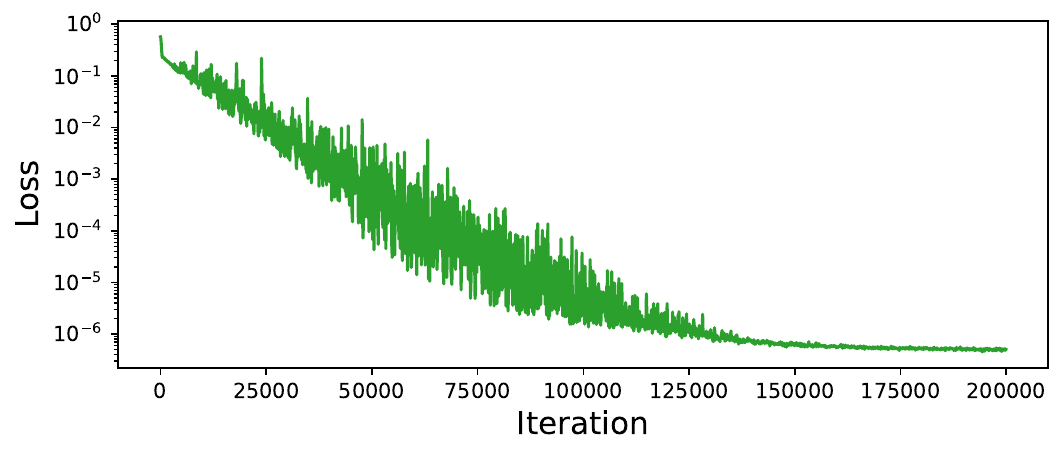}}
\hfill
\subcaptionbox{Lid-driven cavity\label{fig:ldc_loss}}[0.32\textwidth]{%
\includegraphics[width=\linewidth]{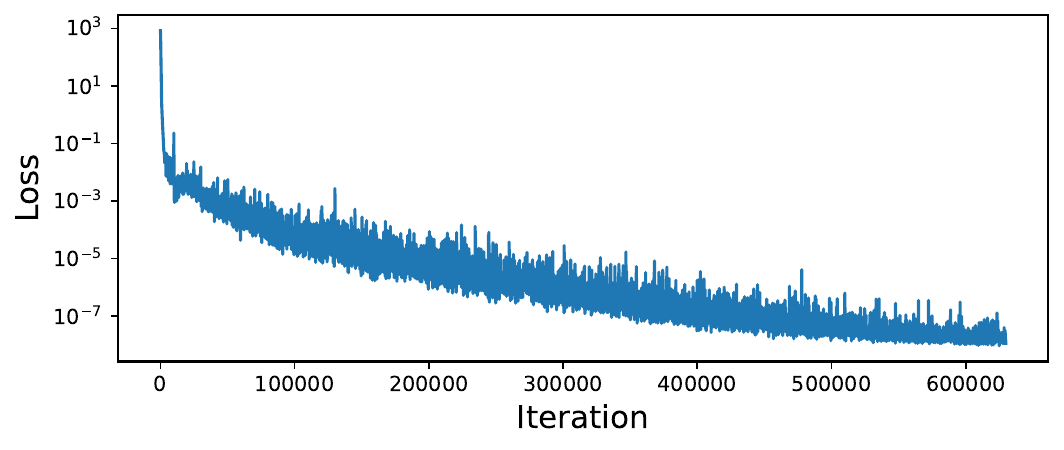}}\\
\subcaptionbox{Allen-Cahn\label{fig:ac_error}}[0.32\textwidth]{%
\includegraphics[width=\linewidth]{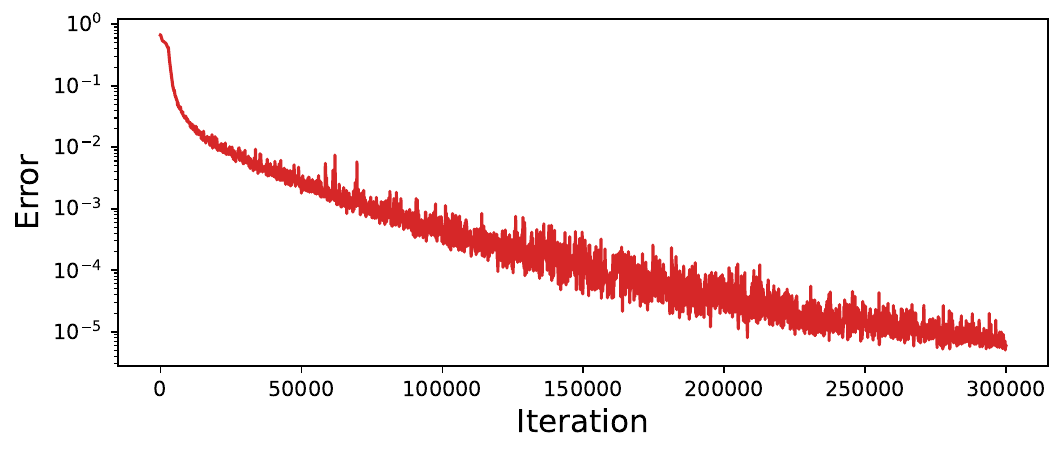}}
\hfill
\subcaptionbox{Korteweg-De Vries\label{fig:kdv_error}}[0.32\textwidth]{%
\includegraphics[width=\linewidth]{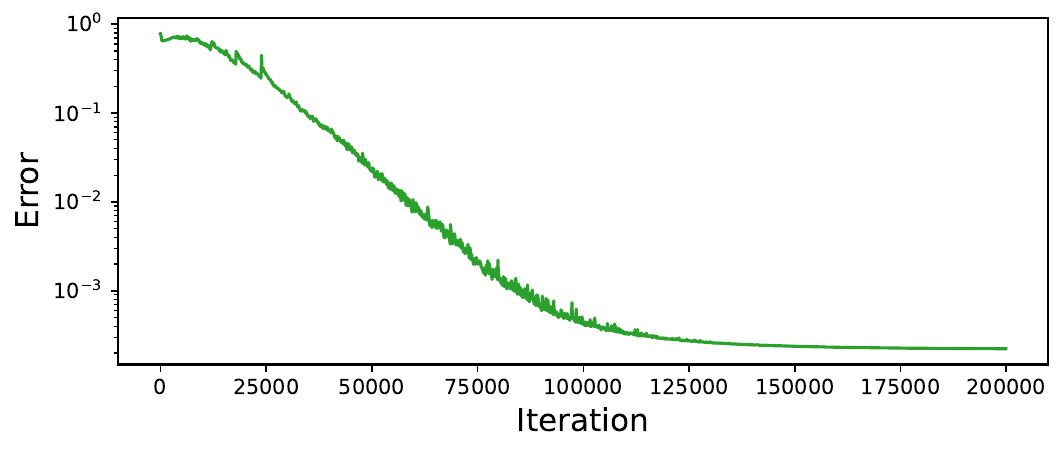}}
\hfill
\subcaptionbox{Lid-driven cavity\label{fig:ldc_error}}[0.32\textwidth]{%
\includegraphics[width=\linewidth]{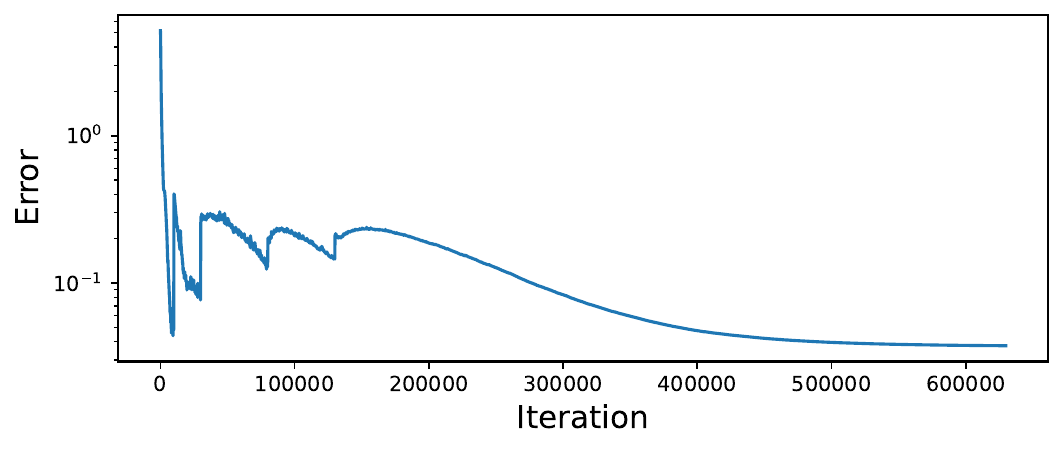}}

\caption{Training loss and relative error dynamics of JAX-PI benchmark for PirateNet+AT+$T_{\mathrm{exp}}^\theta$.}
\label{fig: loss_error_dynamics_jaxpi}
\end{figure}

We use PirateNet~\cite{wang2024piratenets} as the baseline and compare it with PirateNet combined with $T_{\mathrm{exp}}^\theta$. Notably, the results in~\cite{wang2024piratenets} were obtained using PirateNet with several advanced techniques, including random Fourier feature embeddings~\cite{tancik2020fourier}, NTK-based weighting schemes~\cite{wang2022and}, causality~\cite{wang2022respecting}, curriculum training~\cite{shukla2024comprehensive}, and random weight factorization~\cite{wang2022random}. Accordingly, in \cref{tab:result_jaxpi}, 'PirateNet+AT' denotes PirateNet with these advanced techniques. The 'PirateNet+AT+$T_{\mathrm{exp}}^\theta$' implementation also uses these techniques, showing that CTPINNs can be seamlessly integrated with state-of-the-art techniques rather than only with the basic PirateNet architecture. The results in \cref{tab:result_jaxpi} show that CTPINNs consistently improve the performance of PirateNet ('PirateNet+AT') on all three equations. The training dynamics of CTPINNs with PirateNet on these three equations are shown in \cref{fig: loss_error_dynamics_jaxpi}.
\begin{table}[ht]
\centering
\caption{Results on JAX-PI benchmarks}
\label{tab:result_jaxpi}
\begin{tabular}{lccc}
\toprule
 & Allen-Cahn & Korteweg-De Vries & Lid-driven cavity flow \\
\midrule
PirateNet+AT & $2.24 \times 10^{-5}$ & $4.27 \times 10^{-4}$ & $4.21 \times 10^{-2}$ \\
PirateNet+AT+$T_{\mathrm{exp}}^\theta$ & \cellcolor{blue!25}$6.01 \times 10^{-6}$ & \cellcolor{blue!25}$2.23 \times 10^{-4}$ & \cellcolor{blue!25}$3.75 \times 10^{-2}$ \\
Promotion & $73.2\%$ & $47.8\%$ & $10.9\%$ \\
\bottomrule
\end{tabular}
\end{table}

We further conduct experiments on three 3D PDEs from the SPINN benchmark~\cite{cho2023separable}: the Helmholtz equation, the Klein-Gordon equation, and the Navier-Stokes equation. Their formulations and boundary conditions are as follows:

\textbf{Helmholtz equation.} We consider the following steady Helmholtz equation:
\begin{equation}
\begin{aligned}
    \Delta u +  u &= q, && \boldsymbol{x}\in[-1,1]^3,\\
    u(\boldsymbol{x}) &= 0, && \boldsymbol{x}\in\partial[-1,1]^3.
\end{aligned}
\end{equation}
We use the analytical solution
$u(\boldsymbol{x})=\sin(4\pi x_1)\sin(4\pi x_2)\sin(3\pi x_3)$, with $q=-(4\pi)^2u-(4\pi)^2u-(3\pi)^2u+u$.

\textbf{Klein-Gordon equation.} We consider the following inhomogeneous Klein--Gordon equation:
\begin{equation}
\begin{aligned}
    \partial_{tt}u-\Delta u+u^2 &= f, && \boldsymbol{x}\in\Omega=[-1,1]^3,\ t\in[0,10],\\
    u(\boldsymbol{x},0) &= x_1+x_2+x_3, && \boldsymbol{x}\in\Omega=[-1,1]^3,\\
    u(\boldsymbol{x},t) &= g(\boldsymbol{x}), && \boldsymbol{x}\in\partial\Omega=[-1,1]^3,\ t\in[0,10].
\end{aligned}
\end{equation}
We use the analytical solution
$u(\boldsymbol{x},t)=(x_1+x_2+x_3)\cos(t)+x_1x_2x_3\sin(t)$; the source term $f$ and boundary data
$g$ are obtained from this exact solution.

\textbf{Navier-Stokes equation.} We use the vorticity formulation of the incompressible Navier--Stokes equation with periodic boundary conditions:
\begin{equation}
\begin{aligned}
    \partial_t\boldsymbol{\omega}
    +(\boldsymbol{u}\cdot\nabla)\boldsymbol{\omega}
    &= (\boldsymbol{\omega}\cdot\nabla)\boldsymbol{u}
       +0.05\Delta\boldsymbol{\omega}+\boldsymbol{F},
    && \boldsymbol{x}\in[0,2\pi]^3,\ t\in[0,5],\\
    \nabla\cdot\boldsymbol{u} &= 0,
    && \boldsymbol{x}\in[0,2\pi]^3,\ t\in[0,5],\\
    \boldsymbol{\omega}(\boldsymbol{x},0) &= \boldsymbol{\omega}_0(\boldsymbol{x}),
    && \boldsymbol{x}\in[0,2\pi]^3.
\end{aligned}
\end{equation}
The analytical velocity and vorticity fields are
\begin{equation}
\begin{aligned}
    u_x &= 2e^{-0.45 t}\cos(2x)\sin(2y)\sin(z),\\
    u_y &= -e^{-0.45 t}\sin(2x)\cos(2y)\sin(z),\\
    u_z &= -2e^{-0.45 t}\sin(2x)\sin(2y)\cos(z),\\
    \omega_x &= -3e^{-0.45 t}\sin(2x)\cos(2y)\cos(z),\\
    \omega_y &= 6e^{-0.45 t}\cos(2x)\sin(2y)\cos(z),\\
    \omega_z &= -6e^{-0.45 t}\cos(2x)\cos(2y)\sin(z).
\end{aligned}
\end{equation}
The forcing term $\boldsymbol{F}$ is given by
\begin{equation}
\begin{aligned}
    F_x &= -6e^{-0.9 t}\sin(4y)\sin(2z),\\
    F_y &= -6e^{-0.9 t}\sin(4x)\sin(2z),\\
    F_z &= 6e^{-0.9 t}\sin(4x)\sin(4y).
\end{aligned}
\end{equation}

\begin{figure}[ht]
\centering
\subcaptionbox{Reference x-y\label{fig:ns_ref_xy}}[0.32\textwidth]{%
\includegraphics[width=\linewidth]{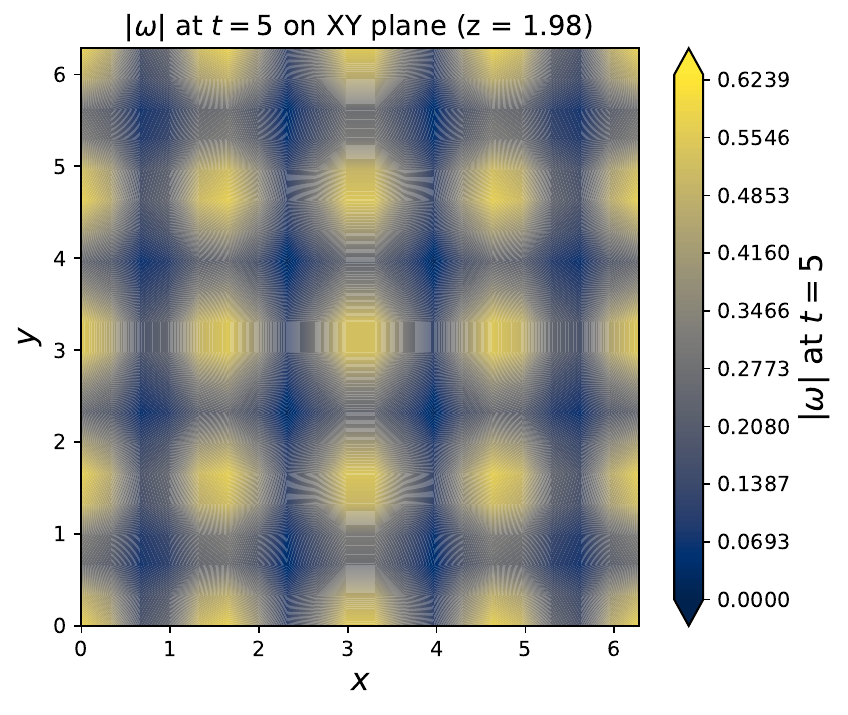}}
\hfill
\subcaptionbox{Reference x-z\label{fig:ns_ref_xz}}[0.32\textwidth]{%
\includegraphics[width=\linewidth]{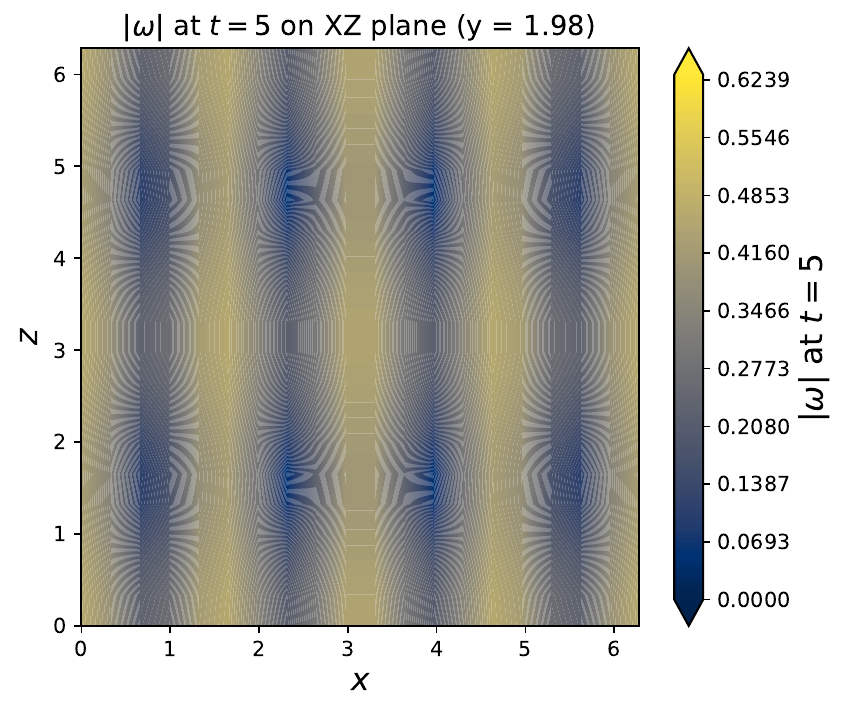}}
\hfill
\subcaptionbox{Reference y-z\label{fig:ns_ref_yz}}[0.32\textwidth]{%
\includegraphics[width=\linewidth]{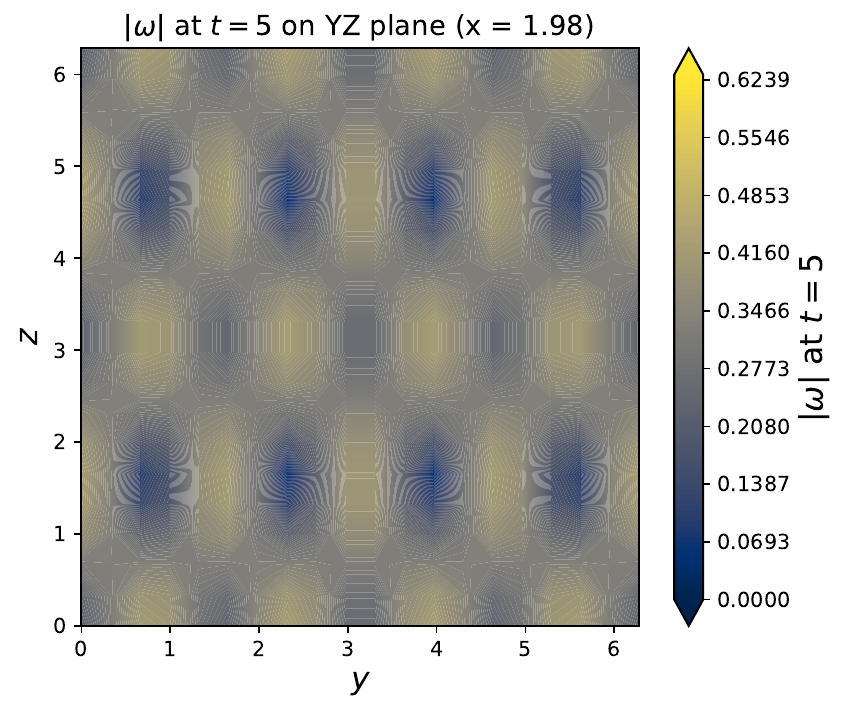}}\\
\subcaptionbox{Predicted x-y\label{fig:ns_pred_xy}}[0.32\textwidth]{%
\includegraphics[width=\linewidth]{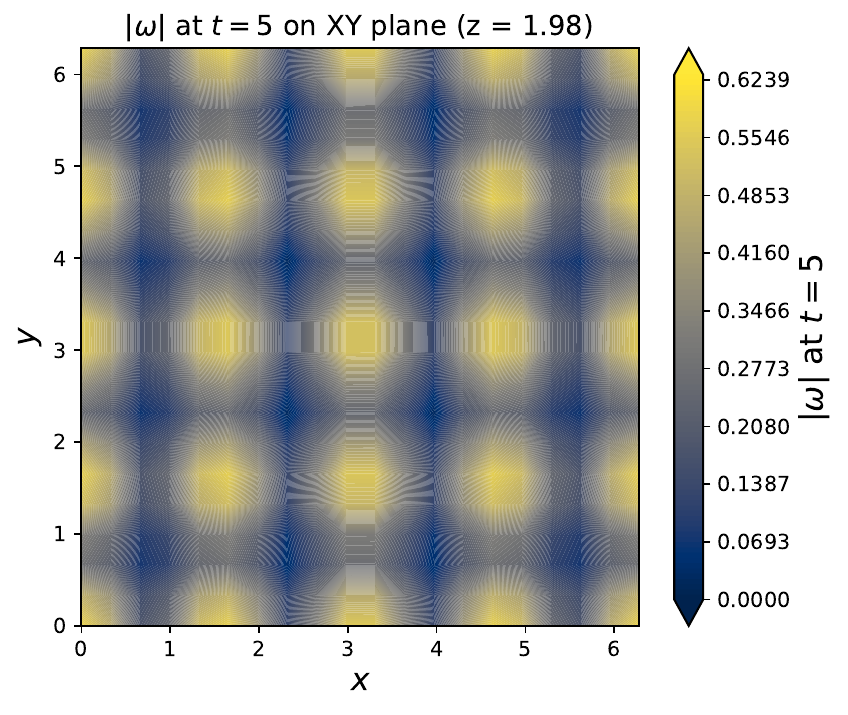}}
\hfill
\subcaptionbox{Predicted x-z\label{fig:ns_pred_xz}}[0.32\textwidth]{%
\includegraphics[width=\linewidth]{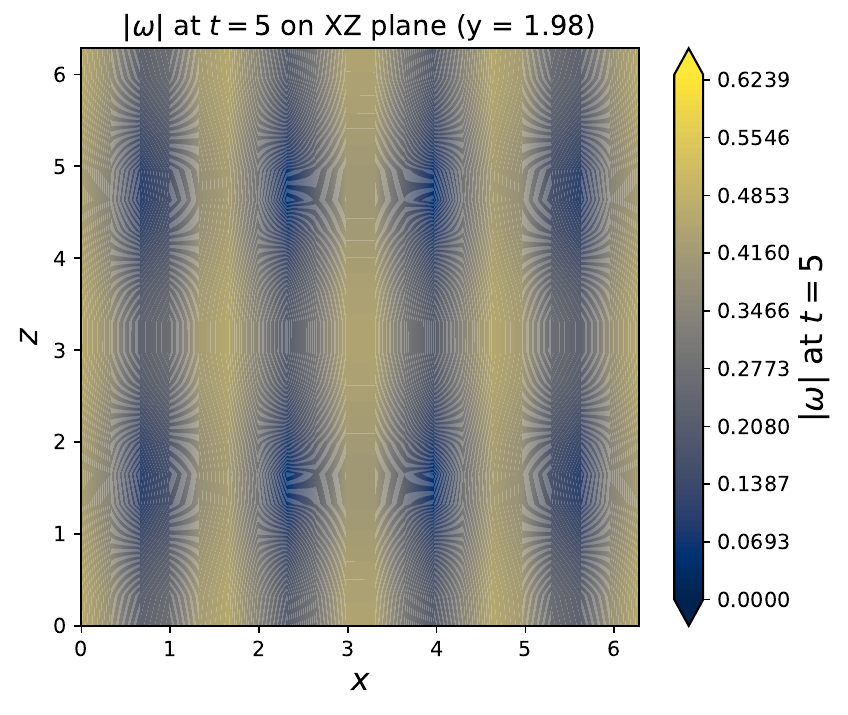}}
\hfill
\subcaptionbox{Predicted y-z\label{fig:ns_pred_yz}}[0.32\textwidth]{%
\includegraphics[width=\linewidth]{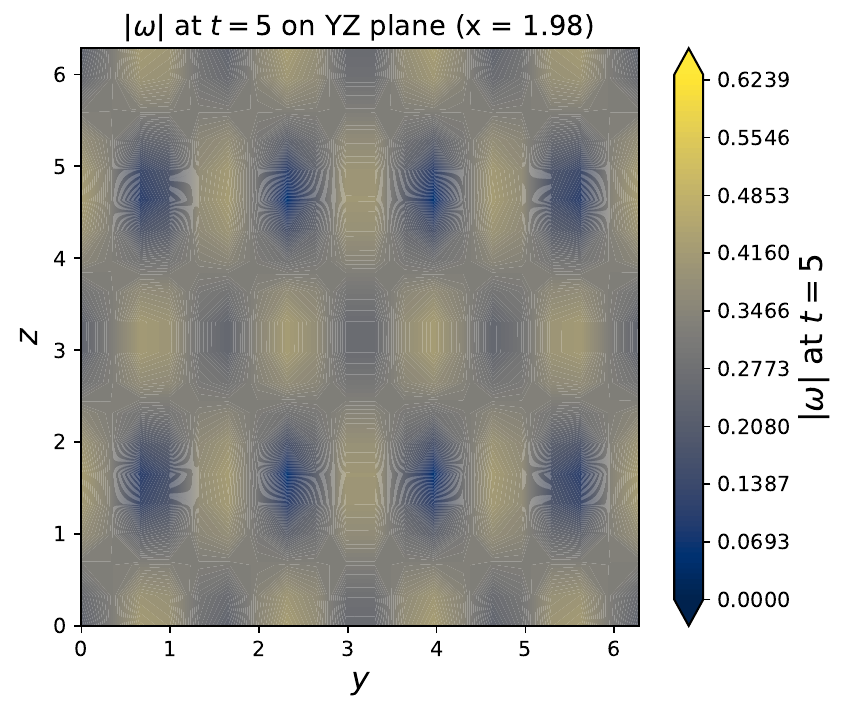}}

\caption{Reference and predicted solution distributions for Navier-Stokes equation in SPINN benchmarks}
\label{fig: solution_spinn}
\end{figure}

We use SPINN~\cite{cho2023separable} as the baseline method and compare it with SPINN+$T_{\mathrm{exp}}^\theta$. The results are summarized in \cref{tab:result_spinn}, and we also show the reference and predicted solution distributions of the Navier-Stokes equation in \cref{fig: solution_spinn}.

\begin{table}[ht]
\centering
\caption{Results on SPINN benchmarks}
\label{tab:result_spinn}
\begin{tabular}{lccc}
\toprule
 & Helmholtz & Klein-Gordon & Navier-Stokes \\
\midrule
SPINN & $3.22 \times 10^{-2}$ & $1.09 \times 10^{-3}$ & $3.33 \times 10^{-3}$ \\
SPINN+$T_{\mathrm{exp}}^\theta$ & \cellcolor{blue!25}$1.60 \times 10^{-2}$ & \cellcolor{blue!25}$8.51 \times 10^{-4}$ & \cellcolor{blue!25}$1.47 \times 10^{-3}$ \\
Promotion & $50.3\%$ & $21.9\%$ & $55.9\%$ \\
\bottomrule
\end{tabular}
\end{table}

\subsubsection{High-dimensional PDEs}\label{sec: high_dimensional_sg_equations} We next demonstrate the performance of the proposed CT module on high-dimensional PDEs. Specifically, we consider the following Sine-Gordon equation in high-dimensional spaces of varying dimension $d$~\cite{shi2024stochastic}.
\begin{equation}
\begin{aligned}
    \Delta u + \sin(u) & = f \quad x\in [0,1]^d,\\
    \quad u & = g \quad x \in \partial[0,1]^d.
\end{aligned}
\end{equation}
We define the reference solution as
\begin{equation}
u(\boldsymbol{x};c)=\left(\frac{1}{d-2}\sum_{i=1}^{d-2}\exp\left(-c x_i x_{i+1}x_{i+2}\right)\right)
\left(1-\frac{1}{d}\sum_{i=1}^{d}x_i^2\right),
\end{equation}
where $c$ is a prescribed parameter, $f=\Delta u+\sin(u)$, and $g=u|_{\partial[0,1]^d}$. Furthermore, to obtain a more uniform set of collocation points in high-dimensional space, we use the Halton low-discrepancy sequence from QRPINNs~\cite{yu2026quasi}.

\begin{figure}[ht]
\centering
\subcaptionbox{$d=3$\label{fig:sg_3}}[0.32\textwidth]{%
\includegraphics[width=\linewidth]{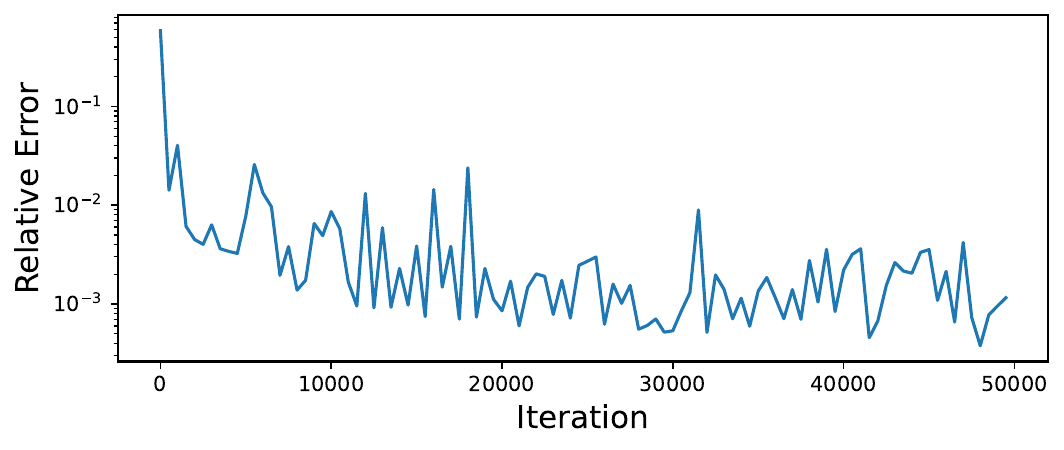}}
\hfill
\subcaptionbox{$d=100$\label{fig:sg_100}}[0.32\textwidth]{%
\includegraphics[width=\linewidth]{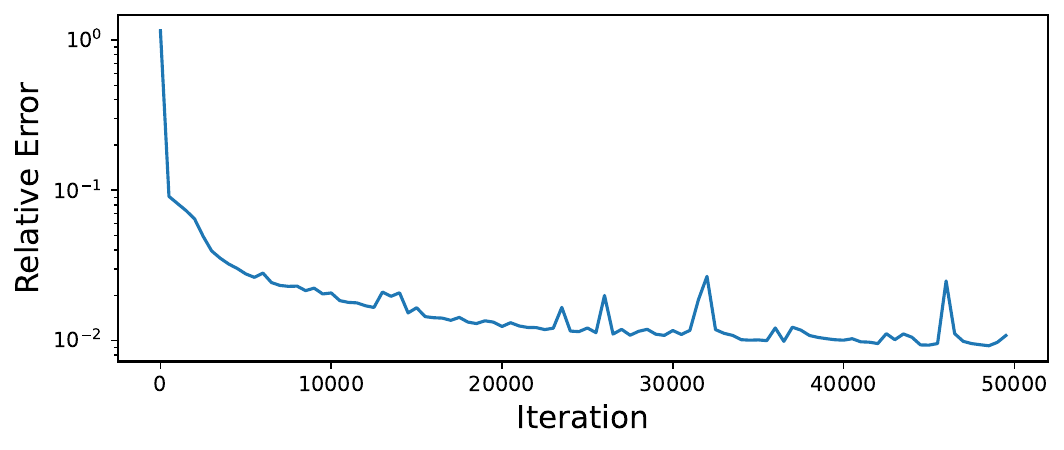}}
\hfill
\subcaptionbox{$d=200$\label{fig:sg_200}}[0.32\textwidth]{%
\includegraphics[width=\linewidth]{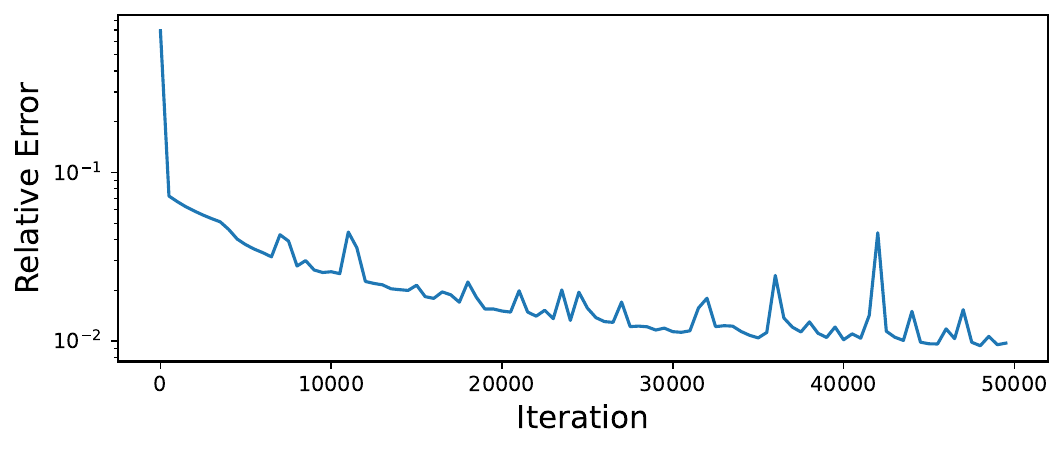}}

\caption{Error dynamics of QRPINN+$T_{\mathrm{exp}}^\theta$ for high-dimensional Sine-Gordon equations.}
\label{fig: error_sg_high}
\end{figure}

The results in \cref{tab: result_sg_high} indicate that CTPINNs consistently achieve better accuracy than QRPINNs, even on high-dimensional problems. We also report the cost analysis in \cref{tab: cost_analysis_sg_high} and show the error dynamics of QRPINN+$T_{\mathrm{exp}}^\theta$ in \cref{fig: error_sg_high}.
\begin{table}[ht]
\centering
\caption{Results on high-dimensional Sine-Gordon equations}
\label{tab: result_sg_high}
\begin{tabular}{lcccc}
\toprule
 & $d=3$ & $d=30$ & $d=100$ & $d=200$ \\
\midrule
QRPINN & $3.25 \times 10^{-4}$ & $4.07 \times 10^{-3}$ & $1.26 \times 10^{-2}$ & $1.52 \times 10^{-2}$ \\
QRPINN+$T_{\mathrm{exp}}^\theta$ & \cellcolor{blue!25}$2.91 \times 10^{-4}$ & \cellcolor{blue!25}$3.91 \times 10^{-3}$ & \cellcolor{blue!25}$9.05 \times 10^{-3}$ & \cellcolor{blue!25}$9.66 \times 10^{-3}$ \\
Promotion & $10.5\%$ & $3.9\%$ & $28.2\%$ & $36.4\%$ \\
\bottomrule
\end{tabular}
\end{table}

\subsection{Cost Analysis}\label{sec: cost_analysis}
In this section, we analyze the computational cost of CTPINNs. All experiments are conducted on a single NVIDIA A40 GPU using JAX. In the tables below, '$|\theta|$' denotes the number of trainable parameters, and 'Training Time (s)' denotes the total wall-clock training time in seconds.

The training cost of CTPINNs is an important consideration: they can improve performance but also incur additional costs. For explicit constructions, the additional training cost is negligible, especially with fixed majorants, since fixed majorants require only an extra forward pass to compute the transformation. For implicit constructions, CTPINNs incur a moderate additional training cost but provide clear accuracy gains, yielding a favorable trade-off in practice. The ablation results in \cref{tab: cost_analysis_bl1d} show that simply increasing the network size does not yield the same improvement as adding a CT module with a comparable number of parameters, underscoring the importance of the CT plug-in module.

Furthermore, \cref{tab: cost_analysis_sg_high} reports the cost for Sine-Gordon problems with different values of $d$. The results demonstrate the efficiency and effectiveness of CTPINNs and indicate that the additional time remains moderate and practically acceptable relative to the baseline across low- to high-dimensional problems.

\begin{table}[ht]
\centering
\caption{Cost analysis for 1D boundary layer problems.}
\label{tab: cost_analysis_bl1d}
\begin{adjustbox}{width=\textwidth}
\begin{tabular}{llcccccc}
\toprule
& & \multicolumn{6}{c}{MLP size (depth$\times$width)} \\
Method & Metric & $5\times 50$ & $7\times 50$ & $10\times 50$ & $5\times 100$ & $7\times 100$ & $10\times 100$ \\
\midrule
\multicolumn{1}{c}{\multirow{3}{*}{PINN}} & $|\theta|$ & $7801$ & $12901$ & $20551$ & $30601$ & $50801$ & $81101$ \\
 & Error & $1.58 \times 10^{1}$ & $1.60 \times 10^{-2}$ & $1.33 \times 10^{1}$ & $1.18 \times 10^{1}$ & $1.17 \times 10^{1}$ & $5.18 \times 10^{-2}$ \\
 & Training Time (s) & $2.98 \times 10^{2}$ & $7.36 \times 10^{2}$ & $1.16 \times 10^{3}$ & $3.33 \times 10^{2}$ & $8.75 \times 10^{2}$ & $1.52 \times 10^{3}$ \\
\midrule
\multicolumn{1}{c}{\multirow{3}{*}{PINN+$T_{\mathrm{exp}}^\theta$}} & $|\theta|$ & $7900$ & $12999$ & $20649$ & $30699$ & $50899$ & $81199$ \\
 & Error & $7.58 \times 10^{-3}$ & $5.81 \times 10^{-3}$ & $6.82 \times 10^{-3}$ & $1.65 \times 10^{-2}$ & $1.41 \times 10^{-2}$ & $6.86 \times 10^{-3}$ \\
 & Training Time (s) & $3.34 \times 10^{2}$ & $1.09 \times 10^{3}$ & $1.28 \times 10^{3}$ & $4.88 \times 10^{2}$ & $1.46 \times 10^{3}$ & $1.90 \times 10^{3}$ \\
\bottomrule
\end{tabular}
\end{adjustbox}
\end{table}

\begin{table}[ht]
\centering
\caption{Cost analysis for Sine-Gordon equations.}
\label{tab: cost_analysis_sg_high}
\begin{tabular}{lcccc}
\toprule
$d$ & \multicolumn{2}{c}{QRPINN} & \multicolumn{2}{c}{QRPINN+$T_{\mathrm{exp}}^\theta$} \\
& $|\theta|$ & Training Time (s) & $|\theta|$ & Training Time (s) \\
\midrule
$3$ & $7901$ & $1.38 \times 10^{2}$ & $8256$ & $2.31 \times 10^{2}$\\
$100$ & $12751$ & $4.24 \times 10^{3}$ & $22903$ & $5.62 \times 10^{3}$\\
$200$ & $17751$ & $1.18 \times 10^{4}$ & $38003$ & $1.57 \times 10^{4}$\\
\bottomrule
\end{tabular}
\end{table}

%% file: main/conclusion.tex
\section{Conclusion}
In this work, we introduced CTPINNs, a framework that improves PINNs by incorporating fixed or learnable coordinate transformations to reduce solution complexity through gradient equidistribution. For the canonical boundary-layer problem, we show that training CTPINNs with the original loss function drives the coordinate parameter toward values that make the solution gradient more uniformly distributed. Explicit constructions effectively resolve sharp gradients, empirically exhibit the tendency toward gradient equidistribution, and improve accuracy on Poisson problems even when an optimal CT is not used. For implicit constructions, integrating the CT module into PirateNet, SPINN, and QRPINN improves accuracy across challenging benchmark PDEs. Across a range of PDEs, CTPINNs demonstrate flexibility and effectiveness, with the coordinate transformation $T$ serving as a plug-in module that offers a favorable trade-off between accuracy and computational cost.

\paragraph{Limitations and future work} Our study suggests a useful module for PINNs: problems with prior knowledge can be effectively handled using explicitly constructed transformations, while more complex cases require flexible implicit constructions. However, neither construction is perfect: explicit transformations rely on prior knowledge of the problem, whereas implicit constructions lack proper normalization and theoretical guarantees. Developing more suitable CT functions is an important direction for future work. Furthermore, for PDEs defined on irregular domains, several studies~\cite{gao2021phygeonet,SahliCostabal2024DeltaPINNs} have used neural networks to map the irregular domain onto a regular one, thereby avoiding direct alterations to the solution distribution on the irregular domain. Building on this idea, a promising direction is to design a CT function that maps an irregular domain to a regular one while simultaneously transforming the solution gradients toward a more uniform distribution.

\section*{Acknowledgements}
The author T. Yu acknowledges support from the China Scholarship Council under Grant No. 202308090008. The authors acknowledge GPT for its assistance with literature reviews and drafting related-work sections; formulating theorems, algorithms, proofs, and other claims; and writing code and conducting computational experiments. The authors assume responsibility for all content.

%% file: main/appendices.tex
\section{Proof of the Gradient Equidistribution Theorem}\label[appendix]{appendix: proof_equidistribution_optimality}
This section provides the proof of \cref{thm:equidistribution_optimality}.
\begin{proof}[Proof of \cref{thm:equidistribution_optimality}]
Since $\xi=T(x)$ and $v(\xi)=u(T^{-1}(\xi))$, the chain rule gives
\begin{equation}
	v'(\xi)=\frac{u'(x)}{T'(x)},
	\qquad \xi=T(x).
\end{equation}
Therefore, by the change of variables $\xi=T(x)$,
\begin{equation}
	\int_0^1 |v'(\xi)|^2\,\mathrm{d}\xi
	= \int_0^1 \left|\frac{u'(x)}{T'(x)}\right|^2 T'(x)\,\mathrm{d}x
	= \int_0^1 \frac{|u'(x)|^2}{T'(x)}\,\mathrm{d}x.
\end{equation}
Applying the Cauchy--Schwarz inequality with
\begin{equation}
	f(x)=\frac{|u'(x)|}{\sqrt{T'(x)}},
	\qquad
	g(x)=\sqrt{T'(x)},
\end{equation}
we obtain
\begin{equation}
	\left(\int_0^1 |u'(x)|\,\mathrm{d}x\right)^2
	= \left(\int_0^1 f(x)g(x)\,\mathrm{d}x\right)^2
	\le \left(\int_0^1 \frac{|u'(x)|^2}{T'(x)}\,\mathrm{d}x\right)
	   \left(\int_0^1 T'(x)\,\mathrm{d}x\right).
\end{equation}
Because $T(1)-T(0)=1$, we have $\int_0^1 T'(x)\,\mathrm{d}x=1$, and thus
\begin{equation}\label{eq: functional_cauch_ge}
	\int_0^1 \frac{|u'(x)|^2}{T'(x)}\,\mathrm{d}x
	\ge \left(\int_0^1 |u'(x)|\,\mathrm{d}x\right)^2.
\end{equation}
Equality in Cauchy--Schwarz holds if and only if $f$ and $g$ are proportional, namely if and only if there exists a constant $c>0$ such that
\begin{equation}
	\frac{|u'(x)|}{\sqrt{T'(x)}} = c\sqrt{T'(x)}.
\end{equation}
Equivalently,
\begin{equation}
	T'(x)=\frac{|u'(x)|}{c}.
\end{equation}
Using again $\int_0^1 T'(x)\,\mathrm{d}x=1$, we find
\begin{equation}
	c=\int_0^1 |u'(s)|\,\mathrm{d}s,
\end{equation}
and hence
\begin{equation}
	T'(x)=\frac{|u'(x)|}{\int_0^1 |u'(s)|\,\mathrm{d}s}.
\end{equation}
Because $u'(x)\neq 0$ a.e. on $[0,1]$, this candidate satisfies $T'(x)>0$ a.e.; together with $T(0)=0$ and $T(1)=1$, it therefore belongs to $\mathcal A$.
For this choice,
\begin{equation}
	|v'(\xi)|=\left|\frac{u'(x)}{T'(x)}\right|
	= \int_0^1 |u'(s)|\,\mathrm{d}s,
\end{equation}
which is independent of $x$ and therefore of $\xi$.
\end{proof}

\section{Proofs of NTK Convergence Analysis}
\label[appendix]{appendix: non_optimal_ntk}
This section provides the proofs of \cref{prop:boundary_layer_ntk_regular}, \cref{prop:boundary_layer_ntk_gain}, \cref{lemma: unique_K_root}, and \cref{thm:boundary_layer_ntk_acceleration}.

\begin{proof}[Proof of \cref{prop:boundary_layer_ntk_regular}]
Let
\begin{equation}
	\Theta(\xi,\zeta)
	:=\lim_{n_1,\ldots,n_L\to\infty}
	\left\langle\nabla_\theta v^\theta(\xi),
	\nabla_\theta v^\theta(\zeta)\right\rangle
\end{equation}
denote the frozen infinite-width output NTK. Under the stated Gaussian initialization, the covariance and NTK recursions consist of Gaussian expectations of $\tanh$ and its derivatives. Since $\tanh\in C_b^\infty(\mathbb{R})$, differentiation under these expectations shows recursively that $\Theta$ is smooth on $[0,1]^2$; in particular, $\Theta\in C^5([0,1]^2)$.

Then, for $x\in[0,1]$, define
\begin{equation}
\begin{aligned}
	A_\alpha(x)&:=\frac{\epsilon}{\sqrt{N_r}}\bigl(T_\alpha'(x)\bigr)^2,\\
	B_\alpha(x)&:=\frac{1}{\sqrt{N_r}}
	\bigl(\epsilon T_\alpha''(x)+T_\alpha'(x)\bigr),\\
	\mathcal{R}_{\alpha,x}^{(\xi)}
	&:=A_\alpha(x)\partial_\xi^2+B_\alpha(x)\partial_\xi.
\end{aligned}
\end{equation}
Commuting parameter and input derivatives, the PDE--PDE, PDE--boundary, and boundary--boundary blocks of the residual NTK are, respectively,
\begin{align}
	(K_\alpha^{\mathrm{pp}})_{ij}
	&=
	\left.
	\mathcal{R}_{\alpha,x_i}^{(\xi)}
	\mathcal{R}_{\alpha,x_j}^{(\zeta)}
	\Theta(\xi,\zeta)
	\right|_{\substack{\xi=T_\alpha(x_i)\\
	\zeta=T_\alpha(x_j)}},\label{eq: residual_ntk_pp_block}\\
	(K_\alpha^{\mathrm{pb}})_{i,s}
	&=
	\left.
	\mathcal{R}_{\alpha,x_i}^{(\xi)}
	\Theta(\xi,s)
	\right|_{\xi=T_\alpha(x_i)},
	\qquad s\in\{0,1\},\label{eq: residual_ntk_pb_block}\\
	(K_\alpha^{\mathrm{bb}})_{s,r}
	&=\Theta(s,r),
	\qquad s,r\in\{0,1\}.\label{eq: residual_ntk_bb_block}
\end{align}
Here we used $T_\alpha(0)=0$ and $T_\alpha(1)=1$.

Since the map $(\alpha,x)\mapsto T_\alpha(x)$ is smooth on
$(0,\infty)\times[0,1]$, for every $\bar\alpha>\epsilon$ and
$\mathcal{I}=[\epsilon,\bar\alpha]$, the functions $T_\alpha$,
$T_\alpha'$, $T_\alpha''$ and their first $\alpha$-derivatives are uniformly
bounded on $\mathcal{I}\times [0,1]$. Differentiating
\cref{eq: residual_ntk_pp_block,eq: residual_ntk_pb_block,eq: residual_ntk_bb_block}
with respect to $\alpha$ involves derivatives of $\Theta$ of total order at
most five. Since the collocation set is finite and
$\Theta\in C^5([0,1]^2)$, it follows that
\begin{equation}
	C_{\mathcal{I}}:=\sup_{\gamma\in\mathcal{I}}
	\left\|\partial_\gamma K_\gamma\right\|_2<\infty.
\end{equation}
Therefore, for all $\alpha,\beta\in\mathcal{I}$, the fundamental theorem of calculus gives
\begin{equation}
	\|K_\alpha-K_\beta\|_2
	\le
	\int_{\min\{\alpha,\beta\}}^{\max\{\alpha,\beta\}}
	\|\partial_\gamma K_\gamma\|_2\,\mathrm{d}\gamma
	\le C_{\mathcal{I}}|\alpha-\beta|.
\end{equation}
\end{proof}

\begin{proof}[Proof of \cref{prop:boundary_layer_ntk_gain}]
Because $K_\epsilon$ is real and symmetric, the Gershgorin circle theorem gives
\begin{equation}
	\lambda_\epsilon
	\ge
	\min_{1\le i\le N_{\mathrm{ntk}}}
	\left((K_\epsilon)_{ii}
	-\sum_{\substack{1\le j\le N_{\mathrm{ntk}}\\j\ne i}}
	|(K_\epsilon)_{ij}|\right)
	=\underline{\kappa}_\epsilon.
\end{equation}
Under \cref{eq: tangent_feature_separation},
$\eta_\epsilon=\underline{\kappa}_\epsilon-\lambda_I>0$, and hence
\begin{equation}
	\lambda_\epsilon
	\ge\underline{\kappa}_\epsilon
	=\lambda_I+\eta_\epsilon
	>\lambda_I.
\end{equation}
\end{proof}

\begin{proof}[Proof of \cref{lemma: unique_K_root}]
	Let $\Theta(\xi,\zeta)=\langle\Phi(\xi),\Phi(\zeta)\rangle$ be a feature
representation of the positive semidefinite output NTK. Since each residual
operator is linear, the entries of $K_\alpha$ satisfy
\begin{equation}
\begin{aligned}
\left(K_\alpha\right)_{ij}
&=\left\langle\mathcal{R}_{\alpha,i}\Phi,
\mathcal{R}_{\alpha,j}\Phi\right\rangle,\\
z^\top K_\alpha z
&=\left\|\sum_{i=1}^{N_{\mathrm{ntk}}}
z_i\mathcal{R}_{\alpha,i}\Phi\right\|^2\ge0,
\qquad z\in\mathbb{R}^{N_{\mathrm{ntk}}}.
\end{aligned}
\end{equation}
Thus $K_\alpha$ is real symmetric and positive semidefinite. By the spectral
theorem, if $K_\alpha=Q\Lambda Q^\top$ with $\Lambda\ge0$, then
$K_\alpha^{1/2}:=Q\Lambda^{1/2}Q^\top$ is its unique symmetric positive
semidefinite square root.
\end{proof}

\begin{proof}[Proof of \cref{cor:boundary_layer_ntk_coherence}]
Using the Gram representation, for each $i$ we have
\begin{equation}
\begin{aligned}
\left(K_\epsilon\right)_{ii}
&=\|g_{\epsilon,i}\|_2^2\ge m_\epsilon^2,\\
\sum_{j\ne i}|(K_\epsilon)_{ij}|
&=\sum_{j\ne i}
|\langle g_{\epsilon,i},g_{\epsilon,j}\rangle|\\
&\le\mu_\epsilon\sum_{j\ne i}
\|g_{\epsilon,i}\|_2\|g_{\epsilon,j}\|_2\\
&\le(N_{\mathrm{ntk}}-1)\mu_\epsilon M_\epsilon^2.
\end{aligned}
\end{equation}
Thus \cref{eq: tangent_feature_coherence_condition} implies
\begin{equation}
\underline{\kappa}_\epsilon
\ge m_\epsilon^2-(N_{\mathrm{ntk}}-1)\mu_\epsilon M_\epsilon^2
>\lambda_I.
\end{equation}
The conclusion follows from \cref{prop:boundary_layer_ntk_gain}.
\end{proof}

\begin{proof}[Proof of \cref{thm:boundary_layer_ntk_acceleration}]
Fix $\alpha\in\mathcal{I}$ and regard $K_\alpha$ as a perturbation of
$K_\epsilon$:
\begin{equation}
	\Delta K_\alpha:=K_\alpha-K_\epsilon.
\end{equation}
Both matrices are real symmetric, so Weyl's inequality for their smallest
eigenvalues yields
\begin{equation}
	\lambda_{\min}(K_\alpha)
	=\lambda_{\min}(K_\epsilon+\Delta K_\alpha)
	\ge \lambda_{\min}(K_\epsilon)-\|\Delta K_\alpha\|_2.
\end{equation}
By \cref{prop:boundary_layer_ntk_regular}, taking $\beta=\epsilon$ gives the
operator-norm estimate
\begin{equation}
	\|\Delta K_\alpha\|_2
	=\|K_\alpha-K_\epsilon\|_2
	\le C_{\mathcal{I}}|\alpha-\epsilon|.
\end{equation}
Recalling that $\lambda_\alpha=\lambda_{\min}(K_\alpha)$, $\lambda_\epsilon=\lambda_{\min}(K_\epsilon)$ and using
$\alpha\ge\epsilon$, we therefore obtain
\begin{equation}\label{eq: ntk_weyl_alpha}
	\lambda_\alpha
	\ge \lambda_\epsilon-C_{\mathcal{I}}(\alpha-\epsilon).
\end{equation}
By Gershgorin's theorem and the definition of $\eta_\epsilon$, we have
\begin{equation}\label{eq: lamda_eps_inequality}
	\lambda_\epsilon
	\ge\underline{\kappa}_\epsilon
	=\lambda_I+\eta_\epsilon.
\end{equation}
Combining \cref{eq: lamda_eps_inequality} with \cref{eq: ntk_weyl_alpha} proves
\cref{eq: boundary_layer_ntk_rate}. The collocation-point separation condition
\cref{eq: tangent_feature_separation} implies $\eta_\epsilon>0$. Therefore,
for every $\alpha\in\mathcal{I}$ satisfying
$0\le\alpha-\epsilon<\eta_\epsilon/C_{\mathcal{I}}$, the rate bound in
\cref{eq: boundary_layer_ntk_rate} yields $\lambda_\alpha>\lambda_I$.

Finally, since an NTK Gram matrix is positive semidefinite, spectral decomposition of the fixed-coordinate dynamics gives
\begin{equation}
	\|\boldsymbol{r}_\alpha(t)\|_2^2
	=\|e^{-K_\alpha t}\boldsymbol{r}_\alpha(0)\|_2^2
	\le e^{-2\lambda_\alpha t}\|\boldsymbol{r}_\alpha(0)\|_2^2.
\end{equation}
Thus $\lambda_\alpha>\lambda_I$ gives a strictly faster worst-case exponential NTK rate than the identity-coordinate PINN.
\end{proof}

\section{Proofs of Implicit Tendency towards Gradient Equidistribution}
\label[appendix]{appendix: proof_coordinate_dynamics}
This section provides the proofs of \cref{thm:affine_alpha_convergence} and \cref{thm:alpha_descent}.
\begin{proof}[Proof of \cref{thm:affine_alpha_convergence}]
Because $T_\alpha$ is normalized so that $T_\alpha(0)=0$ and
$T_\alpha(1)=1$, the boundary loss depends only on $\theta$ and is independent
of $\alpha$. Define
\begin{equation}
	q(\alpha):=1-\frac{\epsilon}{\alpha},
	\qquad
	p_i(\alpha):=T_\alpha'(x_i),
\end{equation}
and
\begin{equation}
	A(\alpha):=\frac1{N_r}\sum_{i=1}^{N_r}p_i(\alpha)^2,
	\qquad
	B(\alpha):=\frac1{N_r}\sum_{i=1}^{N_r}
	p_i(\alpha)\partial_\alpha p_i(\alpha).
\end{equation}
For $v^{a,b}(\xi)=a+b\xi$, substitution into the residual gives
\begin{equation}
	R_i(a,b,\alpha)=bq(\alpha)p_i(\alpha).
\end{equation}
Consequently,
\begin{equation}
	\partial_\alpha\mathcal L(a,b,\alpha)
	=b^2q(\alpha)H(\alpha),
\end{equation}
where $H(\alpha):=q'(\alpha)A(\alpha)+q(\alpha)B(\alpha)$. Since $q(\epsilon)=0$,
\begin{equation}
	\left.\dot\alpha\right|_{\alpha=\epsilon}
	=-\partial_\alpha\mathcal L(a,b,\epsilon)=0,
\end{equation}
so $\{\alpha=\epsilon\}$ is invariant. If a trajectory with
$\alpha(0)>\epsilon$ reaches this set, let
\begin{equation}
	t_*:=\inf\{t\ge0:\alpha(t)=\epsilon\}.
\end{equation}
Uniqueness of the smooth gradient flow then gives
\begin{equation}
	\alpha(t)=\epsilon\quad\text{for all }t\ge t_*.
\end{equation}
Hence the trajectory cannot cross the invariant set, that is,
\begin{equation}
	\alpha(0)\ge\epsilon
	\quad\Longrightarrow\quad
	\alpha(t)\ge\epsilon\quad\text{for all }t\ge0.
\end{equation}
Moreover, \cref{assump:stable_training} implies that the boundary loss tends
to zero, and hence
\begin{equation}
	b(t)\to e^{-1/\epsilon}-1\ne0.
\end{equation}
The residual loss is
\begin{equation}
	\mathcal L_r(a,b,\alpha)
	=\frac{b^2}{2}q(\alpha)^2A(\alpha).
\end{equation}
Each $p_i$ is positive and continuous on $[\epsilon,\infty)$ and satisfies
$p_i(\alpha)\to1$ as $\alpha\to\infty$. Therefore
\begin{equation}
	\underline A
	:=\inf_{\alpha\in[\epsilon,\infty)}A(\alpha)>0.
\end{equation}
Since $b(t)\to b_*:=e^{-1/\epsilon}-1\ne0$, there exist $t_0\ge0$ and
$\underline b>0$ such that
\begin{equation}
	|b(t)|\ge\underline b
	\qquad\text{for all }t\ge t_0.
\end{equation}
Moreover, convergence of the full loss implies
$0\le\mathcal L_r(t)\le\mathcal L(t)\to0$. Hence, for $t\ge t_0$,
\begin{equation}
	0\le q(\alpha(t))^2 =\frac{2\mathcal L_r(t)}{b(t)^2A(\alpha(t))} \le\frac{2\mathcal L_r(t)}{\underline b^2\underline A}
	\longrightarrow0.
\end{equation}
Thus $q(\alpha(t))\to0$. Since $\alpha(t)\ge\epsilon$ and
$q(\alpha)=1-\epsilon/\alpha$,
\begin{equation}
	\alpha(t)=\frac{\epsilon}{1-q(\alpha(t))}
	\longrightarrow\epsilon.
\end{equation}
\end{proof}

\begin{proof}[Proof of \cref{thm:alpha_descent}]
Define
\begin{equation}
	d_\epsilon:=1-e^{-1/\epsilon}>0,
	\qquad
	b(t):=v^{\theta(t)}(1)-v^{\theta(t)}(0).
\end{equation}
By \cref{assump:stable_training}, the boundary loss tends to zero, so
\begin{equation}
	b(t)\to-d_\epsilon.
\end{equation}
Consequently, there exists $t_1\ge t_0$ such that
\begin{equation}
	|b(t)|\ge\frac{d_\epsilon}{2},
	\qquad t\ge t_1.
\end{equation}
We use $q(\alpha)$, $p_i(\alpha)$, and $A(\alpha)$ from the proof of \cref{thm:affine_alpha_convergence}.
For all $\alpha\ge\epsilon$, they satisfy
\begin{equation}
	0\le q(\alpha)<1,
	\qquad
	A(\alpha)\ge\underline A>0.
\end{equation}
Moreover, define 
\begin{equation}
	M_\epsilon:=\frac{1}{\epsilon(1-e^{-1/\epsilon})},
\end{equation}
then since $\alpha(1-e^{-1/\alpha})$ is increasing on $(0,\infty)$, 
\begin{equation}
	p_i(\alpha)
	\le\frac{1}{\alpha(1-e^{-1/\alpha})}
	\le M_\epsilon.
\end{equation}
For $\xi_i(t):=T_{\alpha(t)}(x_i)$, the affine decomposition in
\cref{eq: mlp_affine_decomposition} gives
\begin{equation}
	R(t,x_i)
	=b(t)q(\alpha(t))p_i(\alpha(t))+e_i(t),
\end{equation}
where $e_i(t)=\epsilon p_i(\alpha(t))^2 w_{\xi\xi}^{\theta(t)}(\xi_i(t))+q(\alpha(t))p_i(\alpha(t)) w_\xi^{\theta(t)}(\xi_i(t))$.
Since $\|w_\xi^{\theta(t)}\|_\infty$ and
$\|w_{\xi\xi}^{\theta(t)}\|_\infty$ are bounded by $\delta_3(\theta(t))$,
\begin{equation}\label{eq: multilayer_residual_perturbation_bound}
	|e_i(t)|
	\le C_\epsilon\delta_3(\theta(t))p_i(\alpha(t)),
\end{equation}
where $C_\epsilon=\epsilon M_\epsilon+1$.
Applying the reverse triangle inequality to the residual vector yields
\begin{equation}
\begin{aligned}
	\sqrt{2\mathcal L_r(\theta(t),\alpha(t))}
	&\ge\sqrt{A(\alpha(t))}\bigl(
	|b(t)|q(\alpha(t))\\
	&\qquad\qquad-C_\epsilon\delta_3(\theta(t))\bigr).
\end{aligned}
\end{equation}
A sufficient upper bound for the proportionality constant is
\begin{equation}\label{eq: coordinate_rho_threshold}
	C_{\rho}:=\frac{d_\epsilon}{4C_\epsilon}
	=\frac{(1-e^{-1/\epsilon})^2}{4(2-e^{-1/\epsilon})}>0.
\end{equation}
This conservative choice depends only on $\epsilon$; collocation-point
dependence enters the residual bound through $\underline A$.
Let $\rho$ denote the proportionality constant in the trajectory condition,
with $0<\rho\le C_{\rho}$, so that
\begin{equation}
	\delta_3(\theta(t))\le\rho q(\alpha(t)),
	\qquad t\ge t_0.
\end{equation}
For $t\ge t_1$, this condition and $\rho\le C_{\rho}$ imply
\begin{equation}\label{eq: multilayer_residual_coercivity}
\begin{aligned}
	\sqrt{2\mathcal L_r(\theta(t),\alpha(t))}
	&\ge\sqrt{A(\alpha(t))}
	\left(\frac{d_\epsilon}{2}-C_\epsilon\rho\right)q(\alpha(t))\\
	&\ge\frac{d_\epsilon\sqrt{\underline A}}{4}q(\alpha(t)).
\end{aligned}
\end{equation}
Since $0\le\mathcal L_r(t)\le\mathcal L(t)\to0$, we obtain
\begin{equation}
	0\le q(\alpha(t))
	\le\frac{4\sqrt{2\mathcal L_r(\theta(t),\alpha(t))}}
	{d_\epsilon\sqrt{\underline A}}
	\longrightarrow0.
\end{equation}
Finally,
\begin{equation}
	\alpha(t)=\frac{\epsilon}{1-q(\alpha(t))}
	\longrightarrow\epsilon.
\end{equation}
\end{proof}

\section{Proof of Error Analysis}\label[appendix]{appendix: proof_error_transfer}
This section provides the proof of \cref{thm:error_transfer}.
\begin{proof}[Proof of \cref{thm:error_transfer}]
For brevity, set $\Omega=\hat{\Omega}=[0,1]$ and write
\begin{equation}
\begin{aligned}
	T&:=T_\alpha,\\
	m&:=m_\alpha,\\
	M&:=M_\alpha,\\
	u&:=u^*,\\
	u^\theta&:=u^{\theta,\alpha},\\
	v&:=v^*.
\end{aligned}
\end{equation}
Define $w(\xi) := v(\xi)-v^\theta(\xi)$,
so that
\begin{equation}
	u(x)-u^\theta(x) = w(T(x)).
\end{equation}
We first estimate the $L^2$ error. By the change of variables $\xi=T(x)$,
\begin{equation}
	\|u-u^\theta\|_{L^2(\Omega)}^2
	= \int_{\Omega} |w(T(x))|^2\,\mathrm{d}x
	= \int_{\hat{\Omega}} |w(\xi)|^2\,\left|(T^{-1})'(\xi)\right|\,\mathrm{d}\xi.
\end{equation}
Since $T'(x)\ge m>0$ for all $x\in\Omega$, we have
\begin{equation}
	\left|(T^{-1})'(\xi)\right| \le m^{-1},
	\qquad \forall \xi\in\hat{\Omega},
\end{equation}
which yields
\begin{equation}
	\|u-u^\theta\|_{L^2(\Omega)}^2 \le m^{-1}\|v-v^\theta\|_{L^2(\hat{\Omega})}^2.
\end{equation}
Taking square roots gives
\begin{equation}
	\|u-u^\theta\|_{L^2(\Omega)} \le C_1\|v-v^\theta\|_{L^2(\hat{\Omega})},
\end{equation}
where $C_1:=m^{-1/2}$. Next, by the chain rule,
\begin{equation}
	(u-u^\theta)'(x)
	= T'(x)\,w'(T(x)).
\end{equation}
Therefore,
\begin{equation}
	\|(u-u^\theta)'\|_{L^2(\Omega)}^2
	=\int_{\hat{\Omega}}T'(T^{-1}(\xi))
	|w'(\xi)|^2\,\mathrm{d}\xi
	\le M\|w'\|_{L^2(\hat{\Omega})}^2.
\end{equation}
Taking square roots gives
\begin{equation}
	\|(u-u^\theta)'\|_{L^2(\Omega)}
	\le C_2\|(v-v^\theta)'\|_{L^2(\hat{\Omega})},
\end{equation}
where $C_2:=M^{1/2}$.
Combining the $L^2$ and $H^1$ seminorm bounds yields
\begin{equation}
	\|u-u^\theta\|_{H^1(\Omega)}
	\le \max\{C_1,C_2\}\|v-v^\theta\|_{H^1(\hat{\Omega})}.
\end{equation}
The elementary one-dimensional Sobolev estimate
\begin{equation}
	\|w\|_{L^\infty(\hat{\Omega})}^2
	\le 2\|w\|_{H^1(\hat{\Omega})}^2
\end{equation}
gives
\begin{equation}
\begin{aligned}
	\|u-u^\theta\|_{H^1(\Omega)}^2
	&=\|w\circ T\|_{L^2(\Omega)}^2
	+\|(w\circ T)'\|_{L^2(\Omega)}^2\\
	&\le 2\|w\|_{H^1(\hat{\Omega})}^2
	+M\|w'\|_{L^2(\hat{\Omega})}^2\\
	&\le(2+M)\|w\|_{H^1(\hat{\Omega})}^2.
\end{aligned}
\end{equation}
Taking the better of the two estimates proves the result with
\begin{equation}
	C_T:=\min\left\{\max\{C_1,C_2\},\sqrt{2+M}\right\}.
\end{equation}
\end{proof}

\section{Details of Experiments}\label[appendix]{appendix: experimental_details}

\subsection{Explicit Construction}

\subsubsection{1D Boundary Layer Problems}
We use a fully connected network with five hidden layers of width=50 and $\tanh$ activation. The network takes the transformed coordinate $T_{\psi}(x)$ as input. We minimize the loss with $\lambda_{r}=1,\lambda_{bc}=100$ using AdamW. The initial learning rate is $10^{-3}$ and is decayed exponentially by a factor of $0.95$ every $10^4$ steps; for the adaptive runs with $\epsilon\leq 10^{-4}$, the initial learning rate is $10^{-4}$. We perform $2\times10^5$ optimization steps, resampling 800 interior collocation points without replacement from a uniform grid of 1,000 points every 1,000 steps. 

\subsubsection{2D Boundary Layer Problems}
We use the same five-hidden-layer, width=50 $\tanh$ network, with $\lambda_{r}=1,\lambda_{bc}=100$ and AdamW training objective as for the 1D problem. The initial learning rate is $10^{-3}$ with exponential decay by a factor of $0.95$ every $10^4$ steps; for the adaptive runs with $\epsilon\leq 10^{-4}$, it is set to $10^{-4}$, and the transformation parameters $\alpha$ are frozen for the first $10^4$ steps. We perform $2\times10^5$ optimization steps and resample 10,000 interior collocation points without replacement from a $1{,}000\times1{,}000$ uniform grid every 1,000 steps. 

\subsubsection{Poisson's Equations}
We use a five-hidden-layer, width=50 $\tanh$ network with an adaptive coordinate transformation applied independently to the two input coordinates. The transformation parameters are trained jointly with the network. The loss is with $\lambda_{r}=1,\lambda_{bc}=100$. We use Adam with an initial learning rate of $10^{-3}$, exponentially decayed by a factor of $0.95$ every $10^4$ steps. Training consists of 10 rounds of 50,000 steps; at the start of each round, we resample 2,000 interior points uniformly in $[0,1]^2$. We use 1,000 randomly sampled boundary points on each of the four edges. 

The RAD-specific hyperparameters are $k_{\mathrm{RAD}}=1, c_{\mathrm{RAD}}=1, r_{\mathrm{RAD}}=50$.

The ACLE-specific hyperparameters are $r_{\mathrm{ACLE}}=4, N_{\mathrm{ref}}^{\max}=2000, \varepsilon_{\mathrm{NTK}}=10^{-8}$.

\subsection{Implicit Construction}

\subsubsection{Challenging PDEs}
The hyperparameters of $T_{\mathrm{exp}}^\theta$ are as follows:

For $T_{\mathrm{exp}}^\theta$ in PirateNet+$T_{\mathrm{exp}}^\theta$ on LDC, we use a two-layer network with hidden width 128 and $\tanh$ activation. The scaling parameter is set to $p=0.001$, and the Lipschitz bound $L$ is refreshed every 100 training steps. Other hyperparameters are same as the baseline models from the corresponding repository.

For $T_{\mathrm{exp}}^\theta$ in 'SPINN+$T_{\mathrm{exp}}^\theta$', we use a two-layer width=$64$ network with $\tanh$ activation. The scaling parameter is set to $p=0.1$, and the Lipschitz bound $L$ is refreshed every 5000 epochs. Other hyperparameters are same as the baseline models from the corresponding repository.
\subsubsection{High-dimensional PDEs}
For $T_{\mathrm{exp}}^\theta$, we use a two-layer, width=50 network with $\tanh$ activation. The scaling parameter is initialized at $p=1$ and is learnable. We update the Lipschitz bound $L$ after every 100 optimization steps.

We use a five-hidden-layer, width=50 network with $\tanh$ activation. The transformation parameters are trained jointly with the network. The loss is with $\lambda_{r}=1,\lambda_{bc}=100$. We use Adam with an initial learning rate of $10^{-3}$, exponentially decayed by a factor of $0.95$ every $10^4$ steps. Training consists of 100 rounds of 500 steps; at the start of each round, we resample 2,000 interior points uniformly in $[0,1]^d$. We use $2000\times d$ randomly sampled boundary points.
Other hyperparameters are same as the baseline models from the corresponding repository.
\subsection{Cost Analysis}
For \cref{tab: cost_analysis_bl1d}, we use different iteration budgets for the MLP depths: $1\times10^5$ iterations for depth 5, $2\times10^5$ for depth 7, and $2.5\times10^5$ for depth 10. Although the experiments use MLPs of different sizes, the $T_{\psi}$ has the same architecture in all cases, with depth=2 and width=32. Other hyperparameters are same as provided above.

For \cref{tab: cost_analysis_sg_high}, all hyperparameters are provided above.

%% file: references.bib
@article{lau2024pinnacle,
  title={{PINNACLE}: {PINN} adaptive collocation and experimental points selection},
  author={Lau, Gregory Kang Ruey and Hemachandra, Apivich and Ng, See-Kiong and Low, Bryan Kian Hsiang},
  journal={arXiv preprint arXiv:2404.07662},
  year={2024}
}

@article{tang2023adversarial,
  title={Adversarial adaptive sampling: Unify {PINN} and optimal transport for the approximation of {PDEs}},
  author={Tang, Kejun and Zhai, Jiayu and Wan, Xiaoliang and Yang, Chao},
  journal={arXiv preprint arXiv:2305.18702},
  year={2023}
}

@inproceedings{rahaman2019spectral,
  title={On the spectral bias of neural networks},
  author={Rahaman, Nasim and Baratin, Aristide and Arpit, Devansh and Draxler, Felix and Lin, Min and Hamprecht, Fred and Bengio, Yoshua and Courville, Aaron},
  booktitle={International conference on machine learning},
  pages={5301--5310},
  year={2019},
  organization={PMLR}
}

@article{lagaris1998artificial,
  title={Artificial neural networks for solving ordinary and partial differential equations},
  author={Lagaris, Isaac E and Likas, Aristidis and Fotiadis, Dimitrios I},
  journal={IEEE transactions on neural networks},
  volume={9},
  number={5},
  pages={987--1000},
  year={1998},
  publisher={IEEE}
}

@article{wu2023comprehensive,
  title={A comprehensive study of non-adaptive and residual-based adaptive sampling for physics-informed neural networks},
  author={Wu, Chenxi and Zhu, Min and Tan, Qinyang and Kartha, Yadhu and Lu, Lu},
  journal={Computer Methods in Applied Mechanics and Engineering},
  volume={403},
  pages={115671},
  year={2023},
  publisher={Elsevier}
}

@article{shi2024stochastic,
  title={Stochastic {Taylor} derivative estimator: Efficient amortization for arbitrary differential operators},
  author={Shi, Zekun and Hu, Zheyuan and Lin, Min and Kawaguchi, Kenji},
  journal={Advances in Neural Information Processing Systems},
  volume={37},
  pages={122316--122353},
  year={2024}
}

@article{yu2025spectral,
  title={Spectral informed neural networks},
  journal={Journal of Computational and Applied Mathematics},
  volume={477},
  pages={117178},
  year={2025},
  doi={https://doi.org/10.1016/j.cam.2025.117178},
  author={Tianchi Yu and Yiming Qi and Ivan Oseledets and Shiyi Chen}
}

@article{raissi2019physics,
  title={Physics-informed neural networks: A deep learning framework for solving forward and inverse problems involving nonlinear partial differential equations},
  author={Raissi, Maziar and Perdikaris, Paris and Karniadakis, George E},
  journal={Journal of Computational physics},
  volume={378},
  pages={686--707},
  year={2019},
  publisher={Elsevier}
}

@article{chiu2022can,
  title={{CAN-PINN}: A fast physics-informed neural network based on coupled-automatic--numerical differentiation method},
  author={Chiu, Pao-Hsiung and Wong, Jian Cheng and Ooi, Chinchun and Dao, My Ha and Ong, Yew-Soon},
  journal={Computer Methods in Applied Mechanics and Engineering},
  volume={395},
  pages={114909},
  year={2022},
  publisher={Elsevier}
}

@article{yu2024sinc,
title = {Sinc {Kolmogorov-Arnold} Network and Its Application for Solving {PDEs} with Singularities},
journal = {Neural Networks},
pages = {109033},
year = {2026},
issn = {0893-6080},
author = {Tianchi Yu and Jingwei Qiu and Jiang Yang and Ivan Oseledets},
}

@article{gordon1973construction,
  title={Construction of curvilinear co-ordinate systems and applications to mesh generation},
  author={Gordon, William J and Hall, Charles A},
  journal={International Journal for Numerical Methods in Engineering},
  volume={7},
  number={4},
  pages={461--477},
  year={1973},
  publisher={Wiley Online Library}
}

@inproceedings{fanaskov2023general,
  title={General covariance data augmentation for neural {PDE} solvers},
  author={Fanaskov, Vladimir and Yu, Tianchi and Rudikov, Alexander and Oseledets, Ivan},
  booktitle={International Conference on Machine Learning},
  pages={9665--9688},
  year={2023},
  organization={PMLR}
}

@article{yu2026quasi,
  title={Quasi-random physics-informed neural networks},
  author={Yu, Tianchi and Oseledets, Ivan},
  journal={Neurocomputing},
  pages={132913},
  year={2026},
  publisher={Elsevier}
}

@article{wang2024piratenets,
  title={{Piratenets}: Physics-informed deep learning with residual adaptive networks},
  author={Wang, Sifan and Li, Bowen and Chen, Yuhan and Perdikaris, Paris},
  journal={Journal of Machine Learning Research},
  volume={25},
  number={402},
  pages={1--51},
  year={2024}
}

@article{wang2022and,
  title={When and why {PINNs} fail to train: A neural tangent kernel perspective},
  author={Wang, Sifan and Yu, Xinling and Perdikaris, Paris},
  journal={Journal of Computational Physics},
  volume={449},
  pages={110768},
  year={2022},
  publisher={Elsevier}
}

@article{yang2024moving,
  title={Moving sampling physics-informed neural networks induced by moving mesh {PDE}},
  author={Yang, Yu and Yang, Qihong and Deng, Yangtao and He, Qiaolin},
  journal={Neural Networks},
  volume={180},
  pages={106706},
  year={2024},
  publisher={Elsevier}
}

@article{hou2023enhancing,
  title={Enhancing {PINNs} for solving {PDEs} via adaptive collocation point movement and adaptive loss weighting},
  author={Hou, Jie and Li, Ying and Ying, Shihui},
  journal={Nonlinear Dynamics},
  volume={111},
  number={16},
  pages={15233--15261},
  year={2023},
  publisher={Springer}
}

@article{shen2026matched,
  title={Matched Asymptotic Expansions-Based Transferable Neural Networks for Singular Perturbation Problems},
  author={Shen, Zhequan and Ju, Lili and Zhu, Liyong},
  journal={International Journal for Numerical Methods in Engineering},
  volume={127},
  number={6},
  pages={e70307},
  year={2026},
  publisher={Wiley Online Library}
}

@article{gie2024semi,
  title={Semi-analytic {PINN} methods for boundary layer problems in a rectangular domain},
  author={Gie, Gung-Min and Hong, Youngjoon and Jung, Chang-Yeol and Munkhjin, Tselmuun},
  journal={Journal of Computational and Applied Mathematics},
  volume={450},
  pages={115989},
  year={2024},
  publisher={Elsevier}
}

@article{xu2026physics,
  title={Physics-Informed Neural Networks for Power Systems: A Systematic Review},
  author={Xu, Xianfeng and Qiang, Kang and Wu, Junyang and Xuan, Liang and Lu, Yong and Ma, Wenfei and Zhang, Jiale},
  journal={Sustainable Energy, Grids and Networks},
  pages={102261},
  year={2026},
  publisher={Elsevier}
}

@article{mehmood2026physics,
  title={A Physics Informed Neural Network ({PINN}) framework for fractional order modeling of {Alzheimer's} disease},
  author={Mehmood, Adnan and Farman, Muhammad and Afzal, Farkhanda and Nisar, Kottakkaran Sooppy and Ahmed, Mohammed Altaf and Hafez, Mohamed},
  journal={Frontiers in Neuroinformatics},
  volume={20},
  pages={1748481},
  year={2026}
}

@article{di2026physics,
  title={Physics-informed neural networks in water and wastewater systems: a critical review},
  author={Di Bella, Antonino and Raissi, Maziar and Santoro, Domenico and Roccaro, Paolo},
  journal={Water Research},
  pages={125449},
  year={2026},
  publisher={Elsevier}
}

@book{liseikin2017grid,
  title={Grid Generation Methods},
  author={Liseikin, Vladimir D},
  year={2017},
  publisher={Springer}
}

@article{li2001moving,
  title={Moving mesh methods in multiple dimensions based on harmonic maps},
  author={Li, Ruo and Tang, Tao and Zhang, Pingwen},
  journal={Journal of Computational Physics},
  volume={170},
  number={2},
  pages={562--588},
  year={2001},
  publisher={Elsevier}
}

@book{mali2013accuracy,
  title={Accuracy verification methods: Theory and algorithms},
  author={Mali, Olli and Neittaanm{\"a}ki, Pekka and Repin, Sergey},
  volume={32},
  year={2013},
  publisher={Springer Science \& Business Media}
}

@article{tancik2020fourier,
  title={Fourier features let networks learn high frequency functions in low dimensional domains},
  author={Tancik, Matthew and Srinivasan, Pratul and Mildenhall, Ben and Fridovich-Keil, Sara and Raghavan, Nithin and Singhal, Utkarsh and Ramamoorthi, Ravi and Barron, Jonathan and Ng, Ren},
  journal={Advances in neural information processing systems},
  volume={33},
  pages={7537--7547},
  year={2020}
}

@article{wang2022respecting,
  title={Respecting causality is all you need for training physics-informed neural networks},
  author={Wang, Sifan and Sankaran, Shyam and Perdikaris, Paris},
  journal={arXiv preprint arXiv:2203.07404},
  year={2022}
}

@article{wang2022random,
  title={Random weight factorization improves the training of continuous neural representations},
  author={Wang, Sifan and Wang, Hanwen and Seidman, Jacob H and Perdikaris, Paris},
  journal={arXiv preprint arXiv:2210.01274},
  year={2022}
}

@inproceedings{de2006good,
  title={Good approximation by splines with variable knots. II},
  author={De Boor, Carl},
  booktitle={Conference on the Numerical Solution of Differential Equations: Dundee 1973},
  pages={12--20},
  year={2006},
  organization={Springer}
}

@article{dorfi1987simple,
  title={Simple adaptive grids for 1-D initial value problems},
  author={Dorfi, EA and LO'C, Drury},
  journal={Journal of Computational Physics},
  volume={69},
  number={1},
  pages={175--195},
  year={1987},
  publisher={Elsevier}
}

@article{huang1994moving,
  title={Moving mesh partial differential equations ({MMPDES}) based on the equidistribution principle},
  author={Huang, Weizhang and Ren, Yuhe and Russell, Robert D},
  journal={SIAM Journal on Numerical Analysis},
  volume={31},
  number={3},
  pages={709--730},
  year={1994},
  publisher={SIAM}
}

@article{lei2023s,
  title={What's the situation with intelligent mesh generation: A survey and perspectives},
  author={Lei, Na and Li, Zezeng and Xu, Zebin and Li, Ying and Gu, Xianfeng},
  journal={IEEE transactions on visualization and computer graphics},
  volume={30},
  number={8},
  pages={4997--5017},
  year={2023},
  publisher={IEEE}
}

@article{peng20253dmeshnet,
  title={{3DMeshNet}: A three-dimensional differential neural network for structured mesh generation},
  author={Peng, Jiaming and Chen, Xinhai and Liu, Jie},
  journal={Graphical Models},
  volume={139},
  pages={101257},
  year={2025},
  publisher={Elsevier}
}

@article{zhang2024towards,
  title={Towards universal mesh movement networks},
  author={Zhang, Mingrui and Wang, Chunyang and Kramer, Stephan and Wallwork, Joseph G and Li, Siyi and Liu, Jiancheng and Chen, Xiang and Piggott, Matthew D},
  journal={Advances in Neural Information Processing Systems},
  volume={37},
  pages={14934--14961},
  year={2024}
}

@article{yu2024flow2mesh,
  title={{Flow2Mesh}: A flow-guided data-driven mesh adaptation framework},
  author={Yu, Jian and Lyu, Hongqiang and Xu, Ran and Ouyang, Wenxuan and Liu, Xuejun},
  journal={Physics of Fluids},
  volume={36},
  number={3},
  year={2024},
  publisher={AIP Publishing}
}

@inproceedings{rowbottom2025g,
  title={{G-Adaptivity}: optimised graph-based mesh relocation for finite element methods},
  author={Rowbottom, James and Maierhofer, Georg and Deveney, Teo and M{\"u}ller, Eike Hermann and Paganini, Alberto and Schratz, Katharina and Lio, Pietro and Sch{\"o}nlieb, Carola-Bibiane and Budd, Chris},
  booktitle={International Conference on Machine Learning},
  pages={52173--52197},
  year={2025},
  organization={PMLR}
}

@article{wang2025ugm2n,
  title={{UGM2N}: An Unsupervised and Generalizable Mesh Movement Network via {M-Uniform} Loss},
  author={Wang, Zhichao and Chen, Xinhai and Wang, Qinglin and Gao, Xiang and Zhang, Qingyang and Jia, Menghan and Zhang, Xiang and Liu, Jie},
  journal={arXiv preprint arXiv:2508.08615},
  year={2025}
}

@article{goodfellow2020generative,
  title={Generative adversarial networks},
  author={Goodfellow, Ian and Pouget-Abadie, Jean and Mirza, Mehdi and Xu, Bing and Warde-Farley, David and Ozair, Sherjil and Courville, Aaron and Bengio, Yoshua},
  journal={Communications of the ACM},
  volume={63},
  number={11},
  pages={139--144},
  year={2020},
  publisher={ACM New York, NY, USA}
}

@inproceedings{arjovsky2017wasserstein,
  title={Wasserstein generative adversarial networks},
  author={Arjovsky, Martin and Chintala, Soumith and Bottou, L{\'e}on},
  booktitle={International conference on machine learning},
  pages={214--223},
  year={2017},
  organization={Pmlr}
}

@article{gulrajani2017improved,
  title={Improved training of {Wasserstein} {GANs}},
  author={Gulrajani, Ishaan and Ahmed, Faruk and Arjovsky, Martin and Dumoulin, Vincent and Courville, Aaron C},
  journal={Advances in neural information processing systems},
  volume={30},
  year={2017}
}

@article{miyato2018spectral,
  title={Spectral normalization for generative adversarial networks},
  author={Miyato, Takeru and Kataoka, Toshiki and Koyama, Masanori and Yoshida, Yuichi},
  journal={arXiv preprint arXiv:1802.05957},
  year={2018}
}

@inproceedings{rezende2015variational,
  title={Variational inference with normalizing flows},
  author={Rezende, Danilo and Mohamed, Shakir},
  booktitle={International conference on machine learning},
  pages={1530--1538},
  year={2015},
  organization={PMLR}
}

@inproceedings{behrmann2019invertible,
  title={Invertible residual networks},
  author={Behrmann, Jens and Grathwohl, Will and Chen, Ricky TQ and Duvenaud, David and Jacobsen, J{\"o}rn-Henrik},
  booktitle={International conference on machine learning},
  pages={573--582},
  year={2019},
  organization={PMLR}
}

@article{chen2018neural,
  title={Neural ordinary differential equations},
  author={Chen, Ricky TQ and Rubanova, Yulia and Bettencourt, Jesse and Duvenaud, David K},
  journal={Advances in neural information processing systems},
  volume={31},
  year={2018}
}

@article{wassing2025physics,
  title={Physics-informed neural networks for inviscid transonic flows around an airfoil},
  author={Wassing, Simon and Langer, Stefan and Bekemeyer, Philipp},
  journal={Physics of Fluids},
  volume={37},
  number={8},
  year={2025},
  publisher={AIP Publishing}
}

@article{schmid2026physics,
  title={Physics-informed neural operators for predicting structural intensity from laser Doppler vibrometry measurements of plates},
  author={Schmid, Johannes D and Zettel, Sebastian F and Marburg, Steffen},
  journal={Mechanical Systems and Signal Processing},
  volume={248},
  pages={114013},
  year={2026},
  publisher={Elsevier}
}

@article{jin2021nsfnets,
  title={{NSFnets} ({Navier-Stokes} flow nets): Physics-informed neural networks for the incompressible {Navier-Stokes} equations},
  author={Jin, Xiaowei and Cai, Shengze and Li, Hui and Karniadakis, George Em},
  journal={Journal of Computational Physics},
  volume={426},
  pages={109951},
  year={2021},
  publisher={Elsevier}
}

@article{shukla2024comprehensive,
  title={A comprehensive and {FAIR} comparison between {MLP} and {KAN} representations for differential equations and operator networks},
  author={Shukla, Khemraj and Toscano, Juan Diego and Wang, Zhicheng and Zou, Zongren and Karniadakis, George Em},
  journal={Computer Methods in Applied Mechanics and Engineering},
  volume={431},
  pages={117290},
  year={2024},
  publisher={Elsevier}
}

@article{wang2026gradient,
  title={Gradient alignment in physics-informed neural networks: A second-order optimization perspective},
  author={Wang, Sifan and Li, Bowen and Perdikaris, Paris and others},
  journal={Advances in Neural Information Processing Systems},
  volume={38},
  pages={168482--168532},
  year={2026}
}

@article{yu2022gradient,
  title={Gradient-enhanced physics-informed neural networks for forward and inverse {PDE} problems},
  author={Yu, Jeremy and Lu, Lu and Meng, Xuhui and Karniadakis, George Em},
  journal={Computer Methods in Applied Mechanics and Engineering},
  volume={393},
  pages={114823},
  year={2022},
  publisher={Elsevier}
}

@article{rathore2024challenges,
  title={Challenges in training pinns: A loss landscape perspective},
  author={Rathore, Pratik and Lei, Weimu and Frangella, Zachary and Lu, Lu and Udell, Madeleine},
  journal={arXiv preprint arXiv:2402.01868},
  year={2024}
}

@article{gao2023failure,
  title={Failure-informed adaptive sampling for {PINNs}},
  author={Gao, Zhiwei and Yan, Liang and Zhou, Tao},
  journal={SIAM Journal on Scientific Computing},
  volume={45},
  number={4},
  pages={A1971--A1994},
  year={2023},
  publisher={SIAM}
}

@article{lai2025hard,
  title={The hard-constraint {PINNs} for interface optimal control problems},
  author={Lai, Ming-Chih and Song, Yongcun and Yuan, Xiaoming and Yue, Hangrui and Zeng, Tianyou},
  journal={SIAM Journal on Scientific Computing},
  volume={47},
  number={3},
  pages={C601--C629},
  year={2025},
  publisher={SIAM}
}

@article{berger1984adaptive,
  title={Adaptive mesh refinement for hyperbolic partial differential equations},
  author={Berger, Marsha J and Oliger, Joseph},
  journal={Journal of computational Physics},
  volume={53},
  number={3},
  pages={484--512},
  year={1984},
  publisher={Elsevier}
}

@article{cho2023separable,
  title={Separable physics-informed neural networks},
  author={Cho, Junwoo and Nam, Seungtae and Yang, Hyunmo and Yun, Seok-Bae and Hong, Youngjoon and Park, Eunbyung},
  journal={Advances in Neural Information Processing Systems},
  volume={36},
  pages={23761--23788},
  year={2023}
}

@article{gao2021phygeonet,
  title   = {{PhyGeoNet}: Physics-informed geometry-adaptive
             convolutional neural networks for solving
             parameterized steady-state {PDEs} on irregular domain},
  author  = {Gao, Han and Sun, Luning and Wang, Jian-Xun},
  journal = {Journal of Computational Physics},
  volume  = {428},
  pages   = {110079},
  year    = {2021},
  doi     = {10.1016/j.jcp.2020.110079}
}

@article{SahliCostabal2024DeltaPINNs,
title = {$\Delta$-{PINNs}: Physics-informed neural networks on complex geometries},
journal = {Engineering Applications of Artificial Intelligence},
volume = {127},
pages = {107324},
year = {2024},
issn = {0952-1976},
doi = {https://doi.org/10.1016/j.engappai.2023.107324},
author = {Francisco {Sahli Costabal} and Simone Pezzuto and Paris Perdikaris}
}
